\documentclass[reqno]{amsart}  %LaTeX2e
\theoremstyle{plain}
\newtheorem{theorem}{Theorem}[section]

\newtheorem{remark}[theorem]{Remark}
\newtheorem{proposition}[theorem]{Proposition}

\newtheorem{corollary}[theorem]{Corollary}

\numberwithin{equation}{section}
\allowdisplaybreaks[1]

\theoremstyle{definition}

\newtheorem{definition}[theorem]{Definition}

\theoremstyle{remark}

\newcommand{\B}{{\mathbb B}}
\newcommand{\C}{{\mathbb C}}

\newcommand{\K}{{\mathbb K}}

\newcommand{\bU}{{\mathbf U}}

\newcommand{\bn}{{\mathbf n}}

\newcommand{\bS}{{\mathbf S}}
\newcommand{\cA}{{\mathcal A}}
\newcommand{\cB}{{\mathcal B}}
\newcommand{\cC}{{\mathcal C}}
\newcommand{\cD}{{\mathcal D}}
\newcommand{\cE}{{\mathcal E}}
\newcommand{\cF}{{\mathcal F}}
\newcommand{\cG}{{\mathcal {G}}}
\newcommand{\cH}{{\mathcal H}}
\newcommand{\cK}{{\mathcal K}}
\newcommand{\cL}{{\mathcal L}}
\newcommand{\cM}{{\mathcal M}}

\newcommand{\cS}{{\mathcal S}}
\newcommand{\cT}{{\mathcal T}}
\newcommand{\cU}{{\mathcal U}}
\newcommand{\cX}{{\mathcal X}}
\newcommand{\cY}{{\mathcal Y}}

\newcommand{\cZ}{{\mathcal Z}}
\newcommand{\D}{{\mathbb D}}

\newcommand{\cal}{\mathcal}

\newcommand{\sbm}[1]{\left[\begin{smallmatrix} #1
     \end{smallmatrix}\right]}

\newcommand{\inn}[2]{\ensuremath{\langle #1,\, #2 \rangle}}

\begin{document}

\begin{abstract}
The operator-valued Schur-class is defined to be the set of
holomorphic functions $S$ mapping the unit disk into the space of
contraction operators between two Hilbert spaces.  There are a number
of alternate characterizations:  the operator of multiplication by $S$
defines a contraction operator between two Hardy Hilbert spaces, $S$
satisfies a von Neumann inequality, a certain operator-valued kernel
associated with $S$ is positive-definite, and $S$ can be realized as
the transfer function of a dissipative (or even conservative)
discrete-time linear input/state/output linear system.  Various
multivariable generalizations of this class have appeared recently,
one of the most encompassing being that of Muhly and Solel where the
unit disk is replaced by the strict unit ball of the elements of a
dual correspondence $E^{\sigma}$ associated with a
$W^{*}$-correspondence $E$ over a $W^{*}$-algebra $\cA$ together with a
$*$-representation $\sigma$ of $\cA$.  The main new point which we
add here is the
introduction of the notion of reproducing kernel Hilbert correspondence
and identification of the Muhly-Solel Hardy spaces as reproducing kernel
Hilbert correspondences associated with a completely positive
analogue of the classical Szeg\"o kernel.  In this way we are able to
make the analogy between the Muhly-Solel Schur class and the
classical Schur class more complete.  We also illustrate the theory
by specializing it to some well-studied special cases; in some
instances there result new kinds of realization theorems.
\end{abstract}

\title[Multivariable generalizations of the Schur class]
{Multivariable generalizations of the Schur class: positive kernel
characterization and transfer function realization}
\author[J.A. Ball]{Joseph A. Ball}
\address{Department of Mathematics,
Virginia Tech,
Blacksburg, VA 24061-0123, USA}
\email{ball@math.vt.edu}

\author[A. Biswas]{Animikh Biswas}
\address{Department of Mathematics,
UNCC, NC 24061-0123, USA}
\email{abiswas@uncc.edu}

\author[Q. Fang]{Quanlei Fang}
\address{Department of Mathematics,
Virginia Tech,
Blacksburg, VA 24061-0123, USA}
\email{qlfang@math.vt.edu}

\author[S. ter Horst]{Sanne Ter Horst}
\address{Afdeling  Wiskunde, Faculteit der Exacte Wettenschappen,
Vrije Universiteit, De Boelelaan 1081a, 1081HV Amsterdam, The Netherlands}
\email{terhorst@few.vu.nl}

\subjclass{47A57}
\keywords{contractive multiplier, von Neumann inequality, completely
positive definite kernel, Hilbert $C^{*}$-module,
reproducing kernel Hilbert correspondence}

\maketitle

%%%%%%%%%%%%%%%%%%%%%%%%%%%%%%%%%%%%%%%%%%%%%%%%%%%%%%%%%%%%%%%%%%%%%%%%
%%%%%%%%%%%%%%%%%%%%%%%%%%%%%%%%%%%%%%%%%%%%%%%%%%%%%%%%%%%%%%%%%%%%%%%%
\section{Introduction} \label{S:intro}
\setcounter{equation}{0}

The classical Schur class $\cS$ (consisting of holomorphic
functions mapping the unit disk $\D$ into the closed unit disk
$\overline{\D}$) along with its operator-valued generalization has
been an object of intensive study over the past century (see
\cite{schur} for the original paper of Schur and
\cite{Schur-survey} for a survey of some of the impact and
applications in signal processing). To formulate the definition of
the operator-valued version, we let $\cL(\cU, \cY)$ denote the
space of bounded linear operators acting between Hilbert spaces $\cU$
and $\cY$. We also let $H^2_{\cU}({\D})$ and $H^{2}_{\cY}({\mathbb
D})$
be the standard Hardy spaces of $\cU$-valued (respectively
$\cY$-valued) holomorphic functions on the unit disk
${\D}$. By the Schur class ${\cS}(\cU, \cY)$ we mean the set of
$\cL(\cU, \cY)$-valued functions holomorphic on the unit disk
${\D}$ with values $S(z)$ having norm at most $1$ for each $z \in
{\D}$. The class $\cS(\cU, \cY)$ admits several remarkable
characterizations. The following result is well known and is
formulated as the prototype for the multivariable generalizations
to follow.

\begin{theorem} \label{T:C}
Let $S$ be an $\cL(\cU, \cY)$-valued function defined
on $\D$.  The following are equivalent:
\begin{enumerate}
\item $S \in \cS(\cU, \cY)$, i.e., $S$ is analytic on ${\D}$ with
contractive values in $\cL(\cU, \cY)$.
   \item[(1$'$)] The multiplication operator $M_{S} \colon f(z) \mapsto
S(z) \cdot f(z)$ defines a contraction from
    $H^{2}_{\cU}({\D})$ into $H^{2}_{\cY}({\D})$.
   \item[(1$''$)] $S$ is analytic and satisfies the von Neumann's
inequality: {\em if $T$ is  any strictly contractive operator on a
Hilbert space  ${\cal K}$, i.e., $\|T\| < 1$, then $S(T)$ is a
contraction operator ($\|S(T)\| \le 1$), where $S(T)$ is the
operator defined by
$$
S(T) = \sum_{n=0}^{\infty} S_{n} \otimes T^n \in {\cL}(\cU \otimes
{\cK}, \cY \otimes{\cK})
\quad\text{if}\quad S(z) = \sum_{n=0}^{\infty} S_{n}z^{n}.$$}

   \item The function $K_{S} \colon {\mathbb D} \times {\mathbb D} \to
   \cL(\cY)$ given by
\begin{equation*}  %\label{KS}
    K_{S}(z, w) = \frac{I_{\cY} - S(z) S(w)^{*}}{1 -
        z \overline{w}}
\end{equation*}
is a positive kernel on $\D\times\D$.
\item There exists a Hilbert space $\cH$ and a
coisometric (or even unitary or contractive) connecting operator
(or colligation) $\bU$ of the form
\begin{equation*}
%\label{1.1}
\bU = \begin{bmatrix} A  & B \\ C & D \end{bmatrix}  \colon
\begin{bmatrix} \cH \\ \cU
\end{bmatrix} \to \begin{bmatrix}\cH \\ \cY \end{bmatrix}
\end{equation*}
so that $S(z)$ can be realized in the form
\begin{equation}
S(z)= D + zC (I_{\cH} - zA)^{-1}B. \label{1.2}
\end{equation}
\end{enumerate}
\end{theorem}

   From the point of view of systems theory, the function
\eqref{1.2} is the
{\em transfer function} of the linear system
$$
\boldsymbol{\Sigma} =\boldsymbol{\Sigma}(\bU) \colon \left\{
\begin{array}{rcl}
    x(n+1) & = & A x(n) + B u(n) \\
    y(n) & = & C x(n) + D u(n)
\end{array}.  \right.
$$

The following well-known Proposition gives several equivalent
definitions for the term ``positive kernel'' used in condition (2) in
Theorem \ref{T:C}. The scalar
case ($\cY = {\mathbb C}$) of this result goes back to
the paper of Aronszajn \cite{aron}, but is also often attributed to
E.H.~Moore and Kolmogorov, while the vector-valued case has been well
exploited in the function-theoretic operator theory literature over
the years (see \cite{NF, dBR}).

\begin{proposition} \label{P:posker}
      Let $K \colon \Omega \times \Omega \to \cL(\cY)$
      be a given function.  Then the following conditions are
      equivalent:
\begin{enumerate}
\item
For any finite collection of points $\omega_{1}, \dots, \omega_{N}
\in \Omega$ and of vectors $y_{1}, \dots, y_{N} \in \cY$ ($N=1,2,\dots$) it holds that
\begin{equation}  \label{criterion1}
   \sum_{i,j = 1, \dots, N} \langle K(\omega_{i},\omega_{j}) y_{j}, y_{i}
   \rangle_{\cY} \ge 0.
   \end{equation}

\item There exists an operator-valued function $H \colon \Omega  \to
\cL(\cH,\cY)$
for some auxiliary Hilbert space $\cH$ so that
\begin{equation} \label{criterion2}
K(\omega', \omega) = H(\omega') H(\omega)^{*}.
\end{equation}

\item There exists a Hilbert space $\cH(K)$ of $\cY$-valued functions
$f$ on $\Omega$ so that the function $K(\cdot, \omega) y$ is in $\cH(K)$ for each
$\omega \in \Omega$ and $y \in \cY$ and has the reproducing property
\begin{equation*}  %\label{criterion3}
\langle f, K(\cdot, \omega) y \rangle_{\cH(K)} = \langle f(\omega),y
\rangle_{\cY}.
\end{equation*}
\end{enumerate}
When any (and hence all) of these equivalent conditions hold, we
say that {\em $K$ is a positive kernel on $\Omega \times \Omega$}.
\end{proposition}

We provide a sketch of the proof of
Theorem \ref{T:C} as a model for how extensions to more general
settings may proceed.

\begin{proof}[Sketch of the proof of Theorem \ref{T:C}]
The easy part is  (3) $\Longrightarrow$ (2) $\Longrightarrow$ (1$''$)
$\Longrightarrow$ (1$'$) $\Longrightarrow$ (1):

\textbf{(3) $\Longrightarrow$ (2):} Assume that $S(z)$ is as in
\eqref{1.2} with ${\mathbf U}$ unitary, and hence, in particular,
coisometric.  From the relations arising from the coisometric
property of ${\mathbf U}$:
$$
   \begin{bmatrix} A  & B \\ C & D \end{bmatrix}  \begin{bmatrix} A^* &
C^* \\ B^* & D^* \end{bmatrix} = \begin{bmatrix} I  & 0 \\ 0 & I
\end{bmatrix}
$$
one can verify that
\begin{align*}
I - S(z) S(w)^{*} & =  I - [D + zC (I - zA)^{-1}B][D + wC (I -
wA)^{-1}B]^{*}
\\ & =  C (I - zA)^{-1}[(1-z\overline{w})I_{\cH}] (I -
\overline{w}A^*)^{-1}C^*.
\end{align*}
This implies that $H(z)=C (I - zA)^{-1}$ satisfies
(\ref{criterion2}).

\textbf{(2) $\Longrightarrow$ (1$''$):} Due to $I - S(z) S(w)^{*} =
H(z)[(1-z\overline{w})I_{\cH}]H(w)^*$, we can see that for any $\|T\|
<1$
$$
I - S(T) S(T)^{*} = H(T)[(1-TT^*)\otimes I_{\cH}]H(T)^* \ge 0.
$$

\textbf{(1$''$) $\Longrightarrow$ (1$'$):} Observe that $M_S = S
(\bS) = s-\mathop {\lim }\limits_{r \uparrow 1 }S(r\bS)$
where $\bS$ is the shift operator $M_z$ on $H^2(\D)$. Thus the
fact that $\|S(r \bS)\| \le 1 $ for any $r<1$ implies $\| M_S \|
\le 1$.

\textbf{(1$'$) $\Longrightarrow$ (1):} Note that since $S(z)u = M_S \cdot
u$ for any
$u \in \cU$, we have $\|M_S\| _{op} = \|S\|_{\infty}$. So $\|M_S\|
\le 1$ implies that $S \in \cS(\cU, \cY) $.

The harder part is (1) $\Longrightarrow$ (1$'$) $\Longrightarrow$
(1$''$)
$\Longrightarrow$ (2) $\Longrightarrow$ (3):

\textbf{(1) $\Longrightarrow$ (1$'$):}  We can view $H^2(\D) \subset
L^2(\mathbf T)$. Thus $\| M_{S}u\|_{L^2(\mathbf T)} \le
\|S\|_{\infty} \cdot \| u\|_{L^2(\mathbf T)}$.

\textbf{(1$'$) $\Longrightarrow$ (1$''$):} According to the
Sz.-Nagy dilation theorem, any contraction operator $T$ has a
unitary dilation $U$.  In the strictly contractive case
$\|T\| < 1$, one can show that in fact the unitary dilation is the
bilateral shift with some multiplicity $N$: $ U = {\mathbf S}
\otimes I_{N}$ (if $N = \infty$, we interpret $I_{N}$
as the identity operator on $\ell^{2}$). We then have
$ T^n = P_K(\bS \otimes I_N)^n |_{ {\cal K}}$.  Therefore
$\|S(T)\| = \| P_{\cY \otimes  \cal K} S(\bS \otimes I_N)
|_{\cU \otimes  \cal K}\| \le \|M_S\| \le 1$.

\textbf{(1$''$) $\Longrightarrow$ (2):}
A direct proof of this implication can be done via a rather long,
intricate argument using a
Gelfand-Naimark-Segal construction in conjunction with a Hahn-Banach
separation argument---we refer to this as a GNS/HB argument.
For the polydisk setting, the argument
originates in \cite{agler-hellinger}; the version for a general
semigroupoid setting in \cite{DMM} covers in particular the classical
setting here.
\vskip .1in

Alternatively, one can avoid the GNS/HB
argument via the following shortcut:

\textbf{(1$''$) $\Longrightarrow$ (1$'$) $\Longrightarrow$ (2):} We
have seen that  (1$''$) $\Longrightarrow$ (1$'$) is easy. For (1$'$)
$\Longrightarrow$ (2), we assume $\|M_S \| \le 1$. View
$H^2(\D)$ as the reproducing kernel Hilbert space
$\cH(k_{\text{Sz}})$, where $k_{\text{Sz}}(z,w)=
\frac{1}{1-z\overline{w}}$ is the Szeg\"o kernel. Since
$M_{S}^{*}k_{S_z}(\cdot,w) y = k_{S_z}(\cdot,w)S(w)^{*} y$, we see
that
$$\sum_{i,j = 1, \dots, N} \langle K_{S}(z_{i},z_{j}) y_{j}, y_{i}
\rangle_{\cY} = \| \sum_{j} k_{\text{Sz}}(\cdot,z_j)y_j \|^2 -
\|(M_S)^* \sum_{j}k_{\text{Sz}}(\cdot,z_j)y_j \|^2 \ge 0$$
and it follows (via criterion \eqref{criterion1}) that $K_{S}$ is a
positive kernel on ${\mathbb D} \times {\mathbb D}$.

   \textbf{ (2) $\Longrightarrow$ (3):} This implication can be done by the
   now standard lurking isometry argument---see \cite{Ball-Winnipeg}
   where this coinage was introduced.
    \end{proof}

The purpose of this paper is to study recent extensions of the Schur
class and the associated analogues of Theorem \ref{T:C} to more
general multivariable settings.  In Section \ref{S:genSchur} we
describe two such extensions: the Drury-Arveson space setting and
the free-semigroup setting.  We emphasize how all the ingredients of
the proof of Theorem \ref{T:C} sketched above have direct analogues
in these two settings; hence the proof of the analogues of Theorem
\ref{T:C} for these two settings (see Theorem \ref{T:BTV} and
Theorem \ref{T:NC1} below) directly parallel the proof of Theorem
\ref{T:C} as sketched above.  A far more sophisticated generalized
Schur class has been introduced by Muhly and Solel (see
\cite{MS-Annalen, MS-Schur}).  The main contribution of the present
paper is to introduce the notion of reproducing kernel Hilbert
correspondence and an analogue of the Fourier (or $Z$-) transform
for the Muhly-Solel setting.  The starting point for most of the
constructions is a $W^{*}$-correspondence $E$ over a $W^{*}$-algebra
$\cA$ together with a $*$-representation $\sigma$ of $\cA$.  We show
that the image, denoted in our notation as $H^{2}(E, \sigma)$ which
is an analogue of $H^{2}$,  of a Muhly-Solel Fock space, denoted as
$\cF^{2}(E, \sigma)$ in our notation which is an analogue of
$\ell^{2}({\mathbb Z}_{+})$, under this $Z$-transform is a space of
$\cE$-valued functions ($\cE$ equal to a coefficient Hilbert space)
on the Muhly-Solel generalized unit disk ${\mathbb
D}((E^{\sigma})^{*})$\footnote{In nice cases, the general situation
collapses to this statement; more correctly, the vector-valued
functions are defined on ${\mathbb D}((E^{\sigma})^{*}) \times
\sigma(\cA)'$ where $\sigma(\cA)'$ is the commutant of the image
$\sigma(\cA)$ of $\sigma$ in $\cL(\cE)$.} and that an element $S$ of
the Muhly-Solel Schur class as introduced in \cite{MS-Schur} induces
a bounded multiplication operator on $H^{2}(E, \sigma)$.  We also
obtain analogues of the other parts of Theorem \ref{T:C} for this
setting (see Theorem \ref{T:MS-Schur} in Section 5 below) and thus
obtain a more complete analogy between the Muhly-Solel Schur class
and the classical Schur class than that presented in
\cite{MS-Schur}.   Section \ref{S:RKC} develops required
preliminaries concerning general correspondences, including the
notions of {\em reproducing kernel correspondence} and of {\em
reproducing kernel Hilbert correspondence}; these are natural
elaborations of the Kolmogorov decomposition for a completely
positive kernel found in \cite{BBLS}.  Section \ref{S:Hardyspaces}
introduces the spaces $H^{2}(E, \sigma)$ and $H^{\infty}(E, \sigma)$
which are the analogues of the Hardy spaces $H^{2}$ and $H^{\infty}$
for this setting.  The final section \ref{S:Examples} applies the
general theory to some familiar more concrete special cases.
Specifically we make explicit how the classical  case discussed
above as well as the Drury-Arveson setting and the free-semigroup
algebra setting discussed in Section \ref{S:genSchur} are particular
cases of the Muhly-Solel setting.  The general theory here leads to
more structured versions of these well-studied settings and
corresponding new types of realization theorems. We also discuss one
of the main examples motivating the work  in \cite{MS-JFA,
MS-Annalen, MS-Schur}, namely the setting of analytic
crossed-product algebras.  It is interesting to note that the
realization theorem for a particular instance of this example
amounts to the realization theorem for input-output maps of
conservative time-varying linear systems obtained in \cite{ABP}.

Another class of examples covered by the Muhly-Solel setting are
graph algebras (also known as semigroupoid algebras) \cite{pM97,
MS-CanJ, KP}; we do not discuss these here.
There are still other types of generalized Schur classes which are not
subsumed under the Muhly-Solel Fock space/correspondence setup.
We mention the Schur-Agler class for the polydisk (see
\cite{agler-hellinger,aglmccar-poly, BT} and for more general domains
\cite{at, bb4}), the noncommutative Schur-Agler class (see
\cite{BGM2, B-KV}), and higher-rank graph algebras (see \cite{KP06}).
A differentiating feature of these variants of the
Schur class is a more implicit version of condition (2) in Theorem
\ref{T:C} where the single positive kernel (the Szeg\"o kernel
$\frac{1}{1-z \overline{w}}$) is replaced by a whole family of
positive kernels.  An abstract framework using this feature as the
point of departure is the semigroupoid approach of
Dritschel-Marcantognini-McCullough \cite{DMM} which incorporates all
the aforementioned settings in \cite{agler-hellinger, aglmccar-poly,
BT, BGM2, KP06}.  However the theory in \cite{DMM} does not appear to
include the analytic crossed-product algebras included in the Muhly-Solel
scheme since it does not allow for the action of a $W^{*}$-algebra
$\cA$ acting on the ambient Hilbert space.  It is conceivable that
some sort of synthesis of these two disparate approaches is possible;
the recent work on product decompositions
over general semigroups (see \cite{Solel}) appears to be a start
in this direction.

The notation is mostly standard but we mention here a few
   conventions for reference.  For $\Omega$ any index set,
   $\ell^{2}(\Omega)$ denotes the space of complex-valued functions
   on $\Omega$ which are absolutely square summable:
   $$
    \ell^{2}(\Omega) = \{ \xi \colon \Omega \to {\mathbb C} \colon
    \sum_{\omega \in \Omega} | \xi(\omega) |^{2} < \infty\}.
   $$
   Most often the choice $\Omega = {\mathbb Z}$ (the integers) or
   $\Omega = {\mathbb Z}_{+}$ (the nonnegative
   integers) appears.
   For $\cH$ a Hilbert space, we use $\ell^{2}_{\cH}(\Omega)$ as
   shorthand for $\ell^{2}(\Omega) \otimes \cH$ (the space of
   $\cH$-valued function on $\Omega$ square-summable in norm).
   More general versions where $\cH$ may be a correspondence also
   come up from time to time.

   %%%%%%%%%%%%%%%%%%%%%%%%%%%%%%%%%%%%%%%%%%%%%%%%%%%%%%%%%%%%%%%%%%%%%%%%

   \section{Some multivariable Schur classes}  \label{S:genSchur}

   In this section we introduce two multivariable settings (the
   Drury-Arveson space setting and the free semigroup algebra setting)
   for the Schur class and formulate the analogue of Theorem \ref{T:C}
   for these two settings.

   \subsection{Drury-Arveson space}  \label{S:DA}
     A multivariable generalization of the Szeg\"o kernel
   $k(z,w)=(1-z\bar{w})^{-1}$ much studied of
   late is the positive kernel
   \begin{equation*}
   k_d(z,w)=\frac{1}{1-\langle z,  w  \rangle} \text{ on } {\mathbb B}^{d} \times
   {\mathbb B}^{d},
   %\label{1.1a}
   \end{equation*}
   where $\B^d=\left\{z=(z_1,\dots, z_d)\in\C^d \colon  \langle z,
   z\rangle<1\right\}$ is the unit ball of the  $d$-dimensional Euclidean
   space $\C^d$. By
   $\langle z,  w \rangle=\sum_{j=1}^d z_j \overline{w}_j$
   we mean the standard inner product in $\C^d$.
   The reproducing kernel Hilbert space (RKHS) $\cH(k_d)$ associated
   with
   $k_d$ via Aronszajn's construction \cite{aron} is a natural
   multivariable
   analogue of the Hardy space $H^2$ of the unit disk and coincides with
   $H^2$ if $d=1$.

   For $\cY$ an auxiliary Hilbert space, we consider the tensor product
   Hilbert space $\cH_\cY(k_d):=\cH(k_d)\otimes\cY$ whose
   elements can be viewed as  $\cY$-valued functions in $\cH(k_d)$.
   Then $\cH_\cY(k_d)$ can be characterized as follows:
   \begin{equation*}
   \cH_{\cY}(k_d)=\left\{f(z)=\sum_{\bn \in{\mathbb Z}^d_{+}}f_{\bn}
   z^\bn:\|f\|^{2}=\sum_{\bn \in {\mathbb Z}^d_{+}}
   \frac{\bn!}{|\bn|!}\cdot \|f_{\bn}\|_{\cY}^2<\infty\right\}.
   %\label{char}
   \end{equation*}
   Here and in what follows, we use standard multivariable notations: for
   multi-integers $\bn =(n_{1},\ldots,n_{d})\in \mathbb Z_+^d$ and points
   $z=(z_1,\ldots,z_d)\in\C^d$ we set
   \begin{equation*}
   |\bn| = n_{1}+n_{2}+\ldots +n_{d},\qquad
   \bn!  = n_{1}!n_{2}!\ldots n_{d}!, \qquad
   z^\bn = z_{1}^{n_{1}}z_{2}^{n_{2}}\ldots
   z_{d}^{n_{d}}.
   %\label{mnot}
   \end{equation*}

   By ${\mathcal M}_d(\cU,\cY)$ we denote the space of all
   $\cL(\cU,\cY)$-valued analytic functions  $S$ on $\B^d$ such that the
   induced multiplication operator
   \begin{equation*}
   M_S: \; f(z)\to S(z)\cdot f(z)
   %\label{ms}
   \end{equation*}
   maps $\cH_{\cU}(k_d)$ into $\cH_{\cY}(k_d)$. It follows by the closed
   graph theorem that for every $S\in {\mathcal M}_d(\cU,\cY)$, the
   operator $M_S$ is  bounded. We shall pay particular attention to the
   unit
   ball of  $ {\mathcal M}_{d}(\cU, \cY)$, denoted by
   $$
   {\mathcal S}_{d}(\cU, \cY) = \{ S \in {\mathcal M}_{d}(\cU,\cY) \colon
   \| M_{S} \|_{\text{op}} \le 1 \}.
   $$
   We refer to ${\mathcal S}_{d}(\cU, \cY)$ as a generalized
   ($d$-variable) {\em Schur class} since ${\mathcal S}_{1}(\cU, \cY)$
   collapses to the classical Schur class. Characterizations of
   ${\mathcal S}_{d}(\cU, \cY)$ in terms of realizations originate in
   \cite{aglmccar, BTV, EP}. The following is the analogue of
   Theorem \ref{T:C} for this
   setting; the result with condition (1$^{\prime \prime}$) eliminated
   appeared e.g.~in \cite{BTV, BBF2a}.

   \begin{theorem}
   \label{T:BTV}
   Let $S$ be an $\cL(\cU, \cY)$-valued function defined
   on ${\mathbb B}^{d}$.  The following are equivalent:
   \begin{enumerate}
   \item[(1$'$)] $S$ belongs to $\cS_d(\cU, \cY)$, i.e., the
   multiplication operator $M_{S} \colon f(z) \mapsto S(z)f(z)$
   defines a contraction from $\cH_{\cU}(k_{d})$ into $\cH_{\cY}(k_{d})$.

     \item[(1$''$)] $\|S(\mathbf T)\| \le 1$ for any commutative row
     contraction  $ T = (T_1,
   \cdots,T_d) \in \cL({ \cal K})^d$, i.e., if
   $S$ is given by $S(z)=\sum_{\bn\in{\mathbb
   Z}_{+}^d}S_{\bn}z^{\bn}$ and if
   $(T_{1}, \dots, T_{d})$ is any commuting $d$-tuple of bounded
   linear operators on a Hilbert space $\cK$ such that the row matrix
   $\begin{bmatrix} T_{1} &
   \cdots & T_{d}
   \end{bmatrix}$
   defines a strict contraction operator from $\cK^{d}$  to $\cK$, then
   the operator $S(\mathbf T) \in {\cL}({\cal U}\otimes {\cK}, {\cY}\otimes
     \cK)$
   defined via the operator-norm limit of
   the series
   $S(T):= \sum_n S_n \otimes {\mathbf T}^{\mathbf n}$
   has $\|S(T)\| \le 1$.

   \item[(2)] The function $K_{S} \colon {\mathbb B} \times {\mathbb B}
   \to \cL(\cY)$ given by
   \begin{equation*}  %\label{ball-KS}
      K_{S}(z, w) = \frac{I_{\cY} - S(z) S(w)^{*}}{1 -
      \langle z, w \rangle}
   \end{equation*}
   is a positive kernel (see Proposition \ref{P:posker}).

   \item[(3)] There exists a Hilbert space $\cH$ and a
   unitary (or even coisometric or contractive)
   connecting operator (or colligation) $\bU$ of the form
   \begin{equation*}
   %\label{1.7a}
   \bU = \begin{bmatrix} A  & B \\ C & D \end{bmatrix} =
   \begin{bmatrix} A_{1} & B_{1} \\ \vdots & \vdots \\ A_{d} & B_{d}
   \\ C & D \end{bmatrix} \colon \begin{bmatrix} \cH \\ \cU
   \end{bmatrix} \to \begin{bmatrix}\cH^{d} \\ \cY \end{bmatrix}
   \end{equation*}
   so that $S(z)$ can be realized in the form
   \begin{eqnarray*}
   S(z)&=&D+C\left(I_{\cH}-z_1A_1-\cdots-z_dA_d\right)^{-1}
   \left(z_1B_1+\ldots+z_dB_d\right)\nonumber\\
   &=& D + C (I - Z(z) A)^{-1} Z(z) B
   %\label{1.5a}
   \end{eqnarray*}
   where we set
   \begin{equation}  \label{1.6a}
   Z(z)=\begin{bmatrix}z_1 I_{\cH} & \ldots &
   z_dI_{\cH}\end{bmatrix},\quad
   A=\begin{bmatrix} A_1 \\ \vdots \\ A_d\end{bmatrix},\quad
   B=\begin{bmatrix} B_1 \\ \vdots \\ B_d\end{bmatrix}.
   \end{equation}
   \end{enumerate}
   \end{theorem}

   \begin{proof}[Remarks on the proof : ]  (3) $\Longrightarrow$ (2)
   $\Longrightarrow$ (1$''$) $\Longrightarrow$ (1$'$) follows in the
   same way as in the sketch of the proof of Theorem \ref{T:C} above.
   For (1$'$) $\Longrightarrow$ (2), one can use the same reproducing
   kernel argument as the shortcut discussed in the proof Theorem
   \ref{T:C} above. For (1$'$) $\Longrightarrow$ (1$''$), one can
   follow the corresponding argument sketched above for Theorem
   \ref{T:C} but with the
   Sz.-Nagy dilation theorem replaced with the
   Drury dilation theorem (see \cite{Drury}).
   The implication (2) $\Longrightarrow$ (3) follows
   exactly as in the classical case via the lurking isometry argument
   (see \cite{BTV}).
   Note that (1$''$) $\Longrightarrow$ (2)
   also can be achieved directly by the GNS/HB argument in \cite{DMM}
   specialized to the setting here, but
   this is not usually done since one has the alternative easier route
   (1$''$) $\Longrightarrow$ (1$'$) $\Longrightarrow$ (2).
   \end{proof}

   %%%%%%%%%%%%%%%%%%%%%%%%%%%%%%%%%%%%%%%%%%%%%%%%%%%%%%%%%%%%%%%%%%%%%%%%
   \subsection{Free semigroup algebras} \label{S:free}

   We now discuss the generalization of the Schur class associated with
   free semigroup algebras and models for row contractions (see
   \cite{PopescuNF1, PopescuNF2, Popescu-Nehari,
   Popescu-Memoir, Cuntz-scat}). We follow the formalism and notation as
   used in \cite{Cuntz-scat, NFRKHS}.

   Let $z = (z_{1}, \dots, z_{d})$ and $w = (w_{1}, \dots, w_{d})$
   be two sets of noncommuting indeterminates.  We let $\cF_{d}$ denote
    the free semigroup generated by the $d$ letters $\{1, \dots,
   d\}$.
    A generic element of $\cF_{d}$ is a word $w$ equal to a string of
    letters
    \begin{equation} \label{word}
     \alpha = i_{N} \cdots i_{1}\quad \text{where}\quad i_{k} \in
   \{1,
   \dots, d\} \;    \text{ for } \; k=1, \dots, N.
     \end{equation}
     The product of two words is defined by the usual concatenation.
   The unit element of $\cF_{d}$ is the {\em empty word} denoted by
   $\emptyset$. For $\alpha$ a word of the form \eqref{word}, we let
   $z^{\alpha}$ denote the monomial in noncommuting indeterminates
   $z^{\alpha} = z_{i_{N}} \cdots z_{i_{1}}$ and we let $z^{\emptyset} =
   1$. We extend this noncommutative functional calculus to a $d$-tuple
   of
   operators ${\mathbf A} = (A_{1}, \dots, A_{d})$ on a Hilbert space
   $\cK$:
     \begin{equation*}  %\label{bAv}
    {\mathbf A}^{v} = A_{i_{N}} \cdots A_{i_{1}}\quad \text{if}\quad
   v = i_{N} \cdots i_{1} \in \cF_{d} \setminus \{ \emptyset\};\quad
   {\mathbf
     A}^{\emptyset} = I_{\cK}.
     \end{equation*}
     We will also have need of the {\em transpose operation} on $\cF_{d}$:
   \begin{equation*} %\label{transpose}
      \alpha^{\top} = i_{1} \cdots i_{N} \quad\text{if}\quad \alpha =
   i_{N}\cdots i_{1}.
   \end{equation*}
   Given a coefficient Hilbert space
   $\cY$ we let $\cY\langle z \rangle$ denote the set of all polynomials
   in
   $z = (z_{1}, \dots, z_{d})$ with coefficients in $\cY:$
   $$
        \cY\langle z \rangle = \left\{ p(z) = \sum_{\alpha \in \cF_{d}}
        p_{\alpha} z^{\alpha} \colon p_{\alpha} \in \cY \text{ and }
        p_{\alpha} = 0 \text{ for all but finitely many } \alpha \right\},
        $$
        while $\cY \langle \langle z \rangle \rangle$ denotes the set of
   all
        formal power series in the indeterminates $z$ with coefficients in
        $\cY$:
   $$ \cY\langle \langle z \rangle \rangle = \left\{ f(z) = \sum_{\alpha
   \in
        \cF_{d}} f_{\alpha} z^{\alpha} \colon f_{\alpha} \in \cY \right\}.
        $$
   Note that vectors in $\cY$ can be considered as Hilbert space
   operators between ${\mathbb C}$ and $\cY$.  More generally, if $\cU$
   and $\cY$ are two Hilbert spaces, we let
   $\cL(\cU, \cY)\langle z \rangle$ and $\cL(\cU, \cY)\langle \langle z
   \rangle \rangle$ denote the space of polynomials (respectively,
   formal power series) in the noncommuting indeterminates $z = (z_{1},
   \dots, z_{d})$ with coefficients in $\cL(\cU, \cY)$.
   Given $S = \sum_{\alpha \in \cF_{d}} s_{\alpha}
   z^{\alpha} \in \cL(\cU, \cY)\langle \langle z \rangle
   \rangle$ and $f = \sum_{\beta \in \cF_{d}} f_{\beta} z^{\beta} \in
   \cU\langle \langle z \rangle \rangle$, the product $S(z) \cdot f(z)
   \in \cY \langle \langle z \rangle \rangle$ is defined as an element
   of $\cY\langle \langle z \rangle \rangle$ via the noncommutative
   convolution:
   \begin{equation} \label{multiplication}
   S(z) \cdot f(z) = \sum_{\alpha, \beta \in \cF_{d}} s_{\alpha}
   f_{\beta} z^{\alpha \beta} =
   \sum_{v \in \cF_{d}} \left( \sum_{\alpha, \beta \in \cF_{d} \colon
   \alpha \cdot \beta = v} s_{\alpha} f_{\beta} \right) z^{v}.
   \end{equation}
   Note that the coefficient of $z^{v}$ in \eqref{multiplication}
   is well defined since any
   given word $v \in \cF_{d}$ can be decomposed as a product $v = \alpha
   \cdot \beta$ in only finitely many distinct
   ways.

   In general, given a coefficient Hilbert space $\cC$, we use the
   $\cC$ inner product to generate a pairing
   $$ \langle \cdot, \, \cdot \rangle_{\cC \times \cC\langle \langle w
   \rangle \rangle} \colon \cC \times \cC\langle \langle w \rangle
   \rangle \to {\mathbb C}\langle \langle w \rangle \rangle
   $$
   via
   $$
   \left\langle c, \sum_{\beta \in \cF_{d}} f_{\beta} w^{\beta}
   \right\rangle_{\cC \times
   \cC\langle \langle w \rangle \rangle} = \sum_{\beta \in \cF_{d}}
   \langle
   c, f_{\beta} \rangle_{\cC} w^{\beta^{\top}} \in {\mathbb C} \langle \langle
   w \rangle \rangle.
   $$
   Similarly we can consider
   $\left\langle \sum_{\alpha \in \cF_{d} } f_{\alpha} w^{\alpha}, c
   \right\rangle_{\cC\langle \langle w \rangle \rangle \times \cC}$
   or the more general pairing
   $$
   \left\langle \sum_{\alpha \in \cF_{d}} f_{\alpha} w^{\prime \alpha},
   \sum_{\beta \in \cF_{d}} g_{\beta} w^{\beta}
   \right\rangle_{\cC\langle \langle w' \rangle \rangle \times \cC
   \langle \langle w \rangle \rangle} =
   \sum_{\alpha, \beta \in \cF_{d}} \langle f_{\alpha}, g_{\beta}
   \rangle_{\cC} w^{\beta^{\top}} w^{\prime \alpha}.
   $$
   Suppose that $\cH$ is a Hilbert space whose elements are formal power
   series in
   $\cY \langle \langle z \rangle \rangle$ and that $K(z,w) =
   \sum_{\alpha, \beta \in \cF_{d}} K_{\alpha, \beta} z^{\alpha}
   w^{\beta^{\top}}$ is a formal power series in the two sets of $d$
   noncommuting indeterminates $z = (z_{1}, \dots, z_{d})$ and $w =
   (w_{1}, \dots, w_{d})$.  We say that $\cH$ is a NFRKHS ({\em noncommutative
   formal reproducing kernel Hilbert space}) if for each $\alpha \in
   \cF_{d}$, the linear operator $\Phi_{\alpha} \colon \cH \to \cY$
   defined by $f(z) = \sum_{\beta \in \cF_{d}} f_{\beta} z^{\beta}
   \mapsto
   f_{\alpha}$ is continuous.  In this case there must be a formal power
   series $k_{\alpha}(z) \in \cL(\cY)\langle \langle z \rangle \rangle$
   so that $k_{\alpha}(\cdot)y \in \cH$ for each $\alpha \in \cF_{d} $
   and $y \in \cY$ and
   $$
      \langle f, k_{\alpha }y \rangle_{\cH} = \langle f_{\alpha}, y
      \rangle_{\cY}.
   $$
   If we set $K(z,w) = \sum_{\beta \in \cF_{d}}
   k_{\beta}(z)  w^{\beta^{\top}}$,  then we have the reproducing
   property
   $$
     \langle f, K(\cdot, w) y \rangle_{\cH \times \cH\langle \langle w
     \rangle \rangle} = \langle f(w), y \rangle_{\cY \langle \langle w
     \rangle \rangle \times \cY}.
   $$
   In this case we say that {\em $K(z,w)$ is the reproducing kernel for
   the NFRKHS $\cH$.}  As explained in detail in \cite{NFRKHS}, we have
   the following equivalent characterizations for such kernels which
   parallel the statements of Proposition
   \ref{P:posker}  for the classical case.

   \begin{proposition}  \label{P:NCposker}
        Let $K(z,w) \in \cL(\cY)\langle \langle z, w \rangle \rangle$ be
        a formal power series in two sets of noncommuting indeterminates
        with coefficients $K_{\alpha, \beta}$ equal to bounded operators
        on the Hilbert space $\cY$.  Then the following conditions are
        equivalent:
   \begin{enumerate}
        \item For all finitely supported $\cY$-valued functions $\alpha
     \mapsto y_{\alpha}$ it holds that
     \begin{equation*} %\label{NC-criterion1}
     \sum_{\alpha, \alpha' \in \cF_{d}} \langle K_{\alpha, \alpha'}
     y_{\alpha'}, y_{\alpha} \rangle \ge 0,
     \end{equation*}
      i.e., the function from $\cF_{d} \times \cF_{d}$ to $\cL(\cY)$
      given by $(\alpha, \beta) \mapsto K_{\alpha, \beta}$ is a positive
      kernel in the classical sense of Proposition \ref{P:posker}.

      \item $K(z,w)$ has a factorization
      \begin{equation*} %\label{NC-criterion2}
      K(z,w) = H(z) H(w)^{*}
      \end{equation*}
      for some $H \in \cL(\cH, \cY)\langle \langle z \rangle \rangle$
      where $\cH$ is some auxiliary Hilbert space.  Here
     $$
     H(w)^{*} = \sum_{\beta \in \cF_{d}} H_{\beta}^{*}
   w^{\beta^{\top}} =
     \sum_{\beta \in \cF_{d}} H_{\beta^{\top}}^{*} w^{\beta}\quad
   \text{if}
   \quad
     H(z) = \sum_{\alpha \in \cF_{d}} H_{\alpha} z^{\alpha}.
     $$
        \item  $K(z,w)$ is the reproducing
   kernel for a NFRKHS $\cH(K)$, i.e.,  for each $\beta  \in \cF_{d}$
   and $y \in \cY$ the formal power
   series $k_{\beta}y$ given by
   $ k_{\beta}y(z): = \sum_{\alpha \in \cF_{d}} K_{\alpha, \beta} y
   z^{\alpha}$ is in $\cH(K)$ and has the reproducing property
   \begin{equation*}  %\label{NC-criterion3}
    \langle f, \sum_{\beta \in \cF_{d}}k_{\beta}y w^{\beta}
    \rangle_{\cH(K) \times \cH(K)\langle \langle w
    \rangle \rangle} = \langle f(w), y \rangle_{\cY\langle \langle w
    \rangle \rangle \times \cY} \text{ for every } f\in\cH(K).
   \end{equation*}
   \end{enumerate}
   \end{proposition}

   A natural analogue of the vector-valued Hardy space over the unit disk
   (see e.g.~\cite{PopescuNF1})
   is the Fock space with coefficients in $\cY$ which we denote here by
   $H^{2}_{\cY}(\cF_{d})$ and express the elements in  power series form:
   \begin{equation}  \label{free-Hardy}
        H^{2}_{\cY}(\cF_{d}) = \left\{ f(z) = \sum_{\alpha \in \cF_{d}}
        f_{\alpha} z^{\alpha} \colon  f_{\alpha}\in \cY, \sum_{\alpha \in
   \cF_{d}} \|
        f_{\alpha}\|^{2} < \infty \right\}.
   \end{equation}
   When $\cY = {\mathbb C}$ we write simply $H^{2}(\cF_{d})$.
   As explained in \cite{NFRKHS}, $H^{2}(\cF_{d})$ is a NFRKHS with
   reproducing kernel  equal the following noncommutative analogue of the
   classical Szeg\"o kernel:
   \begin{equation}  \label{classical-kSz}
      k_{\text{Sz, nc}}(z,w) = \sum_{\alpha \in \cF_{d}} z^{\alpha}
      w^{\alpha^{\top}}.
   \end{equation}
   Thus we have in general
   $H^{2}_{\cY}(\cF_{d}) = \cH(k_{\text{Sz}}  \otimes I_{\cY})$.
   We abuse notation and let $S_{j}$ denote the shift operator
   \begin{equation*}  %\label{shift}
   S_{j} \colon f(z) = \sum_{v \in \cF_{d}} f_{v} z^{v} \mapsto
   f(z) \cdot z_{j} = \sum_{v \in \cF_{d}} f_{v} z^{v \cdot j} \text{
   for } j = 1, \dots, d
   \end{equation*}
   on $H^{2}_{\cY}(\cF_{d})$ for any auxiliary space $\cY$.  The adjoint of
   $S_{j}$ on $H^{2}_{\cY}(\cF_{d})$ is then given by
   \begin{equation*} %\label{bs}
     S_{j}^{*} \colon \sum_{v \in \cF_{d}} f_{v}z^{v} \mapsto
   \sum_{v \in
     \cF_{d}} f_{v \cdot j} z^{v}\quad \text{for}\quad j = 1, \dots,
   d.
   \end{equation*}

   We let $\cM_{nc,d}(\cU, \cY)$ denote the set of formal power series
   $S(z) = \sum_{\alpha \in \cF_{d}} s_{\alpha} z^{\alpha}$ with
   coefficients $s_{\alpha} \in \cL(\cU, \cY)$ such that the associated
   multiplication operator $M_{S} \colon f(z) \mapsto S(z) \cdot f(z)$
   (see \eqref{multiplication}) defines a bounded operator from
   $H^{2}_{\cU}(\cF_{d})$ to $H^{2}_{\cY}(\cF_{d})$. The noncommutative
   Schur class ${\mathcal S}_{nc, d}(\cU, \cY)$ is defined to
   consist of
   such multipliers $S$ for which $M_{S}$ has operator norm at most 1:
   \begin{equation*} %\label{ncSchur}
     {\mathcal S}_{nc, d}(\cU, \cY) = \{ S \in \cL(\cU,
\cY)\langle \langle z
    \rangle \rangle \colon
     M_{S} \colon H^{2}_{\cU}(\cF_{d}) \to H^{2}_{\cY}(\cF_{d})
   \text{
     with } \|M_{S}\|_{op} \le 1 \}.
   \end{equation*}
   The following is the noncommutative analogue of Theorem \ref{T:C}
   for this setting. We refer to \cite{PopescuNF1,PopescuNF2} for
   details.

   \begin{theorem}  \label{T:NC1} Let $S(z) \in \cL(\cU, \cY) \langle
    \langle z \rangle \rangle$
        be a formal power series in $z = (z_{1}, \dots, z_{d})$ with
   coefficients in $\cL(\cU, \cY)$.  Then the following are equivalent:
   \begin{enumerate}
      \item[(1$^{\prime}$)]
      $S \in \cS_{nc,d}(\cU, \cY)$, i.e., $M_{S} \colon \cU\langle
      z\rangle \to \cY \langle \langle z \rangle \rangle$ given by
      $M_{S} \colon p(z) \to S(z) p(z)$ extends to define a
      contraction operator from $H^{2}_{\cU}(\cF_{d})$ into
      $H^{2}_{\cY}(\cF_{d})$.

      \item[(1$^{\prime \prime}$)] For each strict row contraction
      $(T_{1}, \dots, T_{d})$, i.e., a $d$-tuple $(T_{1}, \dots,
      T_{d})$ of operators on a Hilbert space $\cK$ (commutative
      or not) such that the row matrix $\begin{bmatrix} T_{1}&
      \cdots & T_{d} \end{bmatrix}$ defines a strict contraction
      operator from $\cK^{d}$ to $\cK$, we have
      $$
        \|S(T)\| \le 1,
       $$
        where
        $$
        S(T) = \sum_{\alpha \in \cF_{d}} s_{\alpha}
        \otimes T^{\alpha} \in \cL(\cU \otimes \cK), \cY \otimes \cK)
        \quad \text{ if } \quad  S(z) = \sum_{\alpha \in \cF_{d}} s_{\alpha}
        z^{\alpha}
        $$
        and where we set
        $$ T^{\alpha} = T_{i_{N}} \cdots T_{i_{1}} \quad \text{ if }
        \quad
        \alpha = i_{N} \cdots i_{1} \in \cF_{d}.
        $$

      \item[(2)] The formal power series given by
      \begin{equation*}  %\label{NC-KS}
    K_{S}(z, w) : = k_{\text{Sz, nc}}(z,w) - S(z) k_{\text{Sz, nc}}(z,w)
      S(w)^{*}
      \end{equation*}
      is a noncommutative positive kernel (see Proposition
      \ref{P:NCposker}).

      \item[(3)] There exists a Hilbert space $\cH$ and a unitary
   connection
      operator $\bU$ of the form
      \begin{equation*} %\label{NCcolligation}
      \bU = \begin{bmatrix} A & B \\ C & D \end{bmatrix} =
      \begin{bmatrix} A_{1} & B_{1} \\ \vdots & \vdots \\ A_{d} &
      B_{d} \\ C & D \end{bmatrix} \colon \begin{bmatrix} \cX \\ \cU
      \end{bmatrix} \to \begin{bmatrix}  \cX \\ \vdots \\ \cX \\ \cY
      \end{bmatrix}
      \end{equation*}
      so that $S(z)$ can be realized as a formal power series in the
      form
      \begin{equation*}   %\label{NCrealization}
        S(z) = D + \sum_{j=1}^{d} \sum_{v \in \cF_{d}} C A^{v}B_{j}
        z^{v}\cdot z_{j}=
        D + C (I - Z(z) A)^{-1} Z(z) B
       \end{equation*}
   where $Z(z)$, $A$ and $B$ are as in \eqref{1.6a} but where now
   $z_{1}, \dots, z_{d}$ are noncommuting indeterminates rather than
   commuting variables.
     \end{enumerate}
   \end{theorem}

   \begin{proof}[Sketch of the proof of Theorem \ref{T:NC1}:]
      The proof of (3) $\Longrightarrow$ (2) $\Longrightarrow$
      (1$^{\prime \prime}$) $\Longrightarrow$ (1$^{\prime}$)
      $\Longrightarrow$ (1) formally goes through in the same was as the
      classical case.  Let us just note that (1$^{\prime \prime}$)
      $\Longrightarrow$ (1$^{\prime}$) involves viewing $M_{S} \colon
      H^{2}_{\cU}(\cF_{d}) \to H^{2}_{\cY}(k_{d})$ as $M_{S} = S({\mathbf
      S})$ where ${\mathbf S} = (S_{1}, \dots, S_{d})$ are the left
      creation operators of multiplication by $z_{j}$ on the left on the
      Fock space $H^{2}(\cF_{d})$.  From the assumption (1$^{\prime
      \prime}$), we know that $\| S(r {\mathbf S})\| \le 1$ for each $r <
      1$ and hence $\| M_{S} \| = \lim_{r \uparrow 1}\| S(r{\mathbf S})\|
      \le 1$ as well.

      We discuss the harder direction (1$^{\prime}$)
      $\Longrightarrow$ (1$^{\prime \prime}$) $\Longrightarrow$ (2)
      $\Longrightarrow$ (3).

     \textbf{(1$^{\prime}$) $\Longrightarrow$ (1$^{\prime \prime}$):}
      One can follow the proof of (1$^{\prime}$) $\Longrightarrow$
      (1$^{\prime \prime}$) for the classical case but substitute the
      Popescu dilation theorem for row contractions (see
      \cite{Popescu-CLT1}) for the Sz.-Nagy dilation theorem for a single
      contraction operator.

      \textbf{(1$^{\prime \prime}$) $\Longrightarrow$ (2):}  This
      implication again can be
      done via an appropriate version of the GNS/HB argument; see
      \cite{BGM2} for a slightly more general version and \cite{DMM} for
      an even more general version.

      Alternatively, one can follow the route \textbf{(1$^{\prime \prime}$)
      $\Longrightarrow$ (1$^{\prime}$) $\Longrightarrow$ (2):} As we have
      already discussed (1$^{\prime \prime}$) $\Longrightarrow$
      (1$^{\prime}$), it suffices to discuss (1$^{\prime}$)
      $\Longrightarrow$ (2). This can be done by an adaptation of
      the argument for the classical case to the present setting of
      formal, noncommutative reproducing
      kernel Hilbert spaces---see \cite[Theorem 3.15]{NFRKHS}.

      \textbf{(2) $\Longrightarrow$ (3):}  The lurking isometry argument
      works in this context as well---see \cite[Theorem 3.16]{NFRKHS}.
        \end{proof}

        %%%%%%%%%%%%%%%%%%%%%%%%%%%%%%%%%%%%%%%%%%%%%%%%%%%%%%%%%%%%%%%%%%%%%%%%
        %%%%%%%%%%%%%%%%%%%%%%%%%%%%%%%%%%%%%%%%%%%%%%%%%%%%%%%%%%%%%%%%%%%%%%%%
        \section{Reproducing kernel $(\cA, \cB)$-correspondences} \label{S:RKC}

        \setcounter{equation}{0} The notion of a vector-valued
        reproducing kernel Hilbert space based on an operator-valued
        positive kernel has been a standard tool in operator theory as
        well as in other applications for some time now. Recently,
        Barreto, Bhat, Liebscher and Skeide \cite{BBLS} introduced a finer
        notion of positive kernel ({\em completely positive kernel}) and
        gave several equivalent characterizations, but did not develop the
        connections with reproducing kernel Hilbert spaces.  The purpose
        of this section is to fill in this gap, as it is the natural tool
        for the discussion to follow.

        Let $\cB$ be a $C^{*}$-algebra and $E$ a linear space.  For some of the
        discussion to follow, it will be convenient to assume that
$\cB$ has a unit.
        However, any $C^{*}$-algebra has an {\em approximate identity} (see
        \cite[Theorem I.4.8]{Davidson}); by making use of such an approximate
        identity, most arguments using a unit element $1_{\cB}$ can be
        adapted to an approximation argument yielding the desired result for
        the general case where $\cB$ is not assumed to possess a unit.  In the
        sequel we usually leave the details of this adaptation to the reader.

        We say that $E$ is a
        {\em (right) pre-Hilbert  $C^{*}$-module} over $\cB$ if $E$ is a
        right
        module over $\cB$ and is endowed with a $\cB$-valued inner product
        $\langle \cdot, \cdot \rangle_{E}$
        satisfying the following axioms for any $\lambda, \mu \in {\mathbb
        C}$, $e,f,g \in E$ and $b \in \cB$:
        \begin{enumerate}
        \item $\langle \lambda e + \mu  f, g \rangle_{E} =
        \lambda \langle e, g \rangle_{E} + \mu \langle f, g \rangle_{E}$;
        \item  $\langle e \cdot b, f \rangle_{E} =
        \langle e,f \rangle_{E}   b$;
        \item $\langle e, f \rangle_{E}^{*} = \langle f, e \rangle_{E}$;
        \item $\langle e, e \rangle_{E} \ge 0$ (as an element of $\cB$);
        \item $\langle e, e \rangle_{E} = 0$ implies that $e = 0$.
        \end{enumerate}
        We also impose that $(\lambda e) \cdot b = e \cdot (\lambda b)$ for
        all $e \in E$, $b \in \cB$ and $\lambda \in {\mathbb C}$. Note that if
        $\cB$ has a unit, this last condition is automatic from the axioms
        for the identification $\lambda \mapsto \lambda \cdot 1_{\cB}$
        and the axioms for $E$ being
        a module over $\cB$.
        (Unlike some other authors, we take the $\cB$-valued inner-product to
        be linear in the first variable and conjugate-linear in the second
        variable as is usually done in the Hilbert-space setting ($\cB =
        {\mathbb C}$) rather than the reverse.) Note that it then follows that
        $$
     \langle e, f \cdot b \rangle_{E} = b^{*} \langle e, f \rangle_{E}.
        $$
        When the inner product is
        clear, we drop the subscript $E$ and write simply $\langle e, f
        \rangle$ for the $\cB$-valued inner product.
        If $E$ is a pre-Hilbert module over $\cB$, then $E$ is a normed
        linear
        space with norm given by
        \begin{equation}  \label{norm}
      \| e \|_{E} =  \|\langle e, e \rangle^{1/2}\|_{\cB}.
        \end{equation}
        Here $\|\ \|_{\cB}$ denotes the norm associated
        with the $C^*$-algebra $\cB$. One can always complete $E$
        to a Banach space in the norm \eqref{norm} to get what we
        shall call a {\em Hilbert $C^*$-module over $\cB$}.  Moreover, $E$
        has additional structure, namely $E$ carries the structure of an
        {\em operator space}, i.e., $E$ is the upper right corner of
        a subalgebra of operators acting on a Hilbert space with a
        representation as $2 \times 2$-block operator matrices (the {\em
        linking algebra})---see \cite{MS-JFA} or \cite{RW}.

        Given two Hilbert $C^{*}$-modules $E$ and $F$ over the same $C^{*}$-algebra
        $\cB$, it is natural to consider the space ${\mathcal L}(E,F)$ of
        bounded linear
        operators $T \colon E \to F$ between the Banach spaces $E$ and $F$.
        Unlike the Hilbert space case, for a linear map $T$ from $E$ to $F$
        it may or may not
        happen that there is
        an {\em adjoint operator} $T^* \in {\mathcal L}(F, E)$ so that
        $$
     \langle Te, f \rangle_{F} = \langle e, T^{*}f \rangle_{E}
        {\mbox{ for all }e \in E\mbox{ and }f \in F.}
        $$
        In case there exists an operator $T^* \in {\mathcal L}(F, E)$
        with this property we say that $T$ is {\em
        adjointable} and we denote the set of all adjointable linear
        operators between $E$ and $F$ as ${\mathcal L}^{a}(E,F)$ (with the
        usual abbreviation ${\mathcal L}^{a}(E)$ in case $E = F$).
        When the mapping $T \colon E \to F$ is adjointable in this sense,
        necessarily
        $T \in {\mathcal L}(E, F)$ with the additional property that $T$ is a
        $\cB$-module map:
        \begin{equation}  \label{module-map}
      T( e \cdot b) = T(e) \cdot b \text{ for all } e \in E \text{ and }
      b \in {\mathcal B}.
        \end{equation}
        However, this additional property \eqref{module-map} alone is not
        sufficient for
        admission of $T$ in the class ${\mathcal L}^{a}(E,F)$ of adjointable
        maps (see \cite[Example 2.19]{RW}).

        Following \cite{MS-JFA, MS-Annalen} (see also the books \cite{MT, RW}
        for more comprehensive treatments), we now introduce the notion of
        an $(\cA, \cB)$-correspondence.
        If $E$ is a right Hilbert
        $C^{*}$-module over $\cB$ and $\cA$ is another $C^{*}$-algebra,
        we say that $E$ is a {\em $(\cA, \cB)$-correspondence} if $E$ is also
        a left module over $\cA$ which makes $E$ an $(\cA, \cB)$-bimodule:
        $$
        (a \cdot  e) \cdot b = a \cdot (e \cdot b)
        \text{ for all } a \in \cA,\, e \in E \text{ and } b \in \cB
        $$
        with the additional compatibility condition
        \begin{equation}  \label{cor-axiom}
        \langle a \cdot e, \, f \rangle_{E} = \langle e,\, a^{*} \cdot f
        \rangle_{E}.
        \end{equation}
        The compatibility condition in (\ref{cor-axiom}) is
        equivalent to requiring that each of the left multiplication
        operators $ \varphi(a) \colon e \mapsto a \cdot e$ on $E$ is a
        bounded linear operator on $E$ for each $a \in \cA$ and $\varphi$
        is a $C^{*}$-homomorphism from $\cA$ into the $C^{*}$-algebra
        ${\mathcal L}^{a}(E)$ of bounded adjointable operators on $E$:
        thus $\varphi(a)$ is adjointable for each $a \in \cA$ with
        $\varphi(a)^{*} = \varphi(a^{*})$.  We shall occasionally write
        $\varphi(a) e$ rather than $a \cdot e$.

        Note the lack of symmetry
        in the roles of $\cA$ and $\cB$:  the identities $\langle e \cdot
        b, f \rangle = \langle e, f \rangle b$ together with $\langle e, f
        \cdot b^{*}\rangle = b \cdot \langle e, f \rangle$ preclude the
        validity in general of the identity $\langle e \cdot b, f \rangle
        = \langle e, f \cdot b^{*}\rangle $ (the would-be $\cB$ analogue
        of \eqref{cor-axiom}) unless $\cB$ is commutative.

        If both $\cA$ and
        $\cB$ have units, we also demand that the scalar multiplication on
        $E$ is compatible with both the identification $\lambda \mapsto
        \lambda 1_{\cA}$ of ${\mathbb C}$ as a subalgebra of $\cA$ and the
        identification $\lambda \mapsto \lambda 1_{\cB}$ of ${\mathbb C}$ as
        a subalgebra of $\cB$.  This is consistent with demanding the
        additional axioms
        $$
     (\lambda a) \cdot e = a \cdot (\lambda e), \qquad (\lambda e) \cdot
     b = e \cdot (\lambda b)
        $$
        for the general case.

        The classical
        case is the one where $E$ is a Hilbert space
        $\cE$, $\cB = {\mathbb C}$ and $\cA =
        \cL(\cE)$ with the operations given by
        \begin{align*}
      & a \cdot e = a e \text{ (the operator $a$ acting on the vector
      $e$)} \\
      & e \cdot b = b e \text{ (scalar multiplication in $\cE$)}, \\
      & \langle e, f \rangle\text{ (the $\cE$ Hilbert-space inner
      product).}
        \end{align*}
        Another easy example
        is to take $E = \cA = \cB$ all equal to a $C^{*}$-algebra with
        $$
        a \cdot e = ae, \qquad e \cdot b = eb, \qquad \langle e,f \rangle_{E}
        = f^{*}e.
        $$
        We encourage the reader to peruse Section \ref{S:Examples} for a
        variety of additional examples and references for more complete
        details.

        We will have need of various constructions for making new
        correspondences out of given correspondences.  We give formal
        definitions as follows.

        \begin{definition} \label{D:tensor}
        \begin{enumerate}
        \item {\rm \textbf{Direct sum:}} Suppose that $E$ and $F$ are two
        $(\cA, \cB)$-correspondences.  Then the direct-sum correspondence
        $E \oplus F$ is defined to be the direct sum vector space
        $E \oplus F$ together with the diagonal left-$\cA$ action and
        right-$\cB$ action and the direct-sum $\cB$-valued
        inner product:
      \begin{align*}
         & a \cdot (e \oplus f) = (a \cdot e) \oplus (a \cdot f), \qquad
          (e \oplus f) \cdot b = (e \cdot b) \oplus  (f \cdot b), \\
         & \langle e \oplus f, e' \oplus f' \rangle_{E \oplus F} =
          \langle e, e' \rangle_{E} + \langle f, f' \rangle_{F}.
          \end{align*}

          \item  {\rm \textbf{ Tensor product:}}
        Suppose that we are given three $C^{*}$-algebras $\cA, \cB$ and
        ${\mathcal C}$ together with an $(\cA, \cB)$-correspondence $E$
        and a $(\cB, {\mathcal C})$-correspondence $F$.  Then we define
        the tensor product correspondence $E \otimes_{\cB} F$
        (sometimes abbreviated to $E \otimes F$) to be the
        completion of the linear span of all tensors $e \otimes f$ (with
        $e \in E$ and $f \in F$) subject to the identification
        \begin{equation}\label{balance}
        (e \cdot b) \otimes f = e \otimes (b \cdot f),
        \end{equation}
     with left $\cA$-action given by
     $$ a \cdot (e \otimes f) = (a \cdot e) \otimes f,
     $$
     with right ${\mathcal C}$-action given by
     $$
      (e \otimes  f) \cdot c = e \otimes (f \cdot c),
     $$
     and with ${\mathcal C}$-valued inner product $\langle \cdot, \cdot
     \rangle_{E \otimes F}$ given by
     $$
     \langle e \otimes f,\, e' \otimes f' \rangle_{E \otimes F}
     = \left\langle \langle e, e' \rangle_{E} \cdot f, \, f'
     \right\rangle_{F}.
     $$
     \end{enumerate}
     \end{definition}

     It is a straightforward exercise to verify that the balanced
        tensor-product construction is
     well-defined.  For example the computation
     \begin{align*}
         \langle (e \cdot b) \otimes f, \, (e' \cdot b') \otimes f'
         \rangle & =
         \left \langle b^{\prime *} \cdot \langle e, e' \rangle \cdot b
        \cdot
         f, f' \right \rangle \\
         & = \left \langle \langle e,e'\rangle \cdot b \cdot f, \, b'
         \cdot f' \right \rangle \\
         & = \langle  e \otimes (b \cdot f), \, e' \otimes  (b' \cdot f')
         \rangle
     \end{align*}
     shows that the $E \otimes F$-inner product is well-defined.

        \begin{remark} {\em Bounded linear operators between direct sum
        correspondences admit operator matrix decompositions in precisely
        the same way as in the Hilbert space case ($\cB=\C$), while
        adjointability of such an operator corresponds to the operators in
        the decomposition being adjointable. For bounded linear operators
        between tensor-product correspondences the situation is slightly
        more complicated. We give an example how operators can be
        constructed. Let $E$ and $E'$ be $(\cA,\cB)$-correspondences and
        $F$ and $F'$ $(\cB,\cC)$-correspondences, for $C^*$-algebras
        $\cA$, $\cB$ and $\cC$. Furthermore, let $X\in\cL(E,E')$ and
        $Y\in\cL(F,F')$ be $\cB$-module maps. Then we write $X\otimes Y$
        for the operator in $\cL(E\otimes_\cB F,E'\otimes_\cB F')$ which
        is determined by
        \begin{equation}\label{tensoroperator}
        X\otimes Y (e \otimes f)=(Xe)\otimes (Yf)
        \mbox{ for each } e \otimes f \in E \otimes_\cB F.
        \end{equation}
        The $\cB$-module map properties are needed to guarantee that for each
        $e \otimes f \in E\otimes_\cB F$ and all $b\in\cB$ we have
\begin{eqnarray*}
X\otimes Y (eb\otimes f)&=&(X(eb))\otimes(Yf)=(Xe)b\otimes(Yf)
        =(Xe)\otimes b(Yf)=(Xe)\otimes(Y(bf))\\&=&X\otimes Y(e\otimes bf).
\end{eqnarray*}
 Thus the balancing in the tensor product (see \eqref{balance})
is respected by
        the operator $X\otimes Y$. Moreover, $X\otimes Y$ is adjointable in case
        $X$ and $Y$ are adjointable operators, with $(X\otimes
Y)^*=X^*\otimes Y^*$.
        Indeed, this is the case since for $f\otimes g\in E\otimes F$ and
        $e'\otimes f'\in E' \otimes F'$ we have
        \begin{eqnarray*}
        \inn{(X\otimes Y)(e\otimes f)}{e'\otimes f'}_{E'\otimes F'}
        &=&\inn{Xe\otimes Yf}{e'\otimes f'}_{E'\otimes F'}\\
        &=&\inn{\inn{Xe}{e'}_{E'}Yf}{f'}_{F'}\\
        &=&\inn{Y\inn{Xe}{e'}_{E'}f}{f'}_{F'}\\
        &=&\inn{\inn{e}{X^*e'}_{E}f}{Y^*f'}_{F}\\
        &=&\inn{e\otimes f}{X^*e'\otimes Y^* f'}_{E\otimes F}\\
        &=&\inn{e\otimes f}{(X^*\otimes Y^*)e'\otimes f'}_{E\otimes F}.
        \end{eqnarray*}
        In particular, the left action on $E\otimes F$ can now be written as
        $a\mapsto\varphi(a)\otimes I_F\in\cL^a(E\otimes F,E\otimes F)$,
        where
        $I_F\in\cL^{a}(F,F)$ is the identity operator on $F$.  We will have
        occasions to use operators constructed in this way in the sequel.}
        \end{remark}

        We now introduce the notion of {\em reproducing kernel $(\cA,
        \cB)$-correspondence}.

\begin{definition} \label{D:RKC}  {\rm
        Let $\cA$ and $\cB$ be $C^*$-algebras.  By an
        $(\cA,\cB)${\em -reproducing kernel correspondence} on a set
        $\Omega $, we mean an $(\cA,\cB)$-correspondence $E$ whose
        elements are $\cB$-valued functions $f \colon (\omega,a) \mapsto
        f(\omega,a) \in
        \cB$  on $\Omega \times \cA$, which is a vector space with respect to
        the usual point-wise vector-space operations and such that for each $\omega\in
        \Omega$ there is a {\em kernel element} $k_{\omega} \in E$
        with
        \begin{equation}  \label{repro-prop}
      f(\omega,a) = \langle a \cdot f, \, k_{\omega} \rangle_{E}.
        \end{equation}
        When this is the case we say that the function $\K \colon \Omega
        \times \Omega \to \cL(\cA, \cB)$ given by
        \begin{equation} \label{RKC-cpkernel}
     \K(\omega, \omega')[a] = k_{\omega'}(\omega,a)
        \end{equation}
        is the {\em reproducing kernel} for the reproducing kernel
        correspondence $E$.}
        \end{definition}

{}From the inner product characterization in (\ref{repro-prop}) of the
point evaluation for elements in an $(\cA,\cB)$-reproducing kernel
correspondence $E$ on $\Omega$ one easily deduces that the left $\cA$-action
and the right $\cB$-action are given by
        \begin{equation} \label{RKC-Aaction}
       (a \cdot f)(\omega', a') = f(\omega', a' a) \text{ and }  (f \cdot
       b)(\omega',a')=f(\omega',a')b.
       \end{equation}

        It is implicit in Definition \ref{D:RKC} that the map $a \mapsto
        k_{\omega'}(\omega,a) \in \cB $ is linear in $a \in \cA$ for each
        $\omega, \omega' \in
        \Omega$.  In fact the mapping from $\cA$ to $\cB$ given by $a \mapsto
        f(\omega,a)$ is $\cA$-linear for each fixed $f \in E$ and $\omega \in
        \Omega$.  If $\cA$ has a
        unit $1_{\cA}$, this follows from
        the general identity $f(\omega,a) = (a \cdot f)(\omega,1_\cA)$ (a
        consequence of \eqref{RKC-Aaction} together with the linearity of
        the point-evaluation map $f \mapsto f(\omega,1_\cA)$ from
        $E$ to $\cB$ for each fixed $\omega \in \Omega$ which in turn is an easy
        consequence of \eqref{repro-prop}). The general case follows by
        adapting this argument to the setting where one has only an
        approximate identity. Note
        also that we recover the element $k_{\omega'}$ from $\K$ by
using formula
        \eqref{RKC-cpkernel} to define $k_{\omega'}$ as a function of
        $(\omega,a)$ for
        each $\omega'\in \Omega$.

        The next proposition gives some elementary observations concerning
        the structure of reproducing kernel correspondences.

        \begin{proposition}  \label{P:RKC}
        If $E$ is a reproducing kernel $(\cA, \cB)$-correspondence with
        kernel elements $k_{\omega}$ for $\omega \in \Omega$, then the bounded
        evaluation map $e_{\omega,a}$
     from $E$ to $\cB$ given by $e_{\omega,a} \colon f \mapsto
        f(\omega,a)$ is adjointable
        for each fixed $(\omega,a) \in \Omega \times \cA$ and we
        have
        \begin{equation}  \label{kernel/pteval}
      a^* k_{\omega}b = e_{\omega,a}^{*}b \text{ for each  $\omega \in
      \Omega$, $a\in\cA$, and  $b\in\cB$.}
        \end{equation}
        Conversely, suppose that $E$ is an $(\cA,\cB)$-correspondence of
        $\cB$-valued functions on the set $\Omega \times \cA$
        satisfying \eqref{RKC-Aaction} and such that the evaluation map
        $$
        e_{\omega,a} \colon f \mapsto f(\omega,a)
        $$
        is a bounded and adjointable map from $E$ to $\cB$ for each $\omega
        \in \Omega$ and $a \in \cA$.
        Then $E$ is a reproducing kernel $(\cA, \cB)$-correspondence with
        reproducing kernel elements determined by \eqref{kernel/pteval}.

        Moreover, in either case, for each fixed $(\omega,a)$ the
point-evaluation
        map $e_{\omega,a} \colon E \to \cB$ is a $\cB$-module map:
        $$
      (f \cdot b)(\omega,a) = f(\omega,a) b \text{ for all } b \in \cB.
        $$
        \end{proposition}

        \begin{proof}
         Suppose $E$ is a reproducing kernel $(\cA,
\cB)$-correspondence with
        kernel elements $k_{\omega}$ for $\omega \in \Omega$.
        If $e_{\omega,a}$ denotes the evaluation map from $E$ to $\cB$
        given by $e_{\omega,a} \colon f \mapsto f(\omega,a)$, we have
        $$
        \langle e_{\omega,a}f, b \rangle_{\cB} = b^* f(\omega,a) =
b^*\langle a\cdot
        f, k_{\omega} \rangle_{E} =\langle f, a^*k_{\omega} \cdot b \rangle.
        $$
        So $e_{\omega,a}$ is adjointable with $e_{\omega,a}^{*} b=
a^*k_\omega b $ for
        any $b \in \cB$.

        On the other hand, if the evaluation map
        $$
        e_{\omega,a} \colon f \mapsto f(\omega,a)
        $$
        is a bounded and adjointable map from $E$ to $\cB$ for each $\omega
        \in \Omega$
        and $a \in \cA$, then there exists an
        $e_{\omega,a}^{*}$ so that
        $$ b^{*} (e_{\omega,a}f) =
        \langle e_{\omega,a}f, b \rangle_{\cB} = \langle f, e_{\omega,a}^{*}b
        \rangle_{E}.
        $$
        If $\cA$ and $\cB$ have identities $1_{\cA}$ and $1_{\cB}$
        respectively, we set
        $ k_{\omega} = e_{\omega,1_{\cA}}^{*}(1_{\cB})$.
        Using the first identity in (\ref{RKC-Aaction}) it follows from a computation
        similar to that above, that $a^* k_{\omega} =
        e_{\omega, a}^{*}(1_{\cB})$.
        We readily see that
     $$f(\omega,a)= e_{\omega,a}f =\langle f, a^*k_{\omega}\rangle_{E}=
     \langle a\cdot f, k_{\omega}\rangle_{E}.$$

        If $\cA$ and/or $\cB$ does not have a unit, one can do an approximate
        version of
        the above argument using an approximate identity for $\cA$ and/or
        $\cB$.  In any
        case, it follows that $E$ is a reproducing kernel $(\cA,
        \cB)$-correspondence with
        reproducing kernel elements determined by \eqref{kernel/pteval}.

        The last part follows from the definition of the right $\cB$-action
        given by \eqref{RKC-Aaction}.
        \end{proof}

        Given a reproducing kernel $(\cA, \cB)$-correspondence as in
        Definition \ref{D:RKC}, one can show that the associated
        reproducing kernel function $\K \colon \Omega \times \Omega \to
        \cL(\cA, \cB)$ defined by \eqref{RKC-cpkernel} is a completely
        positive kernel in the sense of \cite{BBLS}, i.e., the function
        $$
        ((\omega,a),(\omega',a'))\to \K(\omega, \omega')[a^{*}a' ]
        $$
        is a positive kernel in the classical sense of Aronszajn \cite{aron}
        (extended to the $C^*$-algebra-valued case), that is,
        {\em $ \sum_{i,j=1}^{N} b_i^*\K(\omega_{i}, \omega_{j})[a_{i}^{*} a_{j}]b_j$
        is a positive element of $\cB$
        for each choice of finitely many $(\omega_{1},a_{1}), \dots,
        (\omega_{N},a_{N})$ in $\Omega \times \cA$ and $b_1,\ldots,b_N$ in $\cB$.} In fact, by the axioms of
        an $(\cA, \cB)$-correspondence combined with the
        reproducing property of the kernel elements $k_{\omega}$, we have
        \begin{align*}
        \sum_{i,j=1}^{N}b_i^*\K(\omega_i, \omega_j)[{a_i}^*a_j]b_j  & =
        \sum_{i,j=1}^{N} b_i^*\langle {a_i}^*a_j k_{\omega_j},k_{\omega_i}
\rangle_{E}b_j\\
        & = \sum_{i,j=1}^{N} \langle {a_j} k_{\omega_j}b_j,
        a_i k_{\omega_i}b_i\rangle_{E} \\
        & =  \left\langle \sum_{j=1}^{N}
        {a_j} k_{\omega_j}b_j, \sum_{i=1}^{N} a_i k_{\omega_i}b_i \right\rangle_{E}
        \\& \ge 0.
        \end{align*}
        Actually, we have the following equivalent statements.

        \begin{theorem}
        Given a function $\K \colon \Omega \times \Omega \to
        \cL(\cA,\cB)$, the following are equivalent:
        \begin{enumerate}
        \item $\K$ is a completely positive kernel in the sense that the
        function from $(\Omega \times \cA) \times (\Omega \times \cA) \to
        \cL(\cA, \cB)$ given by
        $$
     ((\omega,a),(\omega',a')) \mapsto \K(\omega', \omega)[a^{*}a']
        $$
        is a positive kernel in the sense of Aronszajn:
        $$
        \sum_{i,j=1}^{N}b_i^*\K(\omega_i, \omega_j)[{a_i}^*a_j]b_j \ge 0 \text{ in } \cB
        \text{ for all } (\omega_{1}, a_{1}), \dots, (\omega_{N}, a_{N}) \in
        \Omega \times \cA\text{ and }b_1,\ldots b_N\in\cB.
        $$

        \item $\K$ has a Kolmogorov decomposition in the sense of
        \cite{BBLS}, i.e., there exists an
        $(\cA,\cB)$-correspondence $E$ and a mapping $\omega \mapsto
k_{\omega} $
        from $\Omega$ into $E$ such
        that
        \begin{equation*} %\label{Kolmogorov}
     \K(\omega',\omega)[a] =\langle  a \cdot k_{\omega}, \,
k_{\omega'}\rangle_{E}
     \text{ for all } a\in \cA.
     \end{equation*}

        \item $\K$ is the reproducing kernel for an $(\cA,\cB)$-reproducing
        kernel correspondence $E = E(\K)$, i.e., there is an $(\cA,
        \cB)$-correspondence $E = E(\K)$ whose elements are $\cB$-valued
        functions on $\Omega \times \cA$
        such that the function $k_{\omega} \colon (\omega',a') \mapsto
        \K(\omega', \omega)[a']$ is in $E(\cK)$ for each $\omega \in \Omega$
        and has the reproducing property
        \begin{equation*}  %\label{cor-criterion3}
     \langle a \cdot f, k_{\omega} \rangle_{E(\K)} =
     \langle f, a^{*} \cdot k_{\omega} \rangle_{E(\K)} =
     f(\omega,a)
     \text{ for all } \omega \in \Omega \text{ and } a \in \cA
        \end{equation*}
        where $a^{*} \cdot k_{\omega}$ is given by
        \begin{equation}  \label{concrete-Kol}
        (a^{*} \cdot k_{\omega})(\omega', a') = \K(\omega', \omega)[a'a^{*}].
        \end{equation}
        \end{enumerate}
        \end{theorem}

        \begin{proof}
        For the equivalence of  (1) and (2), we refer to Theorem 3.2.3 in
        \cite{BBLS}. The argument in the paragraph preceding the statement
        of the theorem shows that (3) $\Longrightarrow$ (1).  To see that
        (2) $\Longrightarrow$ (3), assume that $E$ is an $(\cA,
        \cB)$-correspondence as in (2).  Without loss of generality we may
        assume that
        \begin{equation}  \label{nondegen}
      E = \overline{\operatorname{span}} \{ a \cdot k_{\omega} b \colon a
      \in \cA, \omega \in \Omega, b \in \cB\}.
        \end{equation}
        We view elements $f$ of $E$ as
        $\cB$-valued functions on $\Omega \times \cA$ by defining
        \begin{equation*}  %\label{point-eval}
        f(\omega,a) = \langle a \cdot  f, k_{\omega} \rangle_{E} \text{ for
        each  $\omega \in \Omega$ and $a \in \cA$.}
        \end{equation*}
        The nondegeneracy assumption \eqref{nondegen} says that
        $$
      f(\omega,a)b = 0 \text{ for all } a \in \cA,\, b \in \cB, \,
     \omega \in \Omega \Longrightarrow f = 0 \text{ in } E.
        $$
        Hence the map $f \mapsto f(\cdot, \cdot )$ is injective.
        Finally \eqref{concrete-Kol} holds by definition.
        \end{proof}

        We now tailor this general theorem to the case where $\cB= \cL(\cE)$
        for a Hilbert space $\cE$. Note that $\cE$ is a
        $(\cL(\cE),\C)$-correspondence, i.e., a Hilbert space with a
        $*$-representation $b \mapsto \varphi(b) \in {\mathcal L}(\cE)$ of
        $\cL(\cE)$ (namely, the identity representation). Hence, given that
        $E$ is an
        $(\cA, \cL(\cE))$-correspondence, we may form
        the tensor product $ E \otimes_{\cL(\cE)} \cE$ to obtain an
        $(\cA, \C)$-correspondence, i.e., a Hilbert space which we will
        denote by $\cH$  equipped with an
        $\cL(\cH)$-valued $*$-representation $\pi
        \colon \cA \to \cL(\cH)$ of $\cA$.
     Similarly, if we view $\cB = \cL(\cE)$ as a
$(\cL(\cE),\cL(\cE))$-correspondence, we
        may form the tensor product $ \cL(\cE) \otimes \cE$ to arrive at
        the Hilbert space $\cE$, via the balancing
        (\ref{balance}), viewed as a $(\cL(\cE), {\mathbb
        C})$-correspondence.  Let us suppose also that $E$ is a reproducing
        kernel correspondence.  Then via the formula $f \otimes e \in
        E\otimes \cE \mapsto f(\omega,a) \otimes e \in \cL(\cE)\otimes \cE$ for
        each $\omega \in \Omega$ and $a\in\cA$ extended via linearity
and continuity to
        the whole space $E \otimes \cE$, we may view each
        $\mathbf f \in \cH = E\otimes \cE$ as a $\cE$-valued function on
        $\Omega \times \cA$ such that point-evaluation $\mathbf f \mapsto
        \mathbf f(\omega,a)$ is continuous, i.e., $\cH$ is a reproducing kernel
        Hilbert space of vector-valued functions on $\Omega
        \times \cA$, but with the additional wrinkle that there is also a
        representation $a \mapsto \pi (a)$ of $\cA$ on $\cH$ with
        $\pi(a)(f \otimes e) = ( a \cdot f) \otimes e$ such that
        $$
        ( \pi(a)  ( f \otimes e))(\omega',a') =  f
        (\omega', a' a) \otimes e
        $$
        with reproducing kernel (in the sense of a vector-valued reproducing
        kernel Hilbert space) $K(\cdot, \cdot)$ of the
        special form
        $$
     K((\omega',a'), (\omega,a) ) = \K(\omega',\omega)[a^{*}a']
        $$
        for a completely positive kernel $\K \colon \Omega
        \times \Omega \to \cL(\cA, \cL(\cE))$: for ${\mathbf f} \in \cH(\K)$,
        $e \in \cE$  and $(\omega, a) \in \Omega \times \cA$,
        $$\langle \mathbf f , \K (\cdot, \omega)[a]e \rangle_{\cH} =
        \langle \mathbf f(\omega,a), e \rangle_{\cE}
        $$
        where  $\K$  is completely positive.
        This leads us to an
        alternative reproducing-kernel interpretation of a completely
        positive kernel $\K \colon \Omega \times \Omega \to
        \cL(\cA, \cB)$ for the case where $\cB = \cL(\cE)$ for a Hilbert
        space $\cE$.

        \begin{theorem}  \label{T:RKHilcor}
        Suppose that $\cA$ is a $C^{*}$-algebra, $\cE$ is a Hilbert space
        and that a function $\K \colon \Omega \times \Omega \to \cL(\cA,
        \cL(\cE))$ is given.
        The following conditions are equivalent.
        \begin{enumerate}
        \item The function $\K$ is a completely positive
        kernel in the sense that
        $$ \sum_{i,j=1}^{N} \langle \K(\omega_{i},
\omega_{j})[a_{i}^{*}a_{j}]e_j,e_i \rangle\ge 0
        $$
        for all finite collections $\omega_{1}, \dots, \omega_{N} \in
        \Omega$, $a_{1}, \dots, a_{N} \in \cA$ and $e_1,\ldots,e_N\in\cE$ for $N=1,2, \dots$.

        \item The kernel ${\mathbb K}$ has a Kolmogorov decomposition:
        there is a Hilbert
        space $\cH$ together with a $*$-representation $\pi\colon \cA
        \to \cL(\cH)$ of $\cA$ and a mapping $H \colon \Omega \to
        \cL(\cH, \cE)$ so that
        $$
        \K(\omega',\omega)[a] =  H(\omega') \pi(a) H(\omega)^{*}.
        $$

        \item There is a ($\cA, {\mathbb C})$-correspondence,
         i.e., a Hilbert space $\cH = \cH(\K)$
         together with a $*$-representation $a \mapsto \pi(a) \in
        \cL(\cH)$ of $\cA$
         whose elements are
        $\cE$-valued functions on $\Omega \times \cA$
         such that:
        \begin{enumerate}
       \item The $*$-representation $\pi$ is given by
         $$
        ( \pi(a) {\mathbf f})(\omega',a')
         = {\mathbf f}(\omega',a'a).
         $$
         \item The function $ k_{\omega}  \colon \Omega
        \times \cA \to \cL(\cE)$ given by
        $$
        k_{\omega}(\omega',a') \colon e \mapsto
        \K(\omega', \omega)[a'] e
        $$
        is such that $k_{\omega} e \in \cH(\K)$ for each $\omega \in
        \Omega$ and $e \in \cE$ and has the reproducing kernel property:
      $$
      \langle {\mathbf f}, \pi(a)^{*}k_{\omega} e \rangle_{\cH(\K)} =
        \langle {\mathbf f}(\omega,a), e \rangle_{\cE}.
      $$
      \end{enumerate}
      \end{enumerate}
      \end{theorem}

      Let us say that the object described in part (3) of Theorem
      \ref{T:RKHilcor} a {\em reproducing kernel Hilbert correspondence}
      (over the $C^{*}$-algebra $\cA$ with values in the coefficient
      space $\cE$).

      \begin{remark} {\rm  \label{R:pos-ker}
          If $\cH(\K)$ is a reproducing kernel Hilbert correspondence
          space  as
          in part (3) of Theorem \ref{T:RKHilcor}, a special situation
          occurs if the coefficient space $\cE$ is also equipped with a
          $*$-representation $\pi_{\cE} \colon \cA \to \cL(\cE)$.  In
          this case it may or may not happen that point evaluation is an
          $\cA$-module map, i.e., that
       \begin{equation}  \label{special}
       (a \cdot f)(\omega',a') = a \cdot f(\omega',a') \text{ or
equivalently }
       (\pi(a) f)(\omega',a') =
       \pi_{\cE}(a) f(\omega',a').
       \end{equation}
       When \eqref{special} does occur and if also $\cA$ has a unit
       $1_{\cA}$, one can show that the associated
       completely positive kernel $\K(\omega, \omega')[a]$ has
the special property
       \begin{equation}  \label{kerprop0}
        \K(\omega, \omega')[a^{*}a'] = \pi_{\cE}(a)^{*} \K(\omega,
        \omega')[1_{\cA}] \pi_{\cE}(a')
       \end{equation}
       and hence complete positivity of $\K$ reduces to standard
       (Aronszajn) positivity for the kernel
       $K_{0} \colon \Omega \times \Omega \to \cL(\cE)$ given
       by
       $$
       K_{0}(\omega, \omega') = \K(\omega,\omega')[1_{\cA}].
       $$
       Indeed, the computation
       \begin{align*}
           \langle \K(\omega,\omega')[a^{*}a'] e', e \rangle_{\cE} & =
           \langle a^{*}a' \cdot k_{\omega'}e',\, k_{\omega}e
\rangle_{\cH(\K)} \\
           & = \langle (a^{*}a' \cdot k_{\omega'}e')(\omega), \,
e \rangle_{\cE} \\
           & = \langle a^{*} a' \cdot (k_{\omega'}e')(\omega), \, e \rangle_{\cE}
           \text{ (by assumption \eqref{special})} \\
           & = \langle a' (k_{\omega'}e')(\omega), a e \rangle_{\cE} \\
           & = \langle (a' \cdot k_{\omega'} e')(\omega), a e
\rangle_{\cE} \text{
           (by \eqref{special} again)} \\
           & = \langle a' \cdot k_{\omega'}e', k_{\omega} (ae)
\rangle_{\cH(\K)} \\
           & = \langle \K(\omega, \omega')[a']e', ae \rangle_{\cE} \\
           & = \langle a^{*} \K(\omega, \omega') [a'] e', e \rangle_{\cE}
        \end{align*}
        shows that
        \begin{equation}  \label{kerprop1}
       \K(\omega,\omega')[a^{*}a'] = a^{*} \K(\omega, \omega')[a'].
       \end{equation}
        On the
        other hand, the positive-kernel property of the kernel
        $$ ((\omega,a), (\omega',a')) \mapsto K((\omega,a), (\omega',a')):
        =  \K(\omega, \omega')[a^{*}a']
        $$
        implies that $K$ is Hermitian, i.e., $K((\omega,a), (\omega',a')) =
        K((\omega',a'),(\omega,a))^{*}$,
        i.e.,
        $$
        \K(\omega, \omega')[a^{*}a'] = \left( \K(\omega',
\omega)[a^{\prime *} a ]
        \right)^{*}.
        $$
        In particular,
        \begin{align*}
        \K(\omega, \omega')[a'] & = \left( \K(\omega', \omega)[a^{\prime
        *}] \right)^{*} \\
        & = \left( a^{\prime *} \K(\omega',\omega) [1_{\cA}]
\right)^{*} \text{ (by
        \eqref{kerprop1})} \\
        & = \K(\omega, \omega')[1_{\cA}] a'
       \end{align*}
       and hence also
       \begin{equation}  \label{kerprop2}
           \K(\omega, \omega')[a'] = \K(\omega, \omega')[1_{\cA}] a'.
       \end{equation}
       Combining \eqref{kerprop1} and \eqref{kerprop2} gives
       \eqref{kerprop0} as claimed.}
       \end{remark}

%%%%%%%%%%%%%%%%%%%%%%%%%%%%%%%%%%%%%%%%%%%%%%%%%%%%%%%%%%%%%%%%%%%%%%%%

%%%%%%%%%%%%%%%%%%%%%%%%%%%%%%%%%%%%%%%%%%%%%%%%%%%%%%%%%%%%%%%%%%%%%%%%
       \section{Function-theoretic operator theory
       associated with a correspondence $E$}\label{S:Hardyspaces}

       \setcounter{equation}{0}
       In this section we obtain the analogues of Hardy spaces, Toeplitz
       operators, $Z$-transform and Schur class attached to a
       $\cA$-$W^{*}$-correspondence $E$ together with a $*$-representation
       $\sigma$ of $\cA$.  These results flesh out more fully the
       function-theoretic aspects of the work of Muhly-Solel \cite{MS-JFA,
       MS-Annalen, MS-Schur}.

%%%%%%%%%%%%%%%%%%%%%%%%%%%%%%%%%%%%%%%%%%%%%%%%%%%%%%%%%%%%%%%%%%%%%%%%
       \subsection{Hardy Hilbert spaces associated with a
correspondence $E$}

       In this section we shall consider the situation where $\cA
= \cB$; we
       abbreviate the term $(\cA, \cA)$-correspondence to simply
       $\cA$-correspondence.
       We also now restrict our attention to the case where $\cA$ is a
       von Neumann algebra and let $E$ be a $\cA$-$W^{*}$-correspondence.
       This means that $E$ is a
       $\cA$-correspondence which is also {\em
       self-dual} in the sense that any right $\cA$-module map $\rho \colon
       E \to \cA$ is given by taking the inner product against some element
       $e_{\rho}$ of $E$:  $\rho(e) = \langle e, e_{\rho} \rangle_{E} \in
       \cA$. Moreover, the space $\cL^{a}(E)$ of adjointable operators on
       the $W^{*}$-correspondence $E$ is in fact a
$W^{*}$-algebra, i.e., is the
       abstract version of a von Neumann algebra with an
ultra-weak topology
       (see \cite{MS-Annalen}).

       Since $E$ is a $\cA$-correspondence, we may use
       Definition \ref{D:tensor} to define the self-tensor
       product $E^{\otimes 2} = E \otimes_{\cA} E$ which is again an
       $\cA$-correspondence,
       and, inductively, an $\cA$-correspondence
       $E^{\otimes n} = E \otimes_\cA (E^{\otimes (n-1)})$ for each $n =
       1,2, \dots$.  If we use $a \mapsto \varphi(a)$ to denote the left
       $\cA$-action $\varphi(a) e = a \cdot e$ on $E$, we denote the left
       $\cA$-action on $E^{\otimes n}$ by $\varphi^{(n)}$:
       $$
         \varphi^{(n)}(a) \colon \xi_{n} \otimes \xi_{n-1} \otimes \cdots
         \otimes \xi_{1}\mapsto (\varphi(a) \xi_{n}) \otimes
         \xi_{n-1} \otimes \cdots \otimes \xi_{1}.
       $$
       Note that, using the notation in
       \eqref{tensoroperator}, we may write $\varphi^{(n)}(a)=\varphi(a)\otimes
       I_{E^{\otimes n-1}}$. We formally set $E^{\otimes 0} = \cA$. Then
       the Fock space $\cF^{2}(E)$ is defined to be
       \begin{equation}  \label{Fock}
          \cF^{2}(E) = \oplus_{n=0}^{\infty} E^{\otimes n}
       \end{equation}
       and is also an $\cA$-correspondence.
       We denote the left
       $\cA$-action on $\cF(E)$ by $\varphi_{\infty}$:
       \begin{equation}  \label{ArepFock}
         \varphi_{\infty}(a) \colon \oplus_{n=0}^{\infty} \xi^{(n)} \mapsto
         \oplus_{n=0}^{\infty}( \varphi^{(n)}(a) \xi^{(n)}) \text{ for }
         \oplus_{n=0}^{\infty} \xi^{(n)} \in
         \bigoplus_{n=0}^{\infty} E^{\otimes n},
       \end{equation}
       or, more succinctly,
       $$
\varphi_\infty(a)=\textup{diag}(a,\varphi^{(1)}(a),\varphi^{(2)}(a),\ldots).
       $$

       In addition to the von Neumann algebra $\cA$ and the
$\cA$-correspondence $E$,
       suppose that we are also given an auxiliary Hilbert space
$\cE$ and a
       nondegenerate $*$-homomorphism $\sigma \colon \cA \to \cL(\cE)$; as
       this will be the setting for much of the analysis to
follow, we refer
       to such a pair $(E, \sigma)$ as a {\em correspondence-representation
       pair}.
       Then the Hilbert space $\cE$ equipped with $\sigma$ becomes an
       $(\cA,\C)$-correspondence with left $\cA$-action given by $\sigma$:
       $$
          a \cdot y = \sigma(a) y \text{ for all } a \in \cA \text{ and } y
          \in \cE.
       $$
       We let $E \otimes_{\sigma} \cE$ be the
       associated tensor-product $(\cA, {\mathbb C})$-correspondence $E
       \otimes_{\cA} \cE$
       as in Definition \ref{D:tensor}.  As $\cF^{2}(E)$ is also an
       $\cA$-correspondence,
       we may also form the $(\cA, {\mathbb C})$-correspondence
       $$
       \cF^{2}(E, \sigma) := \cF^{2}(E) \otimes_{\sigma} \cE =
       \bigoplus_{n=0}^{\infty}
       ( E^{\otimes n} \otimes_{\sigma} \cE),
       $$
       with left $\cA$-action given by the $*$-representation
       \[
       \varphi_{\infty,\sigma}(a)=\varphi_\infty(a)\otimes I_\cE.
       \]

       It turns out that
       $\cF^{2}(E, \sigma)$ is also a $(\sigma(\cA)', {\mathbb
       C})$-correspondence, where
       $\sigma(\cA)'\subset\cL(\cE)$ denotes the commutant of
       $\sigma(\cA)$:
       \begin{equation}\label{sigmaA'0}
       \sigma(\cA)' = \{ b \in \cL(\cE) \colon b \sigma(a) = \sigma(a) b
       \text{ for all }  a \in \cA\},
       \end{equation}
       and the left $\sigma(\cA)'$-action is given by the
$*$-representation
       $\iota_{\infty,\sigma}$ of $\sigma(\cA)'$ on $\cL(\cF^2(E,\sigma))$:
       \begin{equation}\label{iota}
       \iota_{\infty,\sigma}(b)=I_{\cF^2(E)}\otimes b\mbox{ for each }
       b\in\sigma(\cA)',
       \end{equation}
       using the notation in (\ref{tensoroperator}).  Note that
       $b\in\cL(\cE)$ is in $\sigma(\cA)'$ precisely when $b$ is an
       $\cA$-module map, so that $I_{\cF^2(E)}\otimes b$ is a well defined
       operator on $\cF^{2}(E, \sigma)$.
       Moreover, $\varphi_{\infty,\sigma}(a)$ commutes with
       $\iota_{\infty,\sigma}(b)$ for each $a\in\cA$ and $b\in\sigma(\cA)'$
       since
       \[
       \varphi_{\infty,\sigma}(a)\iota_{\infty,\sigma}(b)
       =\varphi_{\infty}(a)\otimes b
       =\iota_{\infty,\sigma}(b)\varphi_{\infty,\sigma}(a).
       \]
       Thus $\iota_{\infty,\sigma}(b)$ is a $\cA$-module
       map for each $b\in\sigma(\cA)'$ and $\varphi_{\infty,\sigma}(a)$
       is a $\sigma(\cA)'$-module map for each $a\in\cA$.

       We denote by $E^{\sigma}$ the set of all
       bounded linear operators $\mu\colon \cE\to E \otimes_{\sigma} \cE$
       which are
       also $\cA$-module maps:
       \begin{equation} \label{Esigma0}
       E^{\sigma} = \{ \mu\colon \cE\to E \otimes_{\sigma} \cE
       \colon
           \mu \sigma(a)  = ( \varphi(a) \otimes I_{\cE}) \mu \},
       \end{equation}
       and $(E^{\sigma})^*$ for the set of adjoints
       (which are also $\cA$-module maps):
       \begin{equation} \label{Esigma*}
       (E^{\sigma})^*=\{\eta:E \otimes_{\sigma} \cE\to\cE
       \colon \eta^*\in E^{\sigma}\}.
       \end{equation}

       More generally, for a given $\eta \in (E^{\sigma})^{*}$, we may
       define operators
       $\eta^{n} \colon E^{\otimes n} \otimes_{\sigma} \cE \to \cE$
       ({\em generalized powers}) by
       \begin{equation*}  %\label{genpower}
        \eta^{n} = \eta (I_{E} \otimes \eta) \cdots
(I_{E^{\otimes n-1}}
       \otimes \eta)
       \end{equation*}
       where we use the identification
       $$
          E^{\otimes n} \otimes_{\sigma} \cE = E^{\otimes n-1}
          \otimes_\cA (E \otimes_{\sigma} \cE)
       $$
       in these definitions. We also set
       $\eta^0=I_\cE\in\cL(\cE)$. Again the fact that $\eta$
       is an $\cA$-module map ensures that $I_{E^{\otimes k}}\otimes\eta$
       is a well defined operator in $\cL(E^{\otimes k+1}\otimes_\sigma
       \cE, E^{\otimes k}\otimes_\sigma \cE)$. The defining $\cA$-module
       property of $\eta$ in \eqref{Esigma*} then extends to the
       generalized powers $\eta^{n}$ in the form
       \begin{equation}  \label{mod-prop-ext}
        \eta^{n} (\varphi^{(n)}(a) \otimes I_{\cE})
        = \sigma(a) \eta^{n},
       \end{equation}
       i.e., $\eta^{n}$ is also an $\cA$-module map.

       Denote by ${\mathbb D}((E^{\sigma})^{*})$ the set of strictly
       contractive elements of $(E^{\sigma})^{*}$:
       $$
         {\mathbb D}((E^{\sigma})^{*}) = \{ \eta \in
(E^{\sigma})^{*} \colon
         \| \eta \| < 1 \}.
       $$
       Then, for $\eta \in {\mathbb D}((E^{\sigma})^{*})$ and $b \in
       \sigma(\cA)'$, we may define a bounded operator $ f \mapsto
       f^{\wedge}(\eta, b)$
       from $\cF^2(E,\sigma)$ into $\cE$ by
       \begin{equation}  \label{pteval}
        f^{\wedge}(\eta,b) = \sum_{n=0}^{\infty}
        \eta^{n}(\iota_{\infty,\sigma}(b) f)_{n} =
        \sum_{n=0}^{\infty} \eta^{n}(I_{E^{\otimes n}}
        \otimes b) f_{n}  \text{ if } f = \oplus_{n=0}^{\infty} f_{n}.
       \end{equation}

       Note that the fact that $\| \eta\| < 1$ guarantees that
the series in
       \eqref{pteval} converges.
       The $\cA$-module properties of $\iota_{\infty,\sigma}(b)$ and
       each generalized power $\eta^n$ (see \eqref{mod-prop-ext}) for given
       $b \in \sigma(\cA)'$ and $\eta \in {\mathbb D}((E^{\sigma})^{*})$
       imply that the point-evaluation
       $f \mapsto f^{\wedge}(\eta,b)$ is also an
       $\cA$-module map:
       \begin{equation*}  %\label{mod-map}
        (\varphi_{\infty,\sigma}(a) f)^{\wedge}(\eta,b) =
        \sigma(a)  f^{\wedge}(\eta,b).
         \end{equation*}
         However, the point-evaluation $f \mapsto f^{\wedge}(\eta, b)$ is
         not a $\sigma(\cA)'$-module map, i.e., there is no
guarantee for the
         general validity of the identity $(\iota_{\infty,\sigma}(b)
         f)^{\wedge}(\eta',b') = b f^{\wedge}(\eta', b')$, but
rather we have
         the property
         $$
         (\iota_{\infty,\sigma}(b)) f)^{\wedge}(\eta',b') =
f^{\wedge}(\eta',
         b'b).
         $$

       We denote the space of all $\cE$-valued functions
       on ${\mathbb D}((E^{\sigma})^{*}) \times \sigma(\cA)'$
         of the form $(\eta,b) \mapsto
         f^{\wedge}(\eta,b)$ for some $f \in \cF^{2}(E, \sigma)$
by $H^{2}(E, \sigma)$
       with norm $\|f^{\wedge} \|_{H^{2}(E,\sigma)}$ chosen so as to make
       the map $  f \mapsto
         f^{\wedge}$ a coisometry from $\cF^{2}(E, \sigma)$ to
       $H^{2}(E, \sigma)$:
         $$
         H^{2}(E, \sigma) = \{ f^{\wedge} \colon f \in \cF^{2}(E,
\sigma) \}
         \text{ with } \|f^{\wedge}\|_{H^{2}(E, \sigma)} = \|
P_{(\operatorname{Ker}
         \Phi )^{\perp}} f \|_{\cF^{2}(E, \sigma)}
       $$
         where we denote by $\Phi$ (the generalized {\em Fourier} or
       $Z$-transform for this setting) the transformation from $\cF^{2}(e,
       \sigma)$ into $H^{2}(E, \sigma)$ given by
       \begin{equation}  \label{defPhi}
          \Phi \colon f \mapsto f^{\wedge}.
       \end{equation}
       Then we have the following result.

       \begin{theorem}  \label{T:H2Esigma}
        The space $H^{2}(E,\sigma)$ is a reproducing kernel
           Hilbert correspondence  $\widetilde \cH(\K)$ (as in part (3) of
           Theorem \ref{T:RKHilcor})   over $\sigma(\cA)'$
           consisting of
           $\cE$-valued functions on ${\mathbb D}((E^{\sigma})^{*}) \times
           \sigma(\cA)'$
           with the $*$-representation of
           $\sigma(\cA)'$ given by
           \begin{equation} \label{action}
         (b \cdot f^{\wedge})(\eta',b') = (\iota_{\infty, \sigma}(b)
         f)^{\wedge}(\eta',b') \text{ for } b \in \sigma(A)'.
           \end{equation}
           The completely positive kernel $\K$ associated with
$H^{2}(E,\sigma)$
           as in Theorem \ref{T:RKHilcor}
           $$ \K_{E,\sigma}
        \colon {\mathbb D}((E^{\sigma})^{*}) \times {\mathbb
        D}((E^{\sigma})^{*})
           \to \cL( \sigma(\cA)', \cL(\cE))
           $$
           is the {\em Szeg\"o kernel} for our setting given by
          \begin{equation}  \label{E-sigma-ker}
          \K_{E,\sigma}(\eta, \zeta)[b] =
          \sum_{n=0}^{\infty} \eta^{n} (I_{E^{\otimes n}} \otimes  b)
           (\zeta^{ n})^{*} \text{ for } b \in \sigma(\cA)'.
         \end{equation}
        \end{theorem}

        \begin{proof}

        Define $\Phi \colon \cF^{2}(E, \sigma)
        \to H^{2}(E, \sigma)$ as in \eqref{defPhi}.
        By the definition of the norm on $H^{2}(E, \sigma)$,
$\Phi$ is a
       coisometry.
        For each $b \in \sigma(\cA)'$ and $\eta \in {\mathbb
        D}((E^{\sigma})^{*})$,
        define an associated {\em controllability
operator}\footnote{The
        terminology is motivated by connections with system
theory; for a
        systematic account for the Drury-Arveson and free-semigroup
        algebra settings, we refer to \cite{BBF1}.}   ${\mathcal
        C}_{b, \eta} \colon  \cF^{2}(E, \sigma) \to \cE$ by
          $$
           {\mathcal C}_{b,\eta}  \colon f \mapsto
        f^{\wedge}(\eta, b) \text{ if } f\in\cF^{2}(E, \sigma).
         $$
         By definition,
       \begin{equation*}\label{kerPhi}
         \operatorname{Ker} \Phi  = \bigcap_{b \in \sigma(\cA)',
         \eta \in {\mathbb D}((E^{\sigma})^{*})}
         \operatorname{Ker}{\mathcal C}_{b, \eta}.
         \end{equation*}
         The initial space of the coisometry $\Phi$ is the
         orthogonal complement of its kernel, namely
         $$
          (\operatorname{Ker} \Phi)^{\perp} =
\overline{\operatorname{span}}
          \{
          \operatorname{Ran} {\mathcal C}_{b,\eta}^{*} \colon b \in
          \sigma(\cA)',\, \eta \in {\mathbb D}((E^{\sigma})^{*}) \},
         $$
         where the {\em observability operator} ${\mathcal
C}_{b,\eta}^{*}$ is
         given by
         $$
         {\mathcal C}_{b,\eta}^{*} \colon e \mapsto \oplus_{n=0}^{\infty}
          (I_{E^{\otimes n}} \otimes b^{*}) (\eta^{n})^{*} e \
          \in \cF^{2}(E) \otimes_{\sigma} \cE.
         $$
         We compute

         \begin{align*}
         \langle f^{\wedge}(\zeta,b), e \rangle_{\cE} &
         =
          \langle {\mathcal C}_{b,\zeta} f, e \rangle_{\cE} \\
         & = \langle  f, {\mathcal C}_{b, \zeta}^{*} e
           \rangle_{\cF^{2}(E, \sigma)} \\
         & = \langle  f^{\wedge}, \Phi ({\mathcal C}_{b, \zeta}^{*} e)
         \rangle_{H^{2}(E, \sigma)}\\
         & = \langle  f^{\wedge}, b^* \cdot \Phi ({\mathcal
C}_{I_\cE, \zeta}^{*}
         e) \rangle_{H^{2}(E, \sigma)},
          \end{align*}
         where we use the fact seen above that ${\mathcal
C}_{b,\zeta}^{*} e$ is in
         the initial space of $\Phi$ and that $\Phi((I_{\cF^{2}}(E) \otimes
         b)f) = b (\Phi f)$ for each $b \in \sigma(\cA)'$ and $f \in
         \cF^{2}(e, \sigma)$.
         Hence the operator
         \[
          k_{E,\sigma;\zeta}:=\Phi{\mathcal C}_{I_{\cE},\zeta}^{*}
          :\cE\to H^{2}(E, \sigma)
          \]
          has the reproducing property for $H^{2}(E, \sigma)$; see part
          (3.b) in Theorem \ref{T:RKHilcor}. Since
           $\Phi$ is a coisometry and
          $I_\cE\in\sigma(\cA)'$, we obtain that the reproducing
          kernel $\K_{E, \sigma}$ is necessarily given by
          \begin{align*}
        \K_{E,\sigma}(\eta, \zeta)[b] & = b \cdot
k_{E,\sigma;\zeta}(\eta)\\
        & = {\mathcal C}_{b,\eta}\Phi^*\Phi{\mathcal
C}_{I_\cE,\zeta}^{*} \\
        & = \sum_{n=0}^{\infty} \eta^{n}(I_{E^{\otimes n}} \otimes b)
        (\zeta^{n})^{*}
          \end{align*}
          in agreement with \eqref{E-sigma-ker}.
          \end{proof}

       %%%%%%%%%%%%%%%%%%

          From the proof of Theorem \ref{T:H2Esigma} we see that
we have the
          identification
       \begin{equation*} %\label{kSz}
           b^{*} \cdot k_{E,\sigma; \eta} e
           = {\mathcal C}_{b, \eta}^{*} e =
           \oplus_{n=0}^{\infty} (I_{E^{\otimes n}} \otimes
b^{*}) \eta^{n*} e
       \end{equation*}
       and the {\em initial space} for the coisometry $\Phi \colon \cF(E)
       \otimes_{\sigma} \cE \to H^{2}(E, \sigma)$ can be identified as
       \begin{equation}  \label{initial}
        [\cF(E) \otimes_{\sigma} \cE]_{\text{initial}} =
        \overline{\operatorname{span}} \{ b \cdot k_{E,\sigma; \eta} e
           \colon
        b \in \sigma(\cA)', \, \eta \in {\mathbb D}((E^{\sigma})^{*}),
        e\in \cE \}.
       \end{equation}

       %%%%%%%%%%%%%%%%%%%%%%%%%%%%%%%%%%%%%%%%%%%%%%%
       \subsection{Analytic Toeplitz algebras associated with a
       correspondence $E$}

       Given an $\cA-W^{*}$-correspon\-dence $E$, we let
$\cF^{2}(E)$ be the
       associated
       Fock space as in \eqref{Fock}.  We have already defined the
       $*$-representation of $\cA$ to $\cL^{a}(\cF^{2}(E))$ given by
       $a \mapsto \varphi_{\infty}(a)$ as in \eqref{ArepFock}.
       If we view operators on
       $\cF^2(E)$ as matrices induced by the decomposition
       $\cF^{2}(E) = \bigoplus_{n=0}^{\infty} E^{\otimes n}$ of
       $\cF^2(E)$, we see that each $\varphi_{\infty}(a)$
       has a diagonal representation $\varphi_{\infty}(a) =
       \operatorname{diag}_{n = 0,1,\dots} \varphi^{(n)}(a)$. In addition
       to the operators $\varphi_{\infty}(a) \in
       \cL^{a}(\cF^2(E))$, we introduce the so-called {\em
       creation operators} on $\cF^{2}(E)$ given, for each $\xi \in E$,
       by the subdiagonal (or shift) block matrix
       $$
          T_{\xi} = \begin{bmatrix} 0 & 0 & 0 & \cdots \\
T_{\xi}^{(0)} & 0 &
          0  & \cdots \\
          0 & T_{\xi}^{(1)} & 0 & \cdots & \\
          \ddots & \ddots & \ddots  &  \end{bmatrix}
       $$
       where the block entry $T_{\xi}^{(n)} \colon E^{\otimes n} \to
       E^{\otimes n+1}$ is given by
       \begin{equation*}
        T_{\xi}^{(n)} \colon \xi_{n}\otimes \cdots \otimes \xi_{
        1} \mapsto
         \xi \otimes \xi_{n} \otimes \cdots \otimes \xi_{1}.
       \end{equation*}
       The operator $T_\xi$ is also in $\cL^{a}(\cF^2(E))$.  In summary,
       both $T_\xi$ and $\varphi_{\infty}(a)$ are $\cA$-module maps
       with respect to the right $\cA$-action on $\cF^2(E)$ for each
       $\xi\in E$ and $a\in\cA$.  Moreover, one easily
       checks that
       \[
\varphi_{\infty}(a)T_\xi=T_{a\xi}=T_{\varphi(a)\xi}\quad\mbox{and}\quad
       T_\xi\varphi_{\infty}(a)=T_{\xi a}\quad\mbox{for each
}a\in\cA\mbox{ and }
       \xi\in E.
       \]

       We let $\cF^{\infty}(E)$ denote the
       weak-$*$ closed algebra generated by the collection of operators
       $$
         \{ \varphi_{\infty}(a), \, T_{\xi} \colon a \in \cA
\text{ and } \xi
         \in E \}
       $$
       in the $W^{*}$-algebra $\cL^{a}(\cF(E))$---we prefer this notation
       over the notation $H^{\infty}(E)$ used for this object in
       \cite{MS-JFA, MS-Annalen}.

       Suppose now that we are also given a
       $*$-representation $\sigma$ of $\cA$ on a Hilbert
       space $\cE$.  Rather than the algebra $\cF^{\infty}(E)$ of
       adjointable operators on the $ \cA$-correspondence
$\cF^{2}(E)$, our main
       focus of interest will be on the algebra $\cF^{\infty}(E) \otimes
       I_{\cE}$  of all operators on the Hilbert space $\cF^{2}(E,
       \sigma)$ of the form $R = T \otimes I_{\cE}$ with $T \in
\cF^{\infty}(E)$
       acting on the Hilbert space $\cF^2(E,\sigma)$.
       Note that the operator $R=T\otimes I_\cE$ is properly
       defined since $R$ is an $\cA$-module map with respect to the right
       $\cA$-action on $\cF^2(E)$. For convenience we shall use the
       abbreviated notation
       \begin{equation*}  %\label{HinfEsigma}
         \cF^{\infty}(E,\sigma) = \cF^{\infty}(E) \otimes I_{\cE},
       \end{equation*}
       and
       \[
       \varphi_{\infty,\sigma}(a)=\varphi_{\infty}(a)\otimes
I_\cE\quad\mbox{and}\quad
       T_{\xi,\sigma}=T_\xi\otimes I_\cE
       \quad\mbox{for all }a\in\cA\mbox{ and }\xi\in E.
       \]
       The algebra
       $\cF^{\infty}(E, \sigma)$ can also be
       described as the weak-$*$ closed algebra generated by the collection of operators
       \begin{equation}  \label{genset}
       \{ \varphi_{\infty,\sigma}(a), \, T_{\xi,\sigma}
       \colon a \in \cA, \, \xi \in E \}.
       \end{equation}

       The following
       alternative characterization of $\cF^{\infty}(E, \sigma)$
       will be useful.  Here we define $E^{\sigma}$ and $ \sigma(\cA)'$
       as in (\ref{Esigma0}) and (\ref{sigmaA'0}).
       Note that each element $\mu$ of $E^{\sigma}$ induces a {\em dual
       creation operator} $T^{d}_{\mu,\sigma}$ in
$\cL^a(\cF^{2}(E, \sigma))$ given by
       $$
         T^{d}_{\mu,\sigma} =
         \begin{bmatrix} 0 & 0 & 0 & \cdots \\ T_{\mu,\sigma}^{d,(0)} & 0
       &
           0  & \cdots \\
          0 & T_{\mu,\sigma}^{d,(1)} & 0 & \cdots & \\
           \ddots & \ddots & \ddots  &  \end{bmatrix}
       $$
       where $T_{\mu,\sigma}^{d,(n)} \colon E^{\otimes n}
       \otimes_{\sigma} \cE \to E^{\otimes n+1} \otimes_{\sigma} \cE$ is
       given by
       $$
       T_{\mu, \sigma}^{d,(n)} \colon \xi_{n} \otimes \cdots \otimes
       \xi_{1} \otimes e \mapsto
       \xi_{n} \otimes \cdots \otimes
       \xi_{1} \otimes \mu e
       $$
       where as usual we make the identification
       $$ E^{\otimes n} \otimes (E \otimes_{\sigma} \cE) =
       E^{\otimes n+1} \otimes_{\sigma} \cE.
       $$
       Using the notation in (\ref{tensoroperator}) we can
       write $T_{\mu,\sigma}^d=I_{\cF^2(E)}\otimes \mu$, where we identify
       $\cF^2(E)\otimes_{\cA} E$ with $\cF^2(E)$, which makes sense because
       $\mu$ is an $\cA$-module map. Also recall that $\iota_{\infty,
       \sigma}$ in (\ref{iota}) defines  a $*$-representation of
       $\sigma(\cA)'$ on $\cF^2(E,\sigma)$.

       \begin{proposition} \label{P:Toeplitz}

       An operator $R \in
         \cL(\cF^{2}(E, \sigma))$ is in $\cF^{\infty}(E,
\sigma)$ if and
         only if $R$ commutes with each of the operators
         $I_{\cF(E)} \otimes_{\sigma} b$ and $T_{\mu,\sigma}^{d}$
         for $b \in \sigma(\cA)'$ and $\mu \in E^{\sigma}$.
         Consequently, the operator $R \in \cL(\cF^2(E, \sigma))$
        with infinite block-matrix representation
         $$
           R = [ R_{i,j}]_{i,j=0,1,2, \dots} \text{ where } R_{i,j} \colon
           E^{\otimes j} \otimes_{\sigma} \cE \to E^{\otimes i}
           \otimes_{\sigma} \cE
         $$
         is in $\cF^{\infty}(E, \sigma)$ if and only if $R$ is
         lower triangular ($R_{i,j} = 0$ for $i<j$) and for $i \ge j$
       $R_{i,j}$
         satisfies the following compatibility (Toeplitz-like) conditions:
       \begin{align}
         & R_{i,j} (I_{E^{\otimes j}} \otimes b) = (I_{E^{\otimes i}}
         \otimes b) R_{i,j} \text{ for all } b \in \sigma(\cA)',
         \label{Toeplitz1} \\
          & R_{i+1,j+1} (I_{E^{\otimes j}} \otimes \mu) =
         (I_{E^{\otimes i}} \otimes \mu) R_{i,j} \text{ for all }
         \mu \in E^{\sigma}.  \label{Toeplitz2}
         \end{align}
         and hence, inductively,
         \begin{equation}  \label{Toeplitz3}
        R_{i,j}\mu^j=(I_{E^{\otimes i-j}}\otimes \mu^j)R_{i-j,0},
          \end{equation}
          where $\mu^j=((\mu^*)^j)^*$, with $(\mu^*)^j$ the
generalized power of
         $\mu^*\in (E^{\sigma})^{*}$.
         \end{proposition}

         \begin{proof}
        The first part follows from Theorem 3.9 of
\cite{MS-Annalen}. The
       second part
         is then a straightforward translation of these
       commutativity conditions to expressions involving the block entries.
       \end{proof}

       Taking the cue from Proposition \ref{P:Toeplitz}, we view
elements $R$
       of $\cF^{\infty}(E, \sigma)$ as the {\em analytic Toeplitz
       operators} for this Fock-space/correspondence setting.

       While it is in general not the case that $R
       (T_{\mu , \sigma}^{d})^{*} =
       (T_{\mu,\sigma}^{d})^{*} R$ for $R \in
       \cF^{\infty}(E, \sigma)$ and $\mu \in E^{\sigma}$, this
       is almost the case as is made precise in the following proposition.

       \begin{proposition}  \label{P:almost}
        For $R \in \cF^{\infty}(E, \sigma)$ and $\mu \in
        E^{\sigma}$, we have
          \begin{equation}  \label{almost-op}
          R (T_{\mu,\sigma}^{d})^{*}|_{\oplus_{n=1}^{\infty}
          E^{\otimes n} \otimes \cE} = (T_{\mu,\sigma}^{d})^{*}
          R|_{\oplus_{n=1}^{\infty} E^{\otimes n} \otimes \cE},
          \end{equation}
          or, in terms of matrix entries, we have inductively
          \begin{equation}  \label{almost}
          R_{i,j} (I_{E^{\otimes j}} \otimes \eta) = (I_{E^{\otimes i}}
          \otimes \eta ) R_{i+1,j+1} \text{ for all } i,j=0,1, \dots.
          \end{equation}
         for $\eta = \mu^{*} \in (E^{\sigma})^{*}$.
          \end{proposition}

          \begin{proof}
          To prove that \eqref{almost-op} holds for all $R \in  \cF^{\infty}(E,
          \sigma)$, it suffices to show that it holds for each $R$ in the
          generating set \eqref{genset}.  We are thus reduced to
showing that
          \eqref{almost} holds for all $R$ of the special form
          $\varphi_{\infty}(a)$ for an $a \in \cA$ and $T_{\xi}$ for a $\xi
          \in E$.  This in turn is a routine
          calculation which we leave to the reader.
          \end{proof}

          Suppose that we are given $R \in \cF^{\infty}(E, \sigma)$.
          We regard $\cE$ as a subspace of $\cF^{2}(E, \sigma)$  via the
          identification $y \cong y \oplus 0 \oplus  0\oplus \cdots$.
          Then the restriction of $R$ to $\cE$ defines an operator
          from $\cE$ into $\cF^{2}(E, \sigma)$ where we have a point
          evaluation in ${\mathbb D}((E^{\sigma})^{*})\times \sigma(\cA)'$
          defined in \eqref{pteval}.
          We may then define an operator $R^{\wedge}(\eta)
          \in \cL(\cE)$ by
          \begin{equation*}%\label{SCpointev}
          R^{\wedge}(\eta) e = (R e)^{\wedge}(\eta,I_\cE).
          \end{equation*}
          Explicitly, we have
          $$
         R^{\wedge}(\eta) = \sum_{n=0}^{\infty} \eta^{n} R_{n,0} \in
           \cL(\cE).
          $$
          Note that, as a consequence of Proposition \ref{P:Toeplitz},
           the full function ${\mathbf f}(\eta,b) = (R
e)^{\wedge}(\eta, b)$
           is then determined from $R^{\wedge}(\eta)$
          and $e \in \cE$ according to
          \[
          (Re)^{\wedge}(\eta,b)=(\iota_{\infty,\sigma}(b)R e)(\eta,I_\cE)
          =(R\iota_{\infty,\sigma}(b)e)(\eta,I_\cE)
          =(Rbe)(\eta) = R^{\wedge}(\eta) (b e)
          \]
           for $\eta\in{\mathbb D}((E^{\sigma})^{*})$ and
$b\in\sigma(\cA)'$.
          This implies that if we would extend the point evaluation to
          ${\mathbb D}((E^{\sigma})^{*})\times \sigma(\cA)'$ by
          $R^{\wedge}(\eta, b)e = (Re)^{\wedge}(\eta, b)$,
          the result would just give
          $R^{\wedge}(\eta,b)=R^{\wedge}(\eta)b$.

          It is of interest that this transform $R \to
R^{\wedge}(\cdot)$ is
          multiplicative.

          \begin{proposition}  \label{P:multiplicative}
          \begin{enumerate}
          \item
          Suppose that $R$ and $S$ are two elements of
          $\cF^{\infty}(E, \sigma)$. Then
          \begin{equation*}  %\label{operator-mult}
              (R S)^{\wedge}(\eta) = R^{\wedge}(\eta)
              S^{\wedge}(\eta)
          \end{equation*}
          for all $\eta \in {\mathbb D}((E^{\sigma})^{*})$.

          \item Suppose that $R$ is an
          operator in $\cF^\infty(E,\sigma)$
            and that $f$ is an element of
          $\cF^2(E,\sigma)$.  Then
          \begin{equation}  \label{operator-funct-mult}
              (Rf)^{\wedge}(\eta,b)
              = R^{\wedge}(\eta)
              f^{\wedge}(\eta,b)
        \end{equation}
        for all $\eta \in {\mathbb D}((E^{\sigma})^{*})$ and
        $b\in\sigma(\cA)'$.
          \end{enumerate}
          \end{proposition}

          \begin{proof}
          Suppose that $R = [R_{i,j}]_{i,j=0,1,\dots}$ is an
operator in
          $\cF^{\infty}(E, \sigma)$ and that $f = \oplus_{j=0}^{\infty}
          f_{j}$ is an element of $\cF^{2}(E, \sigma)$.
          We first note that
          a special case of \eqref{almost} is
          $$ R_{\ell,0} \eta =
(I_{E^{\otimes\ell}}\otimes\eta) R_{\ell +1,1}.
          $$
          Iteration of \eqref{almost} in turn leads to
          \begin{equation}  \label{almost-iterated}
          R_{\ell,0} \eta^{j}
          = I_{E^{\otimes\ell}}\otimes\eta^{j} R_{\ell + j, j} \colon
          E^{\otimes j} \otimes_{\sigma}\cE \to E^{\otimes \ell}
          \otimes_{\sigma} \cE.
          \end{equation}
          Then we compute for
          $\eta\in{\mathbb D}((E^{\sigma})^{*})$ and
$b\in\sigma(\cA)'$ that
          \begin{align*}
          R^{\wedge}(\eta)  f^{\wedge}(\eta,b) & =
          \left(\sum_{\ell=0}^{\infty} \eta^{\ell} R_{\ell,0}\right)
          \left( \sum_{j=0}^{\infty} \eta^{j}
           (I_{E^{\otimes j}}\otimes b) f_j \right) \\
          & = \sum_{\ell,j = 0}^{\infty} \eta^{\ell}
R_{\ell,0} \eta^{j}
           (I_{E^{\otimes j}}\otimes b) f_j \\
           & = \sum_{\ell,j = 0}^{\infty} \eta^{\ell + j}
R_{\ell + j,j}
           (I_{E^{\otimes j}}\otimes b) f_j \text{
          (by \eqref{almost-iterated}) }\\
          & = \sum_{\ell,j = 0}^{\infty} \eta^{\ell + j}
         (I_{E^{\otimes\ell+j}}\otimes b) R_{\ell + j,j}
           f_j \text{
          (by \eqref{Toeplitz1})} \\
          & = \sum_{n=0}^{\infty} \eta^{n}
          (I_{E^{\otimes n}}\otimes b)
          \left( \sum_{j=0}^{n}
          R_{n,j} f_j \right) \\
          & = \sum_{n=0}^{\infty} \eta^{n}
          (I_{E^{\otimes n}}\otimes b)
          [ R f ]_{n} = (R
          f)^{\wedge}(\eta,b)
          \end{align*}
          and part (2) of the Proposition follows.
         Part (1) follows as the special case where $b=I_{\cE}$ and
          $f=Se$ for arbitrary $e\in\cE$.
          \end{proof}

          \begin{remark}  \label{R:H2modmap} {\em
          We note that a consequence of the formula
          \eqref{operator-funct-mult} is that the operator
$M_{R^\wedge}$
          of multiplication by $R^{\wedge}$ on $H^{2}(E, \sigma)$
          $$
           M_{R^{\wedge}} \colon  f^{\wedge}(\eta, b) \mapsto
           R^{\wedge}(\eta)  f^{\wedge}(\eta, b)
           $$
           commutes with the $\sigma(\cA)'$-left action on
$H^{2}(E, \sigma)$:
           $$
        M_{R^{\wedge}}(b \cdot f^{\wedge}) = b \cdot M_{R^{\wedge}}
        f^{\wedge} \text{ where } (b \cdot f)^{\wedge}(\eta',b') =
         f^{\wedge}(\eta', b'b)
        $$
        for all $b,b' \in \sigma(A)'$ and $\eta' \in {\mathbb
        D}((E^{\sigma})^{*})$.  This can also be seen as a
consequence of
        applying the $Z$-transform to the identity
        $$
         R \iota_{\infty, \sigma}(b) = \iota_{\infty,
\sigma}(b) R \text{
         for all } b \in \sigma(\cA)'
           $$
           given in Proposition \ref{P:Toeplitz}.
           }\end{remark}

          Proposition \ref{P:multiplicative} leads immediately to the
          following corollary.

          \begin{corollary}  \label{C:ker-inv}
        \begin{enumerate}
       \item The kernel of the Fourier transform $\Phi \colon f \to
           f^{\wedge}$ in $\cF^{2}(E, \sigma)$
          $$
          \operatorname{Ker} \Phi = \{ f \in \cF^{2}(E, \sigma) \colon
       f^{\wedge}(\eta,b) = 0 \text{ for all } \eta \in {\mathbb
          D}((E^{\sigma})^{*}) \text{ and } b \in \sigma(\cA)' \}
          $$
          is invariant under the analytic Toeplitz operators:
          \begin{align*}
          & f^{\wedge} (\eta,b) = 0 \text{ for all } \eta \in {\mathbb
          D}((E^{\sigma})^{*}) \text{ and } b \in \sigma(\cA)', \,
          R \in \cF^{\infty}(E, \sigma) \\
          & \qquad
          \Longrightarrow (R f)^{\wedge}(\eta,b) = 0
          \text{ for all } \eta \in {\mathbb D}((E^{\sigma})^{*})
\text{ and }
          b \in \sigma(\cA)'.
          \end{align*}

          \item The initial space
          $[\cF^2(E, \sigma)]_{\text{\rm initial}}$
          of the Fourier transform $\Phi$ is invariant under the
          adjoints of the analytic
          Toeplitz operators:
          $$ f \in [\cF^{2}(E, \sigma)]_{\text{\rm initial}}, \, R \in
       \cF^{\infty}(E,
          \sigma) \Longrightarrow R^{*}f \in [\cF^{2}(E,
\sigma)]_{\text{\rm
       initial}}.
          $$
          Explicitly, the action of $R^{*}$ on a generic vector in
          the spanning set \eqref{initial} for $[\cF^{2}(E)
          \otimes_{\sigma}\cE]_{\text{\rm initial}}$ is given by
          $$
         R^{*} (b^{*} \cdot k_{E,\sigma; \eta})e) = b^{*} \cdot
         k_{E,\sigma; \eta}R^{\wedge}(\eta)^{*}e.
          $$
          \end{enumerate}
          \end{corollary}

          \begin{proof}
          If $f^{\wedge}(\eta,b) = 0$ for all $\eta$ and $b$, then, by
          \eqref{operator-funct-mult} we see immediately that
        $$ (R f)^{\wedge}(\eta,b) =
        R^{\wedge}(\eta)f^{\wedge}(\eta,b) = 0
        $$
        for all $\eta$ and $b$ as well as for any $R \in
\cF^{\infty}(E, \sigma)$.
           The first part of the second statement then follows by simply
           taking adjoints.

           To verify the second part of the second statement, it
suffices to
           verify on the generators $R = T_{\xi}$ and $R = \varphi(a)$ for
           $\xi \in E$ and $a \in \cA$; this in turn is straightforward.
         \end{proof}

         \begin{remark}
        {\em  We note that the definition of
$R^{\wedge}(\eta)$ involves only
         the first column of $R$.  From the relations \eqref{Toeplitz3}
         and \eqref{Toeplitz1} one can see that the first column of
         $R$ already uniquely determines the action of $R$ on all of
         $[\cF^{2}(E,\sigma)]_{\text{initial}}$.}
          \end{remark}

       \begin{remark} \label{R:etabmu} {\em
       Let $\mu \in E^{\sigma}$ and $\eta \in (E^{\sigma})^*$ and $b \in
       (\sigma(\cA))'$.
       Then an easy verification using the relations
       $\mu \sigma(a) = (\varphi(a) \otimes I_{\cE})\mu $ and
$\sigma(a)\eta  =
       \eta (\varphi(a) \otimes I_{\cE})$ shows that
       \begin{equation}  \label{etabmu}
        \eta (I_E \otimes b)\mu \in \sigma(\cA)'.
       \end{equation}
       This observation has several
       consequences.

       \begin{enumerate}
        \item
       Given $\mu \in E^{\sigma}$
       and $\eta \in (E^{\sigma})^{*}$ we may define a mapping
       $\theta_{\eta, \mu}$ on $\sigma(\cA)'$ by
       \begin{equation*}  %\label{theta-eta-mu}
        \theta_{\eta, \mu}(b) = \eta (I_{E} \otimes b) \mu.
       \end{equation*}
       Iteration of this map gives
       $$
         \theta^{2}_{\eta, \mu}(b) = \eta (I_{E }\otimes \eta
(I_{E} \otimes
         b) \mu) \mu = \eta^{2}(I_{E^{\otimes 2}} \otimes b) \mu^{2}
       $$
       and more generally
       $$
       \theta_{\eta, \mu}^{n}(b) = \eta^{n}(I_{E^{\otimes n}} \otimes b)
       \mu^{n}
       $$
       where we make use of the generalized power $\eta^{n}$ for an element
       $\eta$ of $(E^{\sigma})^{*}$ (and set $\mu^{n} = (
(\mu^{*})^{n})^{*}
       \colon \cE \to E^{\otimes n} \otimes_{\sigma} \cE$).  For $\eta,
       \zeta \in {\mathbb D}((E^{\sigma})^{*})$, we may take $\mu =
       \zeta^{*}$ and then we have $\| \theta_{\eta, \zeta^{*}}\| < 1$.
       Then we may use the geometric series to compute the inverse of $I -
       \theta_{\eta, \zeta^{*}}$ to get
       $$ (I - \theta_{\eta, \zeta^{*}})^{-1}(b) = \sum_{n=0}^{\infty}
       (\theta_{\eta, \zeta^{*}})^{n}(b) =
       \sum_{n=0}^{\infty} \eta^{n} (I_{E^{\otimes n}} \otimes b)
       (\zeta^{n})^{*}.
       $$
       We conclude that the Szeg\"o kernel \eqref{E-sigma-ker} can also be
       written as
       \begin{equation*}
        \K_{E,\sigma}(\eta, \zeta)[b] = (I - \theta_{\eta,
        \zeta^{*}})^{-1}(b).
       \end{equation*}
       This is the form of the Szeg\"o kernel used in \cite{MS-Annalen,
       MS-Schur}.

       \item Suppose that we are given two elements $\eta, \zeta
\in E^{\sigma}$.
         The special case of \eqref{etabmu} with $b = I_{\cE}$ and
         $\eta = \mu^{\prime *}$
       for a $\mu' \in E^{\sigma}$ enables us to define a
$\sigma(\cA)'$-valued
       inner product on $E^{\sigma}$:
       $$ \langle \mu, \mu' \rangle_{E^{\sigma}} = \mu^{\prime *} \mu \in
       \sigma(\cA)' \text{ for } \mu, \mu' \in E^{\sigma}.
       $$
       Moreover one can check that $E^{\sigma}$ has a well-defined right
       $\sigma(\cA)'$-action
       $$
          (\mu \cdot b)(e) = \mu( b e)
       $$
       and a well-defined left $\sigma(\cA)'$-action
       $$
          (b \cdot \mu)(e) = (I_{E} \otimes b) \mu(e).
       $$
       It is then straightforward to check that $E^{\sigma}$ is a
       $\sigma(A)'$-correspondence.  This observation plays a key role in
       the duality theory in \cite{MS-Annalen} (see also Proposition
       \ref{P:Toeplitz} above).
       \end{enumerate}
       }
       \end{remark}

       Next we introduce
       the space
       \[
       H^\infty(E,\sigma)=\{R^\wedge\colon R\in\cF^\infty(E,\sigma)\},
       \]
       where we interpret $R^\wedge$ as a function mapping
       ${\mathbb D}((E^{\sigma})^{*})$ into $\cL(\cE)$.
       Then $H^\infty(E,\sigma)$ is closed under addition
       ($(R_1+R_2)^\wedge=R_1^\wedge+R_2^\wedge$), scalar multiplication
       ($(\lambda R)^\wedge=\lambda R^\wedge$) and pointwise multiplication
       (Proposition \ref{P:multiplicative} (1)). Moreover, part (2) of
       Proposition \ref{P:multiplicative} implies that a function $S\in
       H^\infty(E,\sigma)$ defines a multiplication operator $M_S$ on
       $H^2(E,\sigma)$ by
       \begin{equation}\label{multop}
       (M_S f^\wedge)(\eta,b)=S(\eta)f^{\wedge}(\eta,b)\mbox{
       for each }\eta\in{\mathbb D}((E^{\sigma})^{*}),b\in\sigma(\cA)',
       f^\wedge\in H^2(E,\sigma).
       \end{equation}

       In fact, we have the following result.

       \begin{proposition}\label{P:Hinftycar}
       A function $S:{\mathbb D}((E^{\sigma})^{*})\to\cL(\cE)$ is in
       $H^\infty(E,\sigma)$ if and only if $S$ defines a multiplication
       operator $M_S$ on $H^2(E,\sigma)$ by (\ref{multop}). In case $S\in
       H^\infty(E,\sigma)$, we have $\|M_S\|\leq\|R\|$ for each
       $R\in\cF^\infty(E,\sigma)$ with $S = R^{\wedge}$ and there exists a
       $R\in\cF^\infty(E,\sigma)$ with $S = R^{\wedge}$ such that
$\|M_S\|=\|R\|$.
       Moreover, if
       $S\in H^\infty(E,\sigma)$, then $M_S$ is a $\sigma(\cA)$-module map
       that in addition commutes with the operators
       \begin{equation*}%\label{comfuncs}
       \Phi(I_{\cF(E)}\otimes \mu)\Phi^*\mbox{ for each }\mu\in E^\sigma.
       \end{equation*}
       Here $\Phi$ is the coisometry from $\cF^2(E,\sigma)$ into
       $H^2(E,\sigma)$ given by $\Phi:f\mapsto f^\wedge$.
       \end{proposition}

       \begin{proof}
       We already observed that $S\in H^\infty(E,\sigma)$ guarantees that
       $M_S$ in (\ref{multop}) defines a multiplication operator
on $H^2(E,\sigma)$.
       Moreover, for $R\in\cF^\infty(E,\sigma)$ with $S=R^\wedge$ we have
       \[
       (M_S\Phi f)(\eta,b)=(M_S f^\wedge)(\eta,b)=S(\eta)f^\wedge(\eta,b)
       =R^{ \wedge}(\eta)f(\eta,b)=(Rf)^\wedge(\eta,b)=(\Phi
       Rf)(\eta,b)
       \]
       for each $f\in\cF^2(E,\sigma)$, $\eta\in{\mathbb
       D}((E^{\sigma})^{*})$ and $b\in\sigma(\cA)'$. Hence
       \[
       M_S\Phi=\Phi R.
       \]
       In particular we have $M_S=\Phi R\Phi^*$ and thus
       $\|M_S\|\leq\|R\|$ since $\Phi$ is a coisometry.

       Now assume that $S$ defines a multiplication
       operator $M_S$ on $H^2(E,\sigma)$ by (\ref{multop}).
       The definition of $M_{S}$ and of the left action on $H^2(E,\sigma)$
       in (\ref{action})
       shows that, for $b,b'\in\sigma(\cA)'$ and $\eta\in{\mathbb
       D}((E^{\sigma})^{*})$,
       we have
       \[
       (M_S
b'f^\wedge)(\eta,b)=S(\eta)f^\wedge(\eta,bb')=(b'M_Sf^\wedge)(\eta,b)
       \quad\mbox{for each }f^\wedge\in H^2(E,\sigma).
       \]
       Hence $M_S$ is a $\sigma(\cA)'$-module map.

         We now show that there exists $R \in \cF^{\infty}(E,\sigma)$
       with $R^{\wedge} = S$. We first note that
       \begin{eqnarray}
       \nonumber
       ((I_{\cF^2(E)}\otimes\mu)f)^\wedge(\eta,b)&=&\sum_{n=1}^\infty\eta^n
       (I_{E^{\otimes n}}\otimes b)(I_{E^{\otimes
       n-1}}\otimes\mu)f_{n-1}\\
       \label{invcomp}
       &=&\sum_{n=1}^\infty \eta^{n-1}(I_{E^{\otimes n-1}}\otimes
       \eta(I_E\otimes b)\mu)f_{n-1}\\
       &=&f^\wedge(\eta,\eta(I_E\otimes b)\mu)  \label{invcomp1}
       \end{eqnarray}
       where we use the observation from Remark \ref{R:etabmu} that $\eta
       (I_{E} \otimes b) \mu$ is in $\sigma(\cA)'$.

        From (\ref{invcomp1}), it readily follows that
$I_{\cF^2(E)} \otimes \mu$
       on $\cF^2(E,\sigma)$ leaves
       $\operatorname{Ker} \Phi $ invariant. The same holds for
the operator
       $I_{\cF^2(E)} \otimes b$. Consequently, denoting by $P (=
       \Phi^*\Phi)$ the projection
       on ${\cG}=(\ker \Phi)^{\perp}$, we note that
       $$
       PX = PXP \ \mbox{for}\ X=I_{\cF^2(E)} \otimes \mu,
I_{\cF^2(E)} \otimes b.
       $$

       We  show that the operator $\Phi^*M_S\Phi $ commutes with
       $I_{\cF(E)}\otimes b'$ for all $b' \in (\sigma(\cA))'$.
       To see this, let $f \in \cF^2(E,\sigma)$ and
       $\Phi^*M_S\Phi(I_{\cF(E)}\otimes b')f =g$.
       Due to (\ref{invcomp1}), we have
       $g^{\wedge}(\eta,b) = S(\eta, bb')$. Now if we let
$\Phi^*M_S\Phi f = h$,
       it follows that $h^{\wedge}(\eta,b) = S(\eta)f^{\wedge}(\eta,b)$ and
       consequently,
       $$
       ((I_{\cF(E)}\otimes b')\Phi^*M_S \Phi f)^\wedge (\eta,b) =
       S(\eta)\hat{f}(\eta,bb')
       $$
         and the claim follows.
       A similar computation using (\ref{invcomp1}) shows that
         $P(I_{\cF(E)}\otimes \mu)A = A P(I_{\cF(E)}\otimes
\mu)|\cG$ for all
         $\mu \in E^{\sigma}$, where
         $A = \Phi^*M_S \Phi |\cG$.

         We recall now that the maps $\mu \in E^{\sigma}, b \in
\sigma(\cA)'$
         form an isometric covariant representation of the
         $\sigma(\cA)'$-correspondence $E^{\sigma}$
         (see pages 369-370 in \cite{MS-Annalen}---the precise
definition is
         {\em covariant representation} is
         given in the text surrounding formulas
         \eqref{sigma}--\eqref{covariant} below).
         We may now apply the
         commutant lifting theorem for covariant representations of a
         correspondence due to Muhly-Solel (see Theorem 4.4, \cite{MS-JFA})
         to obtain an operator
         that commutes with the operators $I_{\cF^2(E)}\otimes \mu$ and  $
         I_{\cF^2(E)}\otimes b$ (which implies $R \in \cF^{\infty}(E)$ by
         Proposition \ref{P:Toeplitz}) which moreover satisfies $PR = AP$.
         This immediately implies that $R^{\wedge} = S$.
Furthermore, we can choose
         $R$ such that $\|R\| = \|M_S\|$.
         \end{proof}

          We note that any $R \in \cF^{\infty}(E, \sigma)$ is of the form
          $\widetilde R \otimes I_{\cE}$ for a $\widetilde R \in
          \cF^{\infty}(E)$.  Moreover, the map $\widetilde R
\mapsto R = \widetilde
          R \otimes I_{\cE}$ is an $\cL(\cF^{2}(E, \sigma))$-valued
          representation of $\cF^{\infty}(E)$ which actually extends to a
          $*$-representation $T \mapsto T \otimes I_{\cE}$ of all of
          $\cL^{a}(\cF^{2}(E))$---the restriction of $T \mapsto T \otimes
          I_{\cE}$ to $T \in \cF^{\infty}(E)$ is called the {\em induced
       representation}
          of $\cF^{\infty}(E)$ in the terminology of
\cite{MS-JFA, MS-Annalen}.
          The content of Proposition
          \ref{P:multiplicative} is that, for each $\eta \in {\mathbb
          D}((E^{\sigma})^{*})$, the map $R \mapsto R^{\wedge}(\eta)$ is an
          $\cL(\cE)$-valued representation of $\cF^{\infty}(E,
\sigma)$.  It
          follows that the composition
          \begin{equation}  \label{pi-eta}
        \pi_{\eta}(\widetilde R) = (\widetilde R \otimes
        I_{\cE})^{\wedge}(\eta)
          \end{equation}
          is an (even completely contractive) representation of
          $\cF^{\infty}(E)$ (see \cite{MS-Annalen}).
          For {\em some} $\eta \in    (E^{\sigma})^{*}$ of norm
equal to $1$,
          $\pi_{\eta}$ still defines a representation of $\cF(E)$.
          It is the case that every $\eta$ in the closed unit ball of
          $(E^{\sigma})^{*}$ gives rise to a completely
contractive representation
          of $\cT_{+}(E)$ (the norm-closure of the span of left multipliers
          $\varphi_{\infty}(a)$ ($a \in \cA$) and creation
operators $T_{\xi}$
          ($\xi \in E$) in $\cL^{a}(\cF^{2}(E))$), while it is
not clear for
          which such $\eta$ the representation can be extended to
          $\cF^{\infty}(E)$---this is one of the
          open problems in the theory (see \cite{MS-Annalen}).
It is the case
          that each completely contractive representation $\pi$
of $\cF^{\infty}(E)$
          comes from an $\eta \in \overline{\mathbb
D}((E^{\sigma})^{*})$ for
          some weak-$*$  continuous $*$-representation
          $\sigma \colon \cA \to \cL(\cE)$.  Indeed,
          given a completely contractive representation $\pi \colon
          \cF^{\infty}(E) \to \cL(\cE)$, one can construct $\sigma$ and
       $\eta$ as follows.
          Define $\sigma \colon \cA \to \cL(\cE)$ by
          \begin{equation}  \label{sigma}
           \sigma(a) = \pi(\varphi_{\infty}(a)).
          \end{equation}
          Then define $\boldsymbol \eta \colon E \to \cL(\cE)$ by
          \begin{equation}  \label{boldeta}
           \boldsymbol \eta(\xi) = \pi(T_{\xi}).
          \end{equation}
          We wish to verify that
          \begin{equation}  \label{covariant}
           \boldsymbol \eta(\varphi(a) \xi \cdot a') = \sigma(a)
\boldsymbol
           \eta(\xi) \sigma(a'),
          \end{equation}
          i.e., that the pair $(\boldsymbol \eta, \sigma)$ is a
{\em covariant
          representation} of $E$ in the terminology of Muhly-Solel
          \cite{MS-JFA, MS-Annalen}. As a first step for the
verification of
          \eqref{covariant}, one can easily check that
          $$
        T_{\varphi(a) \xi \cdot a'} = \varphi_{\infty}(a) T_{\xi}
        \varphi_{\infty}(a') \text{ for } a \in \cA,\, \xi \in E.
          $$
          We then compute
          \begin{align*}
          \boldsymbol \eta(\varphi(a) \xi \cdot a') & =
\pi(T_{\varphi(a)
          \xi \cdot a'}) \\
          & = \pi (\varphi_{\infty}(a) T_{\xi} \varphi_{\infty}(a')) \\
          & = \pi(\varphi_{\infty}(a)) \pi(T_{\xi})
          \pi(\varphi_{\infty}(a')) \\
          & = \sigma(a) \boldsymbol \eta(\xi) \sigma(a')
          \end{align*}
          and \eqref{covariant} follows.  As in \cite{MS-JFA}, a covariant
          representation $(\boldsymbol \eta, \sigma)$ of $E$ determines an element
          $\eta \colon E \otimes_{\sigma} \cE \to \cE$ of
$(E^{\sigma})^{*}$
          according to the formula
          \begin{equation}  \label{eta}
          \eta(\xi \otimes e) = \boldsymbol \eta(\xi) e.
          \end{equation}
          Here note that the property $\boldsymbol{\eta}(\xi \cdot a') =
          \boldsymbol{\eta}(\xi) \sigma(a')$ is what is needed to
verify that
          \eqref{eta} is well-defined while the property $\boldsymbol
          \eta(\varphi(a) \xi) = \sigma(a) \boldsymbol \eta(\xi)$
is what is
          needed to verify that $\eta$ is in $(E^{\sigma})^{*}$, i.e., that
          $\eta$ has the $\cA$-module-map property
          $$
        \eta(\varphi(a) \otimes I_{\cE}) = \sigma(a) \eta.
          $$

          There is a converse:  given an element $\eta \in {\mathbb
          D}((E^{\sigma})^{*})$, we may use \eqref{eta} to define
       $\boldsymbol \eta$ so that
          $(\boldsymbol \eta, \sigma)$ is a  completely
contractive  covariant
          representation of $E$.
        The mapping $\pi $ given in
\eqref{sigma}---\eqref{boldeta} then
       extends to a representation
          of $\cF^{\infty}(E)$ (see \cite{MS-JFA}).  For our
          situation here where $(\boldsymbol \eta, \sigma)$ is
given in terms
          of a representation $\pi$ via
\eqref{sigma}---\eqref{boldeta}, we see that a
          representation $\pi$ of $\cF^{\infty}(E)$ determines a
representation
          $\sigma = \sigma_{\pi}$ of $\cA$ according to \eqref{sigma} along
       with an element
          $\eta_{\pi}$ of $(E^{\sigma})^{*}$ according to the formula
          \begin{equation*}  %\label{eta-pi}
        \eta_{\pi}(\xi \otimes e) = \pi(T_{\xi}) e.
          \end{equation*}
          It is then straightforward to check that the formula
          \begin{equation}  \label{FinfErep-eta}
        \pi(\widetilde R) = (\widetilde R \otimes
        I_{\cE})^{\wedge}(\eta_{\pi})
          \end{equation}
          holds for the cases where
          $$
          \widetilde R = \varphi_{\infty}(a) \text{ for some } a
\in \cA, \qquad
          \widetilde R = T_{\xi} \text{ for some } \xi \in E.
          $$
          Under the assumption that $\pi$ is continuous with respect to the
          weak-$*$ topologies on $\cF^{\infty}(E)$ and
$\cL(\cE)$, it then follows
          that \eqref{FinfErep-eta} holds for all $\widetilde R \in
          \cF^{\infty}(E)$, i.e., we recover $\pi$ as $\pi =
\pi_{\eta_{\pi}}$
          where in general $\pi_{\eta}$ is given by \eqref{pi-eta}.

          It is of interest to apply this construction to the induced
          representation
          \begin{equation}  \label{induced}
        \pi_{\text{ind}} \colon \widetilde R \mapsto
\widetilde R \otimes I_{\cE}
          \end{equation}
          of $\cF^{\infty}(E)$ into $\cL(\cF^{2}(E, \sigma))$.
We collect this
          result in the following Proposition.

          \begin{proposition}  \label{P:induced}
         Suppose that we are given an $\cA$-correspondence $E$ together
         with a representation $\sigma \colon \cA \to \cL(\cE)$ for a
         Hilbert space $\cE$ and let $\pi_{\text{\rm ind}} \colon
         \cF^{\infty}(E) \to \cL(\cF^{2}(E, \sigma))$ be the induced
       representation as in
         \eqref{induced}.  Define $\boldsymbol
\eta_{\text{\rm ind}} \colon
         E \to \cL(\cF^{2}(E, \sigma))$ and
$\sigma_{\text{\rm ind}} \colon
         \cA \to \cL(\cF^{2}(E, \sigma))$ by
         $$ \boldsymbol \eta_{\text{\rm ind}}(\xi) =
T_{\xi,\sigma}, \qquad
         \sigma_{\text{\rm ind}}(a) = \varphi_{\infty, \sigma}(a).
         $$
         Then $(\boldsymbol \eta_{\text{\rm ind}},
         \sigma_{\text{\rm ind}})$ is an  (isometric) covariant
       representation of $E$
         with element $\eta_{\text{\rm ind}}
         \colon E \otimes \cF^{2}(E, \sigma)
         \to  \cF^{2}(E, \sigma)$ of $(E^{\sigma_{\text{\rm
ind}}})^{*}$
         associated with
         $(\boldsymbol \eta_{\text{\rm ind}},
\sigma_{\text{\rm ind}})$ as
         in \eqref{eta}
         given by
         $$ \eta_{\text{\rm ind}} \colon
         \xi \otimes \left[  \oplus_{n=0}^{\infty} \xi^{(n)}
\otimes e_{n}
         \right] \mapsto
         0 \oplus \left[ \oplus_{n=1}^{\infty} \xi \otimes \xi^{(n-1)}
         \otimes e_{n-1} \right].
         $$
         Moreover,  we
         recover $R = \widetilde R \otimes I_{\cE} \in
\cF^{\infty}(E, \sigma)$
         via the point evaluation
         \begin{equation*}  %\label{recover}
         \widetilde R \otimes I_{\cE}  =  (\widetilde R
\otimes I_{\cF^{2}(E,
         \sigma)})^{\wedge}(\eta_{\text{\rm ind}}).
         \end{equation*}
         \end{proposition}

         \begin{proof}  The proof is a simple specialization
of the general
         construction sketched in the paragraph preceding the statement
         of the proposition.
          \end{proof}

           It will be convenient to work also with the analytic Toeplitz
           operators acting between $H^{2}(E,\sigma)$-spaces of different
           multiplicity.  For this purpose, we suppose that $\cU$ and $\cY$
           are two additional auxiliary Hilbert spaces (to be thought of as
           an {\em input space} and {\em output space} respectively).  We
           consider higher multiplicity versions of $H^{2}(E,\sigma)$ by
           tensoring  with an auxiliary Hilbert space (which is
to be thought
           of as adding multiplicity):
           $$
           H^{2}_{\cU}(E, \sigma):
           = H^{2}(E,\sigma) \otimes_{\mathbb C} \cU, \qquad
           H^{2}_{\cY}(E, \sigma):
           = H^{2}(E, \sigma)\otimes_{\mathbb C} \cY.
           $$
           Here we view $\cU$ and $\cY$ as
           $(\mathbb{C},\mathbb{C})$-correspondences and
           apply the tensor-product construction of Definition
\ref{D:tensor} (2).
          The space $H^{2}_{\cU}(E,\sigma)$ then is a reproducing kernel
           $(\sigma(\cA)', \cL(\cE\otimes\cU))$-correspondence
           on ${\mathbb D}((E^{\sigma})^{*})$ where the point evaluation
           at a point $(\eta,b)\in{\mathbb
D}((E^{\sigma})^{*})\times\sigma(\cA)'$
           of a function $f^\wedge\otimes u\in H^{2}_{\cU}(E,
\sigma)$ (with
           $f^\wedge\in H^{2}(E, \sigma)$ and $u\in\cU$) is given by
           $(f^\wedge\otimes u)(\eta,b)=f^\wedge(\eta,b)\otimes
u\in\cE\otimes\cU$.
        Moreover, note that the left
           $\sigma(\cA)'$-action is given by $b\mapsto b\otimes I_\cU$.
           The completely
           positive kernel $\K_{(E,\sigma) \otimes \cU }$
           associated with it as in Theorem \ref{T:RKHilcor} is given by
           \begin{equation*} %\label{Szego-ker-U}
           \K_{(E,\sigma) \otimes \cU}(\eta, \zeta)[b] =
           \K_{E,\sigma}(\eta, \zeta)[b]\otimes I_\cU,
           \end{equation*}
           where $\K_{E,\sigma}$ denotes the kernel for
$H^2(E,\sigma)$ defined in
           Theorem \ref{T:H2Esigma}.
           Similar statements hold for $H^{2}_{\cY}(E, \sigma)$, where the
           analogous kernel is denoted by $\K_{(E,\sigma) \otimes \cY}$.

           We now define a higher-multiplicity version of the algebra
of analytic
           Toeplitz operators $H^{\infty}(E,\sigma)$ to be
           the linear space
           $$
           H^{\infty}_{\cL(\cU, \cY)}(E, \sigma) : =
           H^{\infty}(E, \sigma) \otimes \cL(\cU, \cY).
           $$
          This space consists of $\cL(\cE\otimes\cU,\cE \otimes\cY)$-valued
          functions on ${\mathbb D}((E^{\sigma})^{*})$, with
point evaluation of an
          element $S\otimes N\in H^{\infty}_{\cL(\cU, \cY)}(E, \sigma)=
           H^{\infty}(E, \sigma) \otimes \cL(\cU, \cY)$ at
           $\eta\in{\mathbb D}((E^{\sigma})^{*})$ given by
           $(S\otimes N)(\eta)=S(\eta)\otimes N$.
           Moreover, the functions in $H^{\infty}_{\cL(\cU,
\cY)}(E, \sigma)$
           define multiplication operators in
           $\cL(H^2_\cU(E,\sigma),H^2_\cY(E,\sigma))$, in the same way as
           $H^{\infty}(E, \sigma)$. For
           $S\otimes N\in H^{\infty}_{\cL(\cU, \cY)}(E, \sigma)$,
           $S\in H^{\infty}(E, \sigma)$ and $N\in\cL(\cU, \cY)$,
the multiplication
           operator $M_{S\otimes N}$ becomes $M_{S\otimes N}=M_S\otimes N$.

           In addition there are Fock space versions of all these spaces,
           namely
           \begin{align*}
          &  \cF^{2}_{\cU}(E, \sigma): = \cF^{2}(E, \sigma)
\otimes \cU,
           \qquad  \cF^{2}_{\cY}(E, \sigma): = \cF^{2}(E,
\sigma) \otimes
           \cY, \\
          & \cF^{\infty}_{\cL(\cU, \cY)}(E, \sigma) =
           \cF^{\infty}(E, \sigma) \otimes \cL(\cU, \cY).
          \end{align*}
           Point evaluation for elements in $\cF^{2}_{\cU}(E, \sigma)$
          and points in ${\mathbb D}((E^{\sigma})^{*})\times\sigma(\cA)'$
          (and similarly for elements in $\cF^{2}_{\cY}(E,
\sigma)$) is determined
          by attaching to $f\otimes u\in\cF^{2}_{\cU}(E, \sigma)$,
          $f\in\cF^{2}(E, \sigma)$ and $u\in\cU$, and
          $(\eta,b)\in{\mathbb
D}((E^{\sigma})^{*})\times\sigma(\cA)'$ the value
          $(f\otimes u)^\wedge(\eta,b)=f^\wedge(\eta,b)\otimes
u$, so that the map
          \[
          \Phi_\cU:f_u\mapsto f_u^\wedge\mbox{ for
}f_u\in\cF^{2}_{\cU}(E, \sigma)
          \]
          defines a coisometry from $\cF^2_\cU(E,\sigma)$ onto
$H^2_\cU(E,\sigma)$.
          The analogous coisometry for $\cF^2_\cY(E,\sigma)$ is
denoted by $\Phi_\cY$.
          Similarly we determine point evaluation for elements in
          $\cF^{\infty}_{\cL(\cU, \cY)}(E, \sigma)$ and points in
          ${\mathbb D}((E^{\sigma})^{*})$ by attaching to
          $R\otimes X \in\cF^{\infty}_{\cL(\cU, \cY)}(E, \sigma)$,
          $R\in\cF^{\infty}(E,\sigma)$ and $X \in\cL(\cU, \cY)$, and
          $\eta\in{\mathbb D}((E^{\sigma})^{*})$ the value
          \begin{equation}  \label{pteval-mult}
    (R\otimes X)^\wedge(\eta):=R^\wedge(\eta)\otimes X
          \in\cL(\cE\otimes\cU,\cE\otimes\cY).
          \end{equation}

          Then $H^{\infty}_{\cL(\cU, \cY)}(E, \sigma)$ is recovered as
          \[
          \{ R^\wedge\colon R\in\cF^{\infty}_{\cL(\cU, \cY)}(E, \sigma)\},
          \]
          where $R^\wedge$ should be interpreted as a function mapping
          ${\mathbb D}((E^{\sigma})^{*})$ into
          $\cL(\cE\otimes\cU, \cE \otimes\cY)$,
          while the space of multiplication operators in
          $\cL(H^2_\cU(E,\sigma),H^2_\cY(E,\sigma))$ defined by
          $H^{\infty}_{\cL(\cU, \cY)}(E, \sigma)$ is given by
          \[
          \{\Phi_\cY R\Phi_\cU^*\colon R\in\cF^{\infty}_{\cL(\cU,
\cY)}(E, \sigma)\}.
          \]
          In fact, it is easy to check that Proposition
\ref{P:Hinftycar} guarantees
          that for $S\in H^{\infty}_{\cL(\cU, \cY)}(E, \sigma)$ we have
           $M_S=\Phi_\cY R\Phi_\cU^*$ whenever
           $R\in\cF^{\infty}_{\cL(\cU, \cY)}(E, \sigma)$
satisfies $S=R^\wedge$, so
           that $\|M_S\|\leq\|R\|$, and that there exists a
          $R\in\cF^{\infty}_{\cL(\cU, \cY)}(E, \sigma)$ with
$S=R^\wedge$ and
          $\|M_S\|=\|R\|$.

           Alternatively,
           $\cF_{\cL(\cU, \cY)}^{\infty}(E,\sigma)$ can be characterized
           as
           bounded operators from
           $\cF^{2}_{\cU}(E, \sigma)$ to
           $\cF^{2}_{\cY}(E, \sigma)$ with
           block-matrix representation
           $$
         R = [R_{i,j}]_{i,j=0,1, \dots} \text{ with } R_{i,j}\colon
         E^{\otimes j} \otimes_{\sigma} \cE  \otimes  \cU \to
E^{\otimes
       i}
         \otimes_{\sigma} \cE \otimes \cY
           $$
           subject to
           \begin{align*}
          & R_{i,j}(I_{E^{\otimes j}} \otimes b \otimes I_{\cU}) =
           (I_{E^{\otimes i}} \otimes b \otimes I_{\cY}) R_{i,j} \text{
           for } b \in \sigma(\cA)', \\%\label{intertwine1} \\
         & R_{i+1, j+1} (I_{E^{\otimes j}} \otimes \eta^{*} \otimes
         I_{\cU})
         = (I_{E ^{\otimes i}} \otimes \eta^{*} \otimes
I_{\cY}) R_{i,j}
         \text{ for } \eta^{*} \in E^{\sigma}.  %\label{intertwine2}
           \end{align*}
         For such $R \in \cF^{\infty}_{\cL(\cU, \cY)}(E, \sigma)$
         point evaluation in $\eta \in {\mathbb D}((E^{\sigma})^{*})$
         can be written as
          $$
        R^{\wedge}(\eta) = \sum_{n=0}^{\infty} (\eta^{n}
        \otimes I_{\cY}) R_{n,0}.
          $$
           In addition it is routine to see that part (1) of Proposition
           \ref{P:multiplicative} can be extended to the
following statement:
           {\em if $S \in \cF^{\infty}_{\cL(\cU, \cY)}(E,\sigma)$ and
           $R \in
           \cF^{\infty}_{\cL(\cY, \cZ)}(E, \sigma)$, then $ R S \in
           \cF^{\infty}_{\cL(\cU, \cZ)}(E, \sigma)$ and}
           \begin{equation}  \label{ext-mult}
           (RS)^{\wedge}(\eta) = R^{\wedge}(\eta) S^{\wedge}(\eta).
        \end{equation}

        \begin{remark}  \label{R:pteval'}
        {\em
            Suppose that $R \in \cF^{\infty}_{\cL(\cU,
            \cY)}(E, \sigma)$ has the form
            \begin{equation} \label{special'}
            R = \widetilde R \otimes_{\sigma} I_{\cE} \otimes X
        \end{equation}
        where $\widetilde R \in \cF^{\infty}(E)$ and $X \in
        \cL(\cU, \cY)$.  In particular the point evaluation
        \eqref{pteval-mult} defines $R^{\wedge}(\eta) \in
        \cL(\cE \otimes \cU, \cE \otimes \cY)$ for each $\eta
        \in {\mathbb D}((E^{\sigma})^{*})$.  Suppose now that
        $\sigma' \colon \cA \to \cL(\cE')$ is another
        $*$-representation of $\cA$ and $\eta' \in {\mathbb
        D}((E^{\sigma'})^{*})$.  Then we may define a related
        function
        $ \eta' \mapsto R^{\wedge \prime}(\eta') \in \cL(\cE'
        \otimes \cU, \cE' \otimes \cY)$ by
        $$
        R^{\wedge \prime}(\eta') = (\widetilde R
        \otimes_{\sigma'}I_{\cE'}\otimes X)^{\wedge}(\eta').
        $$
        While not all elements $R$ of $\cF^{\infty}_{\cL(\cU,
        \cY)}(E, \sigma)$ are of the special form
        \eqref{special'}, finite linear combinations of
        elements of the special form  \eqref{special'} are
        weak-$*$ dense in $\cF_{\cL(\cU, \cY)}(E, \sigma)$.
        By using linearity and a limiting process, one can
        then make sense of $R^{\wedge \prime}(\eta') \in
        \cL(\cE' \otimes \cU, \cE' \otimes \cY)$ for any
        $\eta' \in {\mathbb D}((E^{\sigma'})^{*})$.  This
        fact will be useful for the formulation of condition
        (1$^{\prime \prime}$) in Theorem \ref{T:MS-Schur}
        below.
         }\end{remark}

        %%%%%%%%%%%%%%%%%%%%%%%%%%%%%%%%%%%%%%%%%%%%%%%
        %%%%%%%%%%%%%%%%%%%%%%%%%%%%%%%%%%%%%%%%%%%%%%%
           \section{The Schur class associated with $(E, \sigma)$}
           \setcounter{equation}{0}

           Given a correspondence-representation pair $(E, \sigma)$
           (where
           $\sigma \colon \cA \to \cL(\cE)$) along with
auxiliary Hilbert
           spaces $\cU$ and $\cY$, we define the associated Schur class
           $\cS_{E,\sigma}(\cU, \cY)$ by
           \begin{align}
               \cS_{E,\sigma}(\cU, \cY) = &
               \{ S \colon {\mathbb D}((E^{\sigma})^{*}) \to
               \cL(\cE\otimes \cU,\cE\otimes\cY)
               \colon \, S(\eta) =
               R^{\wedge}(\eta) \text{ for all } \eta \in {\mathbb
               D}((E^{\sigma})^{*}) \notag \\
               & \qquad \text{for some } R \in \cF^{\infty}_{\cL(\cU,
               \cY)}(E, \sigma) \text{ with } \| R \| \le 1 \}.
               \label{HinftyUY}
           \end{align}

           We have the following characterization of the Schur class
             $\cS_{E,\sigma}(\cU, \cY)$ analogous to the
characterization of the
             classical Schur class given in Theorem \ref{T:C}
             and to the multivariable extensions in Theorems
\ref{T:BTV} and
           \ref{T:NC1}.
           When specialized to the classical case (see Section
           \ref{S:classical} below), \eqref{HinftyUY} gives
the classical
           Schur class as defined in the Introduction, but
from a different
           point of view.
           Rather than simply holomorphic,
           contractive, $\cL(\cU, \cY)$-valued function on
the unit disk,
           \eqref{HinftyUY} asks us to think of such
functions as analytic
           functions $F(z) \sim \sum_{n=0}^{\infty} F_{n} z^{n}$ on
           ${\mathbb D}$ whose Taylor coefficients $\{
F_{n}\}_{n \in {\mathbb
           Z}_{+}}$ induce a Toeplitz matrix
           $$
            T_{F} = \begin{bmatrix} F_{0} & 0 & 0 & \dots  \\
F_{1} & F_{0} &
            0  & \dots \\ F_{2}&
            F_{1} & F_{0} & \\ \ddots & \ddots & \ddots &
\ddots \end{bmatrix}
           $$
           which defines a contraction operator from
$\ell^{2}_{\cU}({\mathbb
           Z}_{+})$ to
           $\ell^{2}_{\cY}({\mathbb Z}_{+})$.  Thus the label
           (1) in Theorem \ref{T:MS-Schur}, when specialized
to the classical
           case, corresponds to a somewhat different
statement than (1) in
           Theorem \ref{T:C}.  The other labels
(1$^{\prime}$), (1$^{\prime
           \prime}$), (2) and (3) in Theorem \ref{T:MS-Schur}
correspond
           exactly to the corresponding statements in
Theorems \ref{T:C},
           \ref{T:BTV} and \ref{T:NC1}.

           \begin{theorem}  \label{T:MS-Schur}
               Suppose that we are given a
correspondence-representation pair
               $(E, \sigma)$  (where $\sigma \colon \cA \to
\cL(\cE)$) along with
               auxiliary Hilbert spaces $\cU$ and $\cY$ and
an operator-valued
               function $
             S \colon {\mathbb D}((E^{\sigma})^{*}) \to
\cL(\cE \otimes
             \cU, \cE \otimes \cY)$.
              Then the following conditions are equivalent:

              \begin{enumerate}
              \item[(1)] $S \in \cS_{E, \sigma}(\cU, \cY)$,
i.e., there exists
              an $R \in \cF^{\infty}_{\cL(\cU, \cY)}(E,
\sigma)$ with $\| R \|
              \le 1$ such that $S(\eta) =  R^{\wedge}(\eta)$
for all $\eta
              \in {\mathbb D}((E^{\sigma})^{*})$.

              \item[(1$^{\prime}$)] The multiplication operator
              $$ M_{S} \colon f(\eta,b) \mapsto S(\eta) f(\eta,b)
              $$
              maps $H^{2}_{\cU}(E, \sigma)$ contractively into
              $H^{2}_{\cY}(E, \sigma)$.

              \item[(1$^{\prime \prime}$)]
              $S$ is such that $S(\eta) =
               R^{\wedge}(\eta)$ for all $\eta \in {\mathbb
              D}((E^{\sigma})^{*})$ for an $R \in
\cF^{\infty}_{\cL(\cU,
               \cY)}(E, \sigma)$ with the
               additional property: for each  representation $\sigma'
               \colon \cA \to \cL(\cE')$ and $\eta' \in {\mathbb
               D}((E^{\sigma'})^{*})$ it happens that
               $$
               \|R^{\wedge \prime}(\eta')\| \le 1,
               $$
               where $R^{\wedge \prime}(\eta')$ is defined as
               in Remark \ref{R:pteval'}.

        \item[(2)] The function  $\K_S:{\mathbb
D}((E^{\sigma})^{*})\times
              {\mathbb D}((E^{\sigma})^{*})
              \to\cL(\sigma(\cA)',\cL(\cE\otimes\cY))$   defined by
              \begin{equation*}  %\label{KSkernel}
              \K_{S}(\eta, \zeta)[b]: =
              \K_{(E, \sigma) \otimes \cY}(\eta, \zeta)[b] -
              S(\eta) \K_{(E, \sigma)\otimes \cU}(\eta, \zeta) [b]
              S(\zeta)^{*}
              \end{equation*}
              is completely positive, or more explicitly,
              there exists an auxiliary Hilbert space $\cH$, an
             operator-valued function $H \colon
             {\mathbb D}((E^{\sigma})^{*}) \to \cL(\cH, \cE
                 \otimes
             \cY)$ and a $*$-representation $\pi$ of
$\sigma(\cA)'$ on $\cH$
            so that
            \begin{equation} \label{Agler-decom}
            \left(\sum_{n=0}^{\infty}
\eta^{n}(I_{E^{\otimes n}} \otimes b)
            (\zeta^{n})^{*} \right) \otimes I_{\cY}   -
            S(\eta) \left[\left(  \sum_{n=0}^{\infty}
\eta^{n}(I_{E^{\otimes
            n}} \otimes b )
            (\zeta^{n})^{*} \right)
            \otimes I_{\cU} \right]   S(\zeta)^{*} =
            H(\eta) \pi(b) H(\zeta)^{*}
             \end{equation}
             for all $\eta, \zeta \in {\mathbb
D}((E^{\sigma})^{*})$ and $b
             \in \sigma(\cA)'$.

          \item[(3)] There exists an auxiliary Hilbert space $\cH$,
          a $*$-representation $\pi \colon \sigma(\cA)' \to
\cL(\cH)$, and
               a coisometric colligation
              \begin{equation}  \label{coiscol}
              {\mathbf U} =
              \begin{bmatrix} A & B \\ C & D \end{bmatrix} \colon
              \begin{bmatrix} \cH \\ \cE \otimes \cU
              \end{bmatrix} \to
              \begin{bmatrix} E^\sigma \otimes\cH \\
              \cE \otimes \cY \end{bmatrix}
             \end{equation}
             which is a $\sigma(\cA)'$-module map, i.e.,
             \begin{equation}  \label{U-modulemap2}
              \begin{bmatrix} A & B \\ C & D \end{bmatrix}
\begin{bmatrix}
             \pi(b) h \\ (b \otimes I_{\cU}) {\mathbf u}
\end{bmatrix} =
           \begin{bmatrix} (I_{E} \otimes b) \otimes I_{\cH}
& 0 \\ 0 & b \otimes
               I_{\cY} \end{bmatrix} \begin{bmatrix} A & B \\ C & D
           \end{bmatrix} \begin{bmatrix} h \\ {\mathbf u}
           \end{bmatrix}
           \end{equation}
           for $h \in {\mathcal H}$ and ${\mathbf u} \in \cE
\otimes \cU$,
           so that $S$ can be realized as
            \begin{equation} \label{cor-realization2}
             S(\eta)  = D + C (I - L_{\eta^{*}}^{*}  A)^{-1}
             L_{\eta^{*}}^{*} B.
             \end{equation}
          Here $L_{\eta^{*}} \colon \cH \to E^{\sigma}
\otimes \cH$ is given
          by
        \[
        L_{\eta^{*}} h = \eta^*\otimes h \mbox{ for each } h\in\cH.
        \]
        \end{enumerate}

        \end{theorem}

        \begin{proof}

        Both \textbf{(1) $\Longrightarrow$ (1$^{\prime}$)} and
        \textbf{(1$^{\prime}$) $\Longrightarrow$ (1)} follow
        immediately after extending
        Proposition \ref{P:Hinftycar} to the case with the
added multiplicity
        as mentioned at the end of Section \ref{S:Hardyspaces}.

            \textbf{(1) $\Longrightarrow$ (1$^{\prime
        \prime}$):}    Given $\eta' \in {\mathbb
        D}((E^{\sigma'})^{*})$ (so $\| \eta'\| < 1$), by the
dilation result in
        \cite[Theorem 2.13]{MS-Annalen} (see also \cite{MS-CanJ})
        we know that $\eta'$ has a dilation to an induced
        representation $\eta_{\text{\rm ind}} \colon
\cF^{\infty}(E) \to
        \cL(\cF^{2}(E, \sigma_{\text{\rm ind}}))$ associated
with a representation
        $\sigma_{\text{\rm ind}} \colon \cA \to
        \cL(\cE_{\text{\rm ind}})$.  As $R$ is contractive by
assumption, it then
        follows that $R^{\wedge \text{\rm
        ind}}(\eta_{\text{\rm ind}})$ is also contractive.
        Since $\eta_{\text{\rm ind}}$ is a
        dilation of $\eta$, we then also have
        \begin{align*}
           \| R^{\wedge \prime} (\eta')\|
           & = \| P_{\cE \otimes \cY}
            R^{\wedge \text{\rm ind}}(\eta_{\text{\rm
ind}})|_{\cE \otimes \cU}\| \\
           & \le \| R^{\wedge \text{\rm ind}}(\eta_{\text{\rm
ind}})\| = \|R\| \le 1
        \end{align*}
        and (1$^{\prime \prime}$) follows.

            \textbf{(1$^{\prime \prime}$) $\Longrightarrow$ (2):}
             This implication requires an adaptation of the GNS/HB
             construction to the setting of completely positive (rather
             than classical positive) kernels.  If
$\K(\omega',\omega)[a]$ is a
             completely positive kernel, then
             $$
               K((\omega',a'), (\omega,a)) =
\K(\omega,\omega')[a^{*}a']
             $$
             is a positive kernel in the classical sense on
$\Omega \times \cA$.
             In this way one can reduce to the classical
setting and adapt
             the GNS/HB construction in \cite{DMM} to the
situation here.
             We leave complete details for another occasion.

            \textbf{ (1$^{\prime \prime}$)  $\Longrightarrow$
            (1$^{\prime}$):}
            Assume that $R \in \cF^{\infty}_{\cL(\cU,
\cY)}(E, \sigma)$
            and that $S = R^{\wedge}$.  From Proposition
            \ref{P:induced} extended to the higher
multiplicity setting,
            we see that we recover $R$ via the
point-evaluation calculus as
            $$
            R = R^{\wedge}(\eta_{\text{\rm ind}}).
            $$
            Hence we also recover $R$ as the strong limit
            $$
             R = \lim_{r \uparrow 1} R^{\wedge}( r
\eta_{\text{\rm ind}}).
            $$
            The assumption (1$^{\prime \prime}$) tells us that
            $$
            \|R^{\wedge}( r \eta_{\text{\rm ind}})\| \le 1
               $$
               for each $r < 1$. Hence $\|R\| \le 1$.

            \textbf{(1$^{\prime}$) $\Longrightarrow$ (2)}
            Assume that $M_{S}$ is as in (1$^{\prime}$).
            From the definitions we see that
            $$ (b \cdot M_{S} f^{\wedge})(\eta', b') = S(\eta')
            f^{\wedge}(\eta', b'b) = (M_{S}(b \cdot
f^{\wedge})) (\eta', b')
            $$
            and hence any multiplication operator $M_{S}$ is a
            $\sigma(\cA)'$-module map.
            The computation
            \begin{align*}
               & \langle M_{S} f, b \cdot (k_{E,\sigma; \zeta} \otimes
            I_{\cY}) (e \otimes y) \rangle_{H^{2}_{\cY}(E,\sigma)}
            = \langle b^{*} \cdot M_{S}f, (k_{E,\sigma; \zeta}
            \otimes I_{\cY}) (e \otimes y) \rangle_{H^{2}_{\cY}(E,
            \sigma)} \\
            & \qquad = \langle M_{S}( b^{*} \cdot f),
(k_{E, \sigma; \zeta}
            \otimes I_{\cY}) (e \otimes y)
            \rangle_{H^{2}_{\cY}(E,\sigma)}\\
            & \qquad = \langle S(\zeta) (b^{*} \cdot
f)(\zeta), e \otimes y
            \rangle_{\cE \otimes \cY} \\
            &\qquad  = \langle b^{*} \cdot f,
(k_{E,\sigma; \zeta}\otimes
            I_{\cU}) S(\zeta)^{*} (e \otimes y)
            \rangle_{H^{2}_{\cU}(E, \sigma)} \\
            & \qquad = \langle f, b \cdot (k_{E,\sigma;
\zeta} \otimes
            I_{\cU}) S(\zeta)^{*} (e \otimes y)
\rangle_{H^{}_{\cU}(E,
            \sigma)}
             \end{align*}
             shows that
           \begin{equation}  \label{MSstar}
           M_{S}^{*} \colon b \cdot (k_{E,\sigma;
\zeta}\otimes I_{\cY}) (e
           \otimes y) \mapsto
           b \cdot (k_{E,\sigma; \zeta} \otimes I_{\cU})
S(\zeta)^{*} (e
           \otimes y).
           \end{equation}
           Since $\| M_{S}\| \le 1$ by assumption, for any
finite collection
        of
           $b_{j} \in \sigma(\cA')$, $\zeta_{j} \in {\mathbb
           D}((E^{\sigma})^{*})$ and $e_{j} \otimes y_{j} \in
\cE \otimes \cY$
           ($j = 1, \dots, N$), we have
           \begin{equation}  \label{haveMSstar}
              \left\| \sum_{j=1}^{n} b_{j} \cdot
(k_{E,\sigma; \zeta_{j}}
              \otimes I_{\cY}) (e_{j} \otimes y_{j}) \right\|^{2} -
              \left\| M_{S}^{*} \, \sum_{j=1}^{n} b_{j} \cdot
(k_{E,\sigma;
        \zeta_{j}}
              \otimes I_{\cY}) (e_{j} \otimes y_{j})
\right\|^{2} \ge 0.
           \end{equation}

           Expanding out inner products and using
\eqref{MSstar} and the
           basic general identities
           \begin{align*}
              & \langle b' \cdot (k_{e, \sigma; \zeta}
\otimes I_{\cY}) (e'
        \otimes y'), b \cdot (k_{E, \sigma; \eta} \otimes
I_{\cY}) (e \otimes
               y) \rangle_{H^{2}_{\cY}(E, \sigma)}
              = \notag \\
              & \qquad \qquad
           \langle \K_{(E,\sigma) \otimes \cY}(\eta, \zeta)
[b^{*}b'](e'
               \otimes y'), e\otimes y \rangle_{\cE \otimes
\cY}, \notag \\
               & \langle b' \cdot (k_{E,\sigma, \zeta} \otimes I_{\cU})
               S(\zeta)^{*}(e' \otimes y'),
               b \cdot (k_{E, \sigma,\eta} \otimes I_{\cU})
S(\eta)^{*} (e
               \otimes y) \rangle_{H^{2}_{\cU}(E)} =  \notag \\
               & \qquad \qquad \langle S(\eta) \K_{(E,
\sigma) \otimes \cU}(\eta, \zeta)
               [b^{*}b'] S(\zeta)^{*} (e' \otimes y'), e \otimes y
               \rangle_{\cE \times \cY}
               %\label{basic}
            \end{align*}
            we see that the left hand side of
\eqref{haveMSstar} is equal to
            $$
            \sum_{i,j=1}^{N} \langle \K_{S}(\zeta_{i},
\zeta_{j})[b_{i}^{*}
            b_{j}] (e_{j} \otimes y_{j}), \, e_{i} \otimes
y_{i} \rangle_{\cE
            \otimes \cY}
            $$
            and we conclude that $\K_{S}$ is a completely
positive kernel
            as wanted.  The characterization given in
\eqref{Agler-decom}
            follows from part (2) of Theorem \ref{T:RKHilcor}.

            \textbf{(2) $\Longrightarrow$ (3):}  The argument
here is an
            adaptation of the proof of Theorem 3.5 in
\cite{MS-Schur} to our
            setting.
            Assume that (2) holds. By Remark \ref{R:etabmu},
the equality
            \eqref{Agler-decom} can be rewritten as
            $$
           ( I - \theta_{\eta, \zeta^{*}})^{-1}(b) \otimes I_{\cY}
            - S(\eta) [ ( I - \theta_{\eta,
            \zeta^{*}})^{-1}(b) \otimes I_{\cU} ] S(\zeta)^{*} =
            H(\eta) \pi(b) H(\zeta)^{*}.
            $$
            Replace $b$ by $[I - \theta_{\eta, \zeta^{*}}](b) = b -
            \theta_{\eta, \zeta^{*}}(b)$ to rewrite this last
expression as an
            {\em Agler decomposition} (see \cite{agler-hellinger})
            \begin{equation} \label{Agler-decom2}
             b \otimes I_{\cY} - S(\eta)  (b \otimes I_{\cU})
S(\zeta)^{*} =
             H(\eta) \pi\left( b - \eta (I_{E} \otimes b)
\zeta^{*}\right)
             H(\zeta)^{*}.
            \end{equation}
            Rearranging and conjugating by two generic
vectors ${\mathbf y}$
        and ${\mathbf
            y}'$ in $\cE \otimes \cY$ then gives us
           \begin{equation}  \label{Agler-decom3}
              {\mathbf y}^{*} H(\eta) \pi(b) H(\zeta)^{*}
{\mathbf y}' +
              {\mathbf y}^{*}S(\eta) (b
               \otimes I_{\cU}) S(\zeta)^{*} {\mathbf y}' =
              {\mathbf y}^{*} H(\eta) \pi( \eta (I_{E} \otimes b)
               \zeta^{*}) H(\zeta)^{*}{\mathbf y}' + {\mathbf y}^{*} (b
               \otimes I_{\cY}) {\mathbf y}'.
           \end{equation}
           From Remark \ref{R:etabmu} we know that $E^{\sigma}$ is a
           $\sigma(\cA)'$-correspondence.  We may also view
the Hilbert space
           $\cH$ as a $(\sigma(\cA)', {\mathbb
C})$-correspondence with the
           left $\sigma(\cA)'$-action given by the
representation $\pi$.  We
           may then form the tensor-product $(\sigma(\cA)', {\mathbb
           C})$-correspondence $E^{\sigma} \otimes \cH$ as in
Definition
           \ref{D:tensor}.  Explicitly, the ${\mathbb
C}$-valued inner product
           on $E^{\sigma} \otimes \cH$ is given by
           \begin{equation*}  %\label{EsigmaHIP}
               \langle \mu \otimes h,\, \mu' \otimes h'
\rangle_{E^{\sigma}
               \otimes \cH} = \langle \pi(\mu^{\prime *}\mu)
h, h' \rangle_{\cH}
               = h^{\prime *} \pi(\mu^{\prime *} \mu) h.
           \end{equation*}
           It follows that the first term on the right-hand side of
           \eqref{Agler-decom3} can be written as
           \begin{equation} \label{term3}
           {\mathbf y}^{*} H(\eta) \pi(\eta (I_{E} \otimes b)
\zeta^{*})
           H(\zeta)^{*} {\mathbf y}' =
           \langle (I_{E} \otimes b) \zeta^{*} \otimes
H(\zeta)^{*} {\mathbf y}',
           \eta^{*} \otimes H(\eta)^{*} {\mathbf y} \rangle_{E^{\sigma}
           \otimes \cH}.
           \end{equation}
           If we replace $b$ with $b^{*}b'$ (where $b,b'$ are
two elements of
           $\sigma(\cA)'$), use \eqref{term3} and do some
           rearranging, we see that
            the equality \eqref{Agler-decom3} can be
expressed in terms of
           inner products
           \begin{align}
             &  \langle \pi(b') H(\zeta)^{*}{\mathbf y}',
\pi(b) H(\eta)^{*} {\mathbf
               y} \rangle_{\cH} + \langle (b' \otimes
I_{\cU}) S(\zeta)^{*}
               {\mathbf y}', (b \otimes I_{\cU}) S(\eta)^{*}
{\mathbf y}
        \rangle_{\cE \otimes \cU}
               \notag \\
             & \qquad \qquad
             = \langle (I_{E }\otimes b') \zeta^{*} \otimes
H(\zeta)^{*} {\mathbf
             y}', (I_{E} \otimes b) \eta^{*} \otimes
H(\eta)^{*} {\mathbf y}
        \rangle_{E^{\sigma}
             \otimes \cH} + \langle (b' \otimes I_{\cY}) {\mathbf y}',
             (b \otimes I_{\cY}) {\mathbf y} \rangle_{\cE \otimes \cY}.
            \label{IPiden}
            \end{align}

          Introduce subspaces
             \begin{align*}
             {\mathcal D}_{V} & = \overline{\operatorname{span}}
             \left\{ \begin{bmatrix} (I_{E} \otimes  b) \eta^{*}
             \otimes H(\eta)^{*} {\mathbf y} \\ (b \otimes
I_{\cY}) {\mathbf y}
        \end{bmatrix}  \colon {\mathbf y} \in \cE \otimes \cY, \eta \in
        {\mathbb D}((E^{\sigma})^{*}), b \in \sigma(\cA)'
\right\}  \subset
        \begin{bmatrix}
        E^{\sigma} \otimes \cH \\ \cE \otimes \cY
\end{bmatrix}  \notag \\
          {\mathcal R}_{V} & = \overline{\operatorname{span}} \left\{
          \begin{bmatrix} \pi(b) H(\zeta)^{*} {\mathbf y} \\ (b \otimes
              I_{\cU}) S(\eta)^{*} {\mathbf y} \end{bmatrix}  \colon
              {\mathbf y} \in \cE \otimes \cY, \eta \in {\mathbb D}((
              E^{\sigma})^{*}), b \in \sigma(\cA)' \right\} \subset
              \begin{bmatrix} \cH \\ \cE \otimes \cU \end{bmatrix}.
              %\label{subspaces}
             \end{align*}
            Note that both ${\mathcal D}_{V}$ and ${\mathcal
R}_{V}$ are
           invariant under the left action of $\sigma(\cA)'$
           on $(E^{\sigma} \otimes \cH) \oplus  (\cE \otimes
           \cY)$ and on $\cH \oplus (\cE \otimes \cU)$ respectively,
           i.e. ${\mathcal D}_{V}$ and
           ${\mathcal R}_{V}$ are $\sigma(\cA)'$-submodules of
           $(E^{\sigma} \otimes \cH) \oplus  (\cE \otimes
\cY)$ and $\cH
             \oplus (\cE \otimes \cU)$ respectively.
          The import of \eqref{IPiden} is that the formula
          \begin{equation}  \label{DV}
              V \colon
              \begin{bmatrix} (I_{E} \otimes  b)
              \eta^{*} \otimes H(\eta)^{*} {\mathbf y} \\ (b
\otimes I_{\cY})
              {\mathbf y}
              \end{bmatrix}  \mapsto
              \begin{bmatrix} \pi(b) H(\zeta)^{*}{\mathbf y}
\\ (b \otimes
               I_{\cU}) S(\eta)^{*} {\mathbf y} \end{bmatrix}
           \end{equation}
           extends by linearity and continuity to a
well-defined unitary
           operator from ${\mathcal D}_{V}$ onto ${\mathcal R}_{V}$.
           One easily checks that
           $$
             V (b \cdot d) = b \cdot Vd
           $$
           for $b \in \sigma(\cA)'$ and $d \in {\mathcal D}_{V}$.

           By restricting in (\ref{DV}) to $b=I_\cE\in\sigma(\cA)'$ and
           $\eta=0\in{\mathbb D}((E^{\sigma})^{*})$ we see
that $\{0\} \oplus
           (\cE \otimes \cY) \subset\cD_V$.
           In particular
           \[
           \cX:=((E^\sigma\otimes\cH)\oplus(\cE\otimes\cY))\ominus\cD_V
           \subset (E^\sigma\otimes\cH) \oplus \{0\}.
           \]
           Moreover, because $\cD_V$ is invariant under the left
        $\sigma(\cA)'$-action we can
           see $\cX$ as a $(\sigma(\cA)',\C)$-correspondence,
where the left action is
           obtained by restricting the left action on
$E^\sigma\otimes\cH$ to $\cX$.
           Hence we can form the $(\sigma(\cA)',\C)$-correspondence
           $\cK=\cH\oplus(\cF^2(E^\sigma)\otimes\cX)$. Note that
           \[
E^\sigma\otimes\cK=E^\sigma\otimes(\cH\oplus(\cF^2(E^\sigma)\otimes\cX))
           =(E^\sigma\otimes\cH)\oplus
(E^\sigma\otimes\cF^2(E^\sigma)\otimes\cX).
           \]

           So we can define an operator ${\mathbf U}$ from
$\cK\oplus(\cE\otimes\cU)$ to
           $(E^\sigma\otimes\cK)\oplus(\cE\otimes\cY)$ via
           \begin{equation}\label{dilation}
           {\mathbf U}^*=\begin{bmatrix}VP_{{\mathcal
D}_{V}}&0&0&\cdots\\
                 P_\cX&0&0& \cdots \\

    0&I_{E^\sigma\otimes\cX}&0&\cdots\\0&0&I_{(E^\sigma)^{\otimes2}\otimes\cX}&\ddots\\
           \vdots&\vdots&\ddots&\ddots
           \end{bmatrix}
        :\begin{bmatrix}(E^\sigma\otimes\cH)\oplus(\cE\otimes\cY)\\E^\sigma\otimes\cX\\
           (E^\sigma)^{\otimes2}\otimes\cX\\\vdots\\\vdots\end{bmatrix}
\to\begin{bmatrix}\cH\oplus(\cE\otimes\cU)\\\cX\\E^\sigma\otimes\cX\\
           (E^\sigma)^{\otimes2}\otimes\cX\\\vdots\end{bmatrix}.
           \end{equation}
           Here $P_{{\mathcal D}_{V}}$ and $P_\cX$ stand for
the projections onto
           ${\mathcal D}_{V}$ and $\cX$ respectively. One easily
checks that
           ${\mathbf U}^*$ is an isometric
$\sigma(\cA)'$-module map. In other words,
           ${\mathbf U}$ is a coisometry, and a
$\sigma(\cA)'$-module map.
           The construction in (\ref{dilation}) is closely
related to the
        dilation result in
           \cite{MS-JFA}; see also Section 3 in
\cite{MS-models} for more details.

           Next we decompose ${\mathbf U}$
           as follows:
           \[
           {\mathbf U}=\begin{bmatrix}A&B\\C&D\end{bmatrix}
           :\begin{bmatrix}\cH \\\cE\otimes\cU\end{bmatrix}
\to\begin{bmatrix}E^\sigma\otimes\cH\\\cE\otimes\cY\end{bmatrix}.
           \]
           The definition of $V$ and the construction of
${\mathbf U}$ imply that
           \[
           \begin{bmatrix}A^*&C^*\\B^*&D^*\end{bmatrix}
           \begin{bmatrix} ((I_\cE\otimes b)\eta^*)\otimes
H(\eta)^* {\mathbf y}\\
               (b\otimes I_\cY) {\mathbf y} \end{bmatrix}
               =\begin{bmatrix} \pi(b)H(\eta)^* {\mathbf y} \\
               (b\otimes I_\cU)S(\eta)^* {\mathbf y}  \end{bmatrix}.
           \]
           By specifying this for $b=I_\cE$ and observing that
           \[
           \eta^*\otimes H(\eta)^* {\mathbf y} =L_{\eta^{*}} H(\eta)^*
           {\mathbf y},
           \]
           we get
           \begin{equation}\label{eqns}
           A^*L_{\eta^{*}} H(\eta)^*+C^*=H(\eta)^*\quad\mbox{and}\quad
           B^*L_{\eta^{*}} H(\eta)^*+D^*=S(\eta)^*.
           \end{equation}
           Moreover, for $h \in \cH$ we have
           \begin{eqnarray*}
           \|L_{\eta^{*}} h\|^2&=&\|\eta^*\otimes h\|^2
           =\langle \eta^*\otimes h,\eta^*\otimes h  \rangle
           =\langle \pi(\eta\eta^*)\otimes h, h \rangle
           =\|\pi(\eta\eta^*)^{\frac{1}{2}}h\|^2\\
           &&\leq\|\pi(\eta\eta^*)^{\frac{1}{2}}\|^2\|h\|^2
           \leq\|(\eta\eta^*)^{\frac{1}{2}}\|^2\|h\|^2
           =\|\eta\|^2\|h\|^2.
           \end{eqnarray*}
           This proves that $\|L_\eta\|\leq\|\eta\|<1$.
           Hence $I -  A^{*}L_{\eta}$ is invertible and
(\ref{eqns}) shows that
           \[
           H(\eta)^*=(I_\cK-A^*L_\eta)^{-1}C^*,
           \]
           and thus,
           \[
           S(\eta)^*=D^*+B^*L_\eta(I_\cK-A^*L_\eta)^{-1}C^*.
           \]
           By taking adjoints we arrive at (\ref{cor-realization2}).

        \textbf{(3) $\Longrightarrow$ (2):}
        Assume that (3) holds. We prove that $\K_S$ admits an
Agler decomposition
        as in \eqref{Agler-decom2} with
        $H(\eta):=C(I-L_{\eta^{*}}^* A)^{-1}$.
        That this is equivalent to the complete positivity of
the kernel
        $\K_{S}$ can be seen via the change of variable used in the
        derivation of \eqref{Agler-decom2}.
        The fact that ${\mathbf U}$ is a coisometric
$\sigma(\cA)'$-module
        map can also be written as
        \begin{align*}
        D(b\otimes I_\cY)D^*+C\pi(b)C^*=b\otimes I_\cY, \quad
        &\quad B(b\otimes I_\cU)D^*=-A\pi(b)C^*, \\
        A\pi(b)A^*+B(b\otimes I_\cU)B^*=(I_E\otimes b)\otimes
I_\cK, \quad
        &\quad D(b\otimes I_\cU)B^*=-C\pi(b)A^*.
        \end{align*}
        Note that
        \[
        H(\eta)=C(I-L_{\eta^{*}}^*A)^{-1}=C+C(I-L_{\eta^{*}}^*A)^{-1}
        L_{\eta^{*}}^*A=C+H(\eta)L_{\eta^{*}}^*A,
        \]
        and
        \[
        S(\eta)=D+H(\eta)L_{\eta^{*}}^{*}B.
        \]
        Hence
        \begin{eqnarray*}
        H(\eta)\pi(b)H(\zeta)^*
    &=&(C+H(\eta)L_{\eta^{*}}^{*}A)\pi(b)(C^*+A^*L_{\zeta^{*}}
H(\zeta)^*)\\
        &=&C\pi(b)C^*+C\pi(b)A^* L_{\zeta^{*}} H(\zeta)^*+
        H(\eta)L_{\eta^{*}}^*A\pi(b)C^*\\
&&        +H(\eta)L_{\eta^{*}}^*A\pi(b)A^*L_{\zeta^{*}} H(\zeta)^*\\
        &=&b\otimes I_\cY-D(b\otimes I_\cU)D^*-D(b\otimes
        I_\cU)B^*L_{\zeta^{*}} H(\zeta)^*\\
        &&-H(\eta)L_{\eta^{*}}^*B(b\otimes I_\cU)D^*
        -H(\eta)L_{\eta^{*}}^*B(b\otimes
I_\cU)B^*L_{\zeta^{*}} H(\zeta)^*\\
        &&+H(\eta)L_{\eta^{*}}^*((I_E\otimes b)\otimes
I_\cH)L_{\zeta^{*}} H(\zeta)^*\\
        &=&b\otimes I_\cY-D(b\otimes I_\cU)S(\zeta)^*
        -H(\eta)L_{\eta^{*}}^*B(b\otimes I_\cU)S(\zeta)^*\\
        &&+H(\eta)\pi(\eta(I_E\otimes b)\zeta^*)H(\zeta)^*\\
        &=&b\otimes I_\cY-S(\eta)(b\otimes I_\cU)S(\zeta)^*
        +H(\eta)\pi(\eta(I_E\otimes b)\zeta^*)H(\zeta)^*.
        \end{eqnarray*}
        In this way we have proved that (\ref{Agler-decom2}) holds.

          \textbf{(2) $\Longrightarrow$ (1$^{\prime}$):}
          Assume that (2) holds.  Consider the formula
          \begin{equation}  \label{formula}
              (M_{S})^{*} \colon b^{*} \cdot (k_{E,\sigma;
\zeta} \otimes
              I_{\cY}) {\mathbf y} \mapsto
              b^{*} \cdot (k_{E,\sigma; \zeta} \otimes
I_{\cU}) S(\zeta)^{*}
              {\mathbf y}
           \end{equation}
           for $b \in \sigma(\cA)'$, $\zeta \in {\mathbb
D}((E^{\sigma})^{*})$
          and ${\mathbf y} \in \cE \otimes \cY$. Then the
complete positivity
           of the kernel $\K_{S}$ is exactly what is needed
to see that the
           formula \eqref{formula} can be extended by
linearity and continuity
           to define a contraction operator $(M_{S})^{*}$
           from $H^{2}_{\cY}(E, \sigma)$ into $H^{2}_{\cU}(E,
\sigma)$ which
           is also a $\sigma(\cA)'$-module map:
           \begin{equation*} %\label{MS*modmap}
           b^{*} \cdot (M_{S}^{*}f) = M_{S}^{*} (b^{*} \cdot
f) \text{ for all
           } b \in \sigma(\cA)' \text{ and } f \in
H^{2}_{\cY}(E,\sigma).
           \end{equation*}
           Here we are using that the span of the collection of kernel
           functions
           $$ \{ b^{*} \cdot (k_{E, \sigma; \zeta} \otimes
           I_{\cY}) {\mathbf y}
            \colon b \in \sigma(\cA)',\, \zeta \in {\mathbb
           D}((E^{\sigma})^{*}), {\mathbf y} \in \cE \otimes \cY \}
           $$
           is dense in $H^{2}_{\cY}(E, \sigma)$.  Then the computation
           \begin{align*}
               \langle (M_{S}f)(\zeta, b), {\mathbf y}
\rangle_{\cE \otimes
               \cY} & =
               \langle M_{S} f, b^{*} \cdot (k_{E,\sigma;
\zeta} \otimes
               I_{\cY}) {\mathbf y} \rangle_{H^{2}_{\cY}(E, \sigma)} \\
               & = \langle f, M_{S}^{*} (k_{E,\sigma; \zeta} \otimes
               I_{\cY}) {\mathbf y} \rangle_{H^{2}_{\cU}(E, \sigma)} \\
               & = \langle f, b^{*} \cdot (k_{E, \sigma; \zeta} \otimes
               I_{\cY}) S(\zeta)^{*} {\mathbf y}
\rangle_{H^{2}_{\cU}(E,
               \sigma)} \\
               & = \langle f(\zeta, b), S(\zeta)^{*} {\mathbf
y} \rangle_{\cE
               \otimes \cU} \\
               & = \langle S(\zeta) f(\zeta,b), {\mathbf y}
\rangle_{\cE
               \otimes \cY}
            \end{align*}
            shows that $M_{S}$ is indeed the operator of
multiplication by $S$.
        \end{proof}

        \section{Examples}  \label{S:Examples}

        In this section we illustrate the general theory for some more
        concrete special cases.  For simplicity we consider
here only examples of the theory developed in
        Sections 3, 4 and 5 with $\cU = \cY = {\mathbb C}$.
Unlike what one might expect, this does not lead
        to scalar versions of the results discussed in
Sections 1 and 2, but rather to square versions, that is,
        we regain Theorems 1.1, 2.1 and 2.3 for the case
$\cU=\cY$, but not necessarily equal to ${\mathbb C}$.

        \subsection{The classical case}  \label{S:classical}

        In this example, we take $\cA=\cL(\cG)$ for a given
Hilbert space
          $\cG$. Let $E=\cL(\cG)$ viewed a correspondence over
itself in the
          standard way:
          \begin{align*}
         & a \cdot \xi = a \xi , \quad \xi \cdot a' =  \xi a'
\text{ (the operator
          multiplication in $\cL(\cG)$) for $a,a' \in \cA,\, \xi
\in E$}, \\
         & \langle \xi', \xi \rangle = \xi^{*} \xi'\text{ (the
          $\cL(\cG)$-inner product
        when considered as a correspondence over itself) for
$\xi', \xi \in E$.}
          \end{align*}
          Since
          $$
        \xi_{n} \otimes \xi_{n-1} \otimes \cdots
        \otimes \xi_{1} = 1\otimes
          1 \otimes \cdots \otimes 1\otimes  \xi_{n} \cdots \xi_{1},
          $$
          we can identify $E^{\otimes n}$ with $E=\cL(\cG)$ and
then the Fock
          space $\cF^{2}(E)$ has the form
          $$
          \cF^{2}(E) = \oplus_{n=0}^{\infty} E^{\otimes n} =
          \ell^2_{\cL(\cG)}(\mathbb{Z}_{+}).
          $$
          The abstract analytic Toeplitz algebra $\cF^{\infty}(E)$ is the
          collection of all
          lower  triangular Toeplitz
          matrices with $\cL(\cG)$-block entries acting as
bounded operators on
          $\ell^2_{\cL(\cG)}(\mathbb{Z}_{+})$.

          Now suppose we are given a Hilbert space $\cE_{0} $, let $\cE =
          \mathcal{G} \otimes \cE_{0}$ and  $\sigma$ be the representation
          of $\cA = \cL(\cG)$ on $\cL(\cE)$
          given by  $\sigma(a)= a \otimes I_{\cE_{0}}$.Then
          \begin{align*}
          \sigma(\cA)' &= \{ b \in \cL(\cE) \colon b \sigma(a) =
\sigma(a) b
          \text{ for all }  a \in \cA\}\\
          & =\{ b \in \cL(\cE) \colon b ( a^{0} \otimes
I_{\cE_{0}}) =  (a^{0} \otimes
          I_{\cE_{0}})b
          \text{ for all }  a^{0} \in \cL(\mathcal{G})\}\\
          &=\{I_{\mathcal{G}}\otimes b^{0} : b^{0} \in \cL(\cE_0)  \}
          \end{align*}
          and hence $\sigma(\cA)'$ can be identified with $\cL(\cE_{0})$.

          We next note that
          $$
          \cF^2(E,\sigma)= \cF^2(E) \otimes_{\sigma} \cE
          =l^2_{\cL(\cG)}(\mathbb{Z}_{+})\otimes_{\sigma}
          (\cG\otimes\cE_{0})=l^2_{\cG}(\mathbb{Z}_{+}) \otimes \cE_0 =
          \ell^{2}_{\cE}({\mathbb Z}_{+}).
          $$
          The representations $\varphi_{\infty,\sigma} \colon \cA
= \cL(\cG)
          \to \cL(\ell^{2}_{\cE}({\mathbb Z}_{+}))$ and
          $\iota_{\infty, \sigma} \colon \sigma(\cA)'  = \cL(\cE_{0}) \to
          \cL(\ell^{2}_{\cE}({\mathbb Z}_{+}))$ are given by
          \begin{align*}
           & \varphi_{\infty, \sigma}(a) =
           I_{\ell^{2}({\mathbb Z}^{+})} \otimes a \otimes I_{\cE_{0}},
           \qquad
            \iota_{\infty,\sigma}( b^{0}) =
           I_{\ell^{2}({\mathbb Z}^{+})} \otimes I_{\cG} \otimes b^{0}.
          \end{align*}

          We next compute
          \begin{align*}
          (E^{\sigma})^{*} &= \{ \eta \colon E \otimes_{\sigma} \cE \to \cE
          \colon
          \eta ( a \otimes I_{\cE}) = \sigma(a)\eta,  a
\in\cL(\cG) \} \\
          &=\{ \eta \colon \cL(\cG) \otimes_{\sigma}
\cG\otimes \cE_{0} \to
          \cG\otimes \cE_{0} \colon
          \eta ( a \otimes I_{\cE}) = ( a \otimes I_{\cE_{0}})\eta,  a
          \in\cL(\cG) \} \\
          &= \{ \eta \colon \cG\otimes \cE_{0} \to \cG\otimes
\cE_{0} \colon
          \eta ( a \otimes I_{\cE_{0}}) = ( a \otimes
I_{\cE_{0}})\eta,  a
          \in\cL(\cG) \} \\
          &= \{I_{\mathcal{G}}\otimes \eta^{0} : \eta^{0} \in
\cL(\cE_0)  \}.
          \end{align*}
          We conclude that $(E^{\sigma})^{*}$ can be identified
with $\cL(\cE_0)$.

          The creation operators and dual creation operators then
have the form
          \begin{align*}
           & T_{\xi, \sigma} ={\mathbf S} \otimes  \xi
\otimes I_{\cE_{0}}
           \text{ for } \xi \in \cA = \cL(\cG), \\
           & T^{d}_{\mu^{0},\sigma} = {\mathbf S} \otimes
I_{\cG} \otimes \mu^{0}
            \text{ for } \mu^{0} \in \cL(\cE^{0}) \cong E^{\sigma}
          \end{align*}
          where ${\mathbf S}$ is the standard shift operator on
          $\ell^{2}({\mathbb Z}_{+})$:
          $$
         {\mathbf S} \colon \{ c_{n}\}_{n \in {\mathbb
Z}_{+}} \mapsto \{
         c'_{n}\}_{n \in {\mathbb Z}_{+}} \text{ where } c'_{0} = 0, \,
         c'_{n} = c_{n-1} \text{ for } n \ge 1.
          $$
          Note that the commutativity properties laid out in Proposition
          \ref{P:Toeplitz} are now transparent for this example.

          Then, for $f = \oplus_{n=0}^{\infty} f_{n} \in
\cF^{2}(E, \sigma)$,
          the Fourier transform $\Phi f = f^{\wedge}$ is given by
          $$
          f^{\wedge}(\eta^{0}, b^{0}) = \sum_{n=0}^{\infty}(I_{\cG} \otimes
          (\eta^{0})^{n} b^{0}) f_{n} \in \cE
          $$
          for $\eta \in {\mathbb B}(\cL(\cE_{0}))$ (the open unit ball of
          $\cL(\cE_{0})$) and $b^{0} \in \cL(\cE_{0})$.
          One can check that $\Phi$ is injective. It follows that
$\Phi$ is a unitary
          transformation from $\ell^{2}_{\cE}({\mathbb Z}_{+})$
onto a Hilbert
          space $H^{2}(E, \sigma)$ of $\cE$-valued functions on
          ${\mathbb B}(\cL(\cE_{0})) \times \cL(\cE_{0})$
          carrying a $\cL(\cE_{0})$-representation:
          $$
        \pi_{H^{2}(E, \sigma)}(b^{0}) \colon
f^{\wedge}(\eta^{\prime 0}, b^{\prime 0})
        \mapsto f^{\wedge}(\eta^{\prime 0}, b^{\prime 0} b^{0}).
          $$
          In fact $f^{\wedge}(\eta^{0}, I_{\cE_{0}}) = 0$ for all
$\eta^{0} \in
          \cB(\cL(\cE_{0}))$ already forces $f$ to be zero in
          $\ell^{2}_{\cE}({\mathbb Z}_{+})$ so the function $f^{\wedge}$ is
          determined completely by its single-variable restriction $f^{\wedge
          1}: = f^{\wedge}(\cdot, I_{\cE_{0}})$ and one can work with the
          space $\widetilde H^{2}(E, \sigma) = \{ f^{\wedge 1} \colon f \in
          \ell^{2}_{\cE}\}$ instead.  One can identify $\widetilde H^{2}(E,
          \sigma)$ with functions of the form
          $$ g(\eta^{0}) = \sum_{n=0}^{\infty} (I_{\cG} \otimes
(\eta^{0})^{n})
          g_{n} \text{ where } \oplus_{n=0}^{\infty} g_{n} \in
          \ell^{2}_{\cE}({\mathbb Z}_{+}) \text{ with } \| g\|_{\widetilde
          H^{2}(E, \sigma)} = \| \oplus_{n=0}^{\infty} g_{n}
          \|_{\ell^{2}_{\cE}({\mathbb Z}_{+})}
          $$
          and with the $\sigma(\cA)' \cong \cL(\cE_{0})$-left
action given by
          $$
        (b^{0} \cdot g)(\eta) = \sum_{n=0}^{\infty} (I_{\cG} \otimes
        (\eta^{0})^{n}) (I \otimes b^{0}) g_{n} \text{ if }
g(\eta^{0}) =
        \sum_{n=0}^{\infty} (I_{\cG} \otimes (\eta^{0})^{n}) g_{n}.
          $$

          An element $S$ of $\cF^{\infty}(E, \sigma)$ is an operator on
          $\ell^{2}_{\cE}({\mathbb Z}_{+})$ having a
lower-triangular Toeplitz matrix
          representation
          $$  R = [ R_{i-j}]_{i,j = 0,1, \dots}
          $$
          where each $R_{n}$ is an operator on $\cE$ of the form $R_{n} =
          R_{n}^{0} \otimes I_{\cE_{0}}$   for an operator $R_{n}^{0} \in
          \cL(\cG)$ with $R_{n}^{0} = 0$ for $n<0$.  Given such an $R$, the
          associated $\cL(\cE)$-valued function $R^{\wedge}
\colon \cL(\cE_{0})
          \to \cL(\cE)$ is then given by
          $$
        R^{\wedge}(\eta^{0}) = \sum_{n=0}^{\infty} R_{n}^{0} \otimes
        (\eta^{0})^{n}.
          $$
          The Schur class $\cS(E, \sigma)$ for this case can be
identified with
          the set of functions $S \colon {\mathbb B}( \cL(\cE_{0}) )
          \to \cL(\cE)$ with a
          presentation of the form
          \begin{equation}  \label{PS-class}
        S(\eta^{0}) = \sum_{n=0}^{\infty} S^{0}_{n} \otimes
(\eta^{0})^{n}
          \end{equation}
          for  which the associated Toeplitz matrix
          $$  [S^{0}_{i-j}]_{i,j=0,1, \dots}
          $$
          defines a contraction operator on
$\ell^{2}_{\cG}({\mathbb Z}_{+})$.
          If we use the single-variable version $\widetilde
H^{2}(E, \sigma)$
          of the Hardy space, the condition in part (2) of Theorem
          \ref{T:MS-Schur} means not only that
          $$
        M_{S} \colon f^{\wedge 1}(\eta^{0}) \mapsto
S(\eta^{0}) f^{\wedge
        1}(\eta^{0})
          $$
          maps $\widetilde H^{2}(E, \sigma)$ contractively into $\widetilde
          H^{2}(E, \sigma)$, but also that $M_{S}$ is a
$\cL(\cE_{0})$-module map:
          $$
         M_{S} (b \cdot f^{\wedge 1}) = b \cdot M_{S} f^{\wedge 1}.
          $$

          The realization formula
          \eqref{U-modulemap2} and \eqref{cor-realization2} from
part (3) of
          Theorem \ref{T:MS-Schur} tells us that such
          functions $S$ are characterized by having a realization
of the form
          \begin{equation}  \label{real-class}
         S(\eta^{0}) = D + C (I - \pi(\eta^{0})A)^{-1} \pi(\eta^{0}) B
          \end{equation}
          where
          $$
         {\mathbf U} = \begin{bmatrix} A & B \\ C & D
\end{bmatrix} \colon
         \begin{bmatrix} \cH \\ \cE \end{bmatrix} \to
\begin{bmatrix} \cH \\
             \cE \end{bmatrix}
          $$
          is a unitary operator and $\pi $ is a $*$-representation of
          $\cL(\cE_{0})$ to $\cL(\cH)$ which is also a
$\cL(\cE_{0})$-module map:
          \begin{equation}  \label{U=mod-map}
         \begin{bmatrix} A & B \\ C & D \end{bmatrix} \begin{bmatrix}
             \pi(b^{0}) & 0 \\ 0 & I_{\cG} \otimes b^{0} \end{bmatrix}
           = \begin{bmatrix} \pi(b^{0}) & 0 \\ 0 & I_{\cG}
\otimes b^{0}
          \end{bmatrix} \begin{bmatrix} A & B \\ C & D \end{bmatrix}.
          \end{equation}
          Here we use that $E^{\sigma} \otimes_{\pi} \cH$ can be identified
          with $\cH$ since $(I_{\cG} \otimes (\eta^{0})^{*}) \otimes h =
          I_{\cE} \otimes \pi((\eta^{0})^{*}) h$.

          We note that is easy to see that a realization as in
          \eqref{real-class} implies that $S$ has a presentation
of the form
          \eqref{PS-class}.  Indeed, if $U$ is unitary and satisfies
          \eqref{U=mod-map}, since $A$ commutes with $\pi(\eta^{0})$ we see that
          $ (\pi(\eta^{0}) A)^{n} = A^{n} \pi(\eta^{0})^{n}$.
Hence expansion
          of the inverse in \eqref{real-class} as a geometric series and
          repeated usage of \eqref{U=mod-map} gives
          $$
        S(\eta^{0}) =   \sum_{n=0}^{\infty} S_{n} (I \otimes
        (\eta^{0})^{n})
          $$
          where
          $$
         S_{0} = D, \qquad S_{n} = C A^{n-1}B \text{ for } n \ge 1.
          $$
          Additional usage of \eqref{U=mod-map} gives us
          \begin{equation*}  %\label{Snmodmap}
        S_{n} (I \otimes \eta^{0}) = (I \otimes \eta^{0})
S_{n} \text{ for
        all } \eta^{0} \in \cL(\cE_{0})
          \end{equation*}
          from which we conclude that $S_{n}$ has the form $S_{n}
= S^{0}_{n}
          \otimes I_{\cE_{0}}$ for operators $S^{0}_{n}$ acting
on $\cG$, and
          hence $S(\eta^{0})$ has the form as in \eqref{PS-class}.

          Conversely, if $S \colon {\mathbb B}(\cL(\cE_0)) \to
          \cL(\cE)$ is of the form \eqref{PS-class},
          it follows that $S^0(\lambda) = \sum_{n=0}^{\infty}
S^0_n \lambda^n$
          is in the classical Schur class $\cS(\cG, \cG)$. By the
classical realization
          theorem we can write
          $$
           S^0(\lambda) = D^0 + \lambda C^0 (I - \lambda A^0)^{-1} B^0
           $$
           where
           $$
           {\mathbf U}^0 = \begin{bmatrix} A^0 & B^0 \\ C^0 & D^0
           \end{bmatrix}
           \colon \begin{bmatrix} \cH^0 \\ \cG \end{bmatrix} \to
           \begin{bmatrix} \cH^0 \\ \cG \end{bmatrix}
           $$
           is coisometric (or even unitary).  Then
           $$
           {\mathbf U} = {\mathbf U}^0 \otimes I_{\cE^0} =
\begin{bmatrix} A^0 \otimes
           I_{\cE_0} & B \otimes I_{\cE_0} \\ C^0 \otimes
I_{\cE_0} & D^0 \otimes I_{\cE_0}
       \end{bmatrix} \colon \begin{bmatrix} \cH \\ \cE
\end{bmatrix} \to \begin{bmatrix}
       \cH  \\  \cE \end{bmatrix}
           $$
           (where we set $\cH = \cH^0 \otimes \cE_0$)  with
           $$
         \pi(b) = I_{\cH^0} \otimes b \in \cL(\cH) \text{ for
} b \in \cL(\cE^0)
        $$
        provides a realization for $S$ as in
\eqref{real-class}. Thus the general
        theory provides a new kind of realization result, but
one can easily derive
        this result directly from the classical realization theorem.

          Two special cases of the above analysis are of interest.

          \begin{enumerate}
          \item If we take $\mathcal {G}=\cE, \cE_{0}=\mathbb{C}$ in the
          example, we have
          $\cF^2(E) = l^2_{\cL(\cE)}(\mathbb{Z}_{+})$ with
$\cF^{\infty}(E)$
          equal to the
          collection of all lower  triangular Toeplitz matrices with
          $\cL(\cE)$-block entries acting on
$l^2_{\cL(\cE)}(\mathbb{Z}_{+})$.
          In this case $\sigma(\cA)' =\mathbb{C}I_{\cE}$.
          and
          $(E^{\sigma})^*=I_{\cE}\otimes \mathbb{C} $ is isomorphic to
          $\mathbb{C}$; thus $\mathbb{D}((E^{\sigma})^{*})$ may
be identified
          with the open unit disk $\mathbb{D}$ of $\mathbb{C}$.
          Moreover $\cF(E) \otimes_{\sigma} \cE =l^2_{\cE}
(\mathbb{Z}_{+})$
          and for a given $\lambda \in \mathbb{D}$, we have the bounded
          point-evaluation:
          $$
          f=\oplus_{n=0}^{\infty}f_n \in\cF(E) \otimes_{\sigma}
          \cE=l^2_{\cE}(\mathbb{Z}_{+}) \to \widehat{f}(\lambda) =
          \sum_{n=0}^{\infty}f_n\lambda^n \in H^2_{\cE}(\mathbb{D}).
          $$
          Then $$H^{2}(E, \sigma) =H^2_{{\cE}}(\mathbb{D})$$ and
          $$H^{\infty}(E, \sigma) =
H^{\infty}_{{\cL(\cE})}(\mathbb{D})\otimes
          I_{\cE} = H^{\infty}_{\cL(\cE)}(\mathbb{D}).$$ Hence $S\in
          H^{\infty}(E, \sigma) $ means, for $\lambda \in \mathbb{D}$, that
          $S(\lambda)=\sum_{n=0}^{\infty}S_n{\lambda}^n$ with $S_n \in
          \cL(\cE)$. The operators in $H^{\infty}(E, \sigma) $ with norm at
          most equal to $1$ form the classical Schur class.  If
we apply the
          general theorem \ref{T:MS-Schur} for this case, we simply recover
          Theorem \ref{T:C} (where $\cU = \cY = \cG$).

          \medskip

          \item  If we take $\mathcal {G}=\mathbb{C},  \cE_{0}=\cE$, then
          $\cF^2(E) = \ell^2(\mathbb{Z}_{+})$, $\cF^{\infty}(E)$
is the collection
          of all lower  triangular Toeplitz matrices  acting on
          $\ell^2(\mathbb{Z}_{+})$. In this case $\sigma(\cA)'
=\cL(\cE)$ and
          \begin{equation*}
          (E^{\sigma})^{*} = \{ \eta \colon \mathbb{C}
\otimes_{\sigma} \cE \to
          \cE \colon
          \eta ( a \otimes I_{\cE}) = a\eta,  a \in\mathbb{C} \}.
          \end{equation*}
          Since $ \mathbb{C} \otimes_{\sigma} \cE$ can be
identified with $\cE$
          in the obvious way, $(E^{\sigma})^*$ amounts to
$\cL(\cE)$. We also have
          $\cF^2(E) \otimes_{\sigma} \cE =l^2_{\cE}(\mathbb{Z}_{+})$ for a
          given $\eta \in\mathbb{D}((E^{\sigma})^*)=
\mathbb{B}(\cL(\cE))$, i.e.,
          $\eta \in \cL(\cE)$ with $\|\eta\| < 1$, we have the
bounded point
          evaluation:
          $$
          f=\oplus_{n=0}^{\infty}f_n \in\cF(E) \otimes_{\sigma}
          \cE= \ell^2_{\cE}(\mathbb{Z}_{+}) \mapsto \left( (\eta,b)
      \in {\mathbb B}( \cL(\cE)) \times \cL(\cE) \mapsto
      f^{\wedge}(\eta, b) = \sum_{n=0}^{\infty} \eta^{n}b f_{n}
      \right).
      $$
    We may view this $H^{2}(E, \sigma)$ simply as functions of
    the form $\eta \mapsto \sum_{n=0}^{\infty} \eta^{n} f_{n}$
    with $\oplus_{n=0}^{\infty} f_{n} \in \ell^{2}_{\cE}$ which
    also carries an $\cL(\cE)$-action:
    $$
    f(\eta) \mapsto (b \cdot f)(\eta) = \sum_{n=0}^{\infty}
    \eta^{n} b f_{n} \text{ if } f(\eta) = \sum_{n=0}^{\infty}
    \eta^{n} f_{n}
    $$
    or
    $$
    b \cdot f = \sum_{n=0}^{\infty} {\mathbf S}^{n} b P_{\cE}
    {\mathbf S}^{*n} f
    $$
    where $P_{\cE}$ is the projection onto the constant functions
    and ${\mathbf S}$ is the operator-argument shift operator
    $({\mathbf S}f)(\eta) =
    \eta \cdot f(\eta)$.
    If we identify $S = S^{0} \otimes I_{\cE} \in H^{\infty}(E, \sigma)
       = H^{\infty} \otimes I_{\cE}$ with the scalar-valued function $S^{0}
     \in H^{\infty}({\mathbb D})$, then  the associated function with
    operator argument
    \begin{equation}  \label{S0form}
     \eta \mapsto
    S^{0}(\eta)=\sum_{n=0}^{\infty}\eta^n (S^{0}_n  \otimes I_{\cE})
    \in \cL(\cE)
    \end{equation}
    corresponds to the functional calculus for scalar holomorphic
    functions with operator argument usually defined via
       the holomorphic functional calculus (see e.g. \cite{Rudin}).
    The positivity of the kernel
    $$
    K_{S}(\eta, \zeta) = (I - \eta \zeta^{*})^{-1} - S(\zeta) (I
    - \eta \zeta^{*})^{-1} S(\eta)^{*}
    $$
    guarantees that the multiplication operator
    $$
    M_{S} \colon f(\eta) \mapsto S(\eta) f(\eta)
    $$
    is contractive on $H^{2}(E, \sigma)$ while complete
    positivity
    of the enlarged kernel
    $$
    \K_{S}(\eta, \zeta)[b] = \sum_{n=0}^{\infty} \eta^{n} b
    \zeta^{*n} - S(\eta) \left( \sum_{n=0}^{\infty} \eta^{n} b
    \zeta^{*n}\right) S(\zeta)^{*}
    $$
    guarantees in addition that $S$ has the form \eqref{S0form}  and that
    the associated multiplication operator
    $M_{S}$ commutes with the $\cL(\cE)$-action:
    $$
    b \cdot (M_{S}f) = M_{S}(b \cdot f) \text{ for } b \in
    \cL(\cE),\, f \in H^{2}(e, \sigma).
    $$
    The realization result (the equivalence of (1) and (3) in
          Theorem \ref{T:MS-Schur}) follows from the classical
          realization result for scalar-valued Schur-class
          functions in the same way as was explained above for
          the general case of this example.
        \end{enumerate}

          \subsection{Free semigroup algebras}  \label{S:free'}

          In this example, we take $\cA=\cL(\cG)$ for a given Hilbert space
          $\cG$ and $E$ to be the $d$-fold column space
          $\operatorname{col}_{j=1}^{d} \cL(\cG)$ over $\cL(\cG)$
viewed as a
          correspondence over $\cL(\cG)$ in the standard way (see
\cite{MS-JFA,
          RW}):
          \begin{align}
         & a\cdot \begin{bmatrix} \xi_{1} \\ \vdots \\
\xi_{d} \end{bmatrix}
          =\begin{bmatrix} a \xi_{1} \\ \vdots \\ a \xi_{d} \end{bmatrix},
          \qquad
          \begin{bmatrix} \xi_{1} \\ \vdots \\ \xi_{d} \end{bmatrix} \cdot
          a=\begin{bmatrix} \xi_{1} a \\ \vdots \\ \xi_{d} a
      \end{bmatrix},
         \qquad \left\langle\begin{bmatrix} \xi^{\prime}_{1}
\\ \vdots \\
         \xi^{\prime}_{d} \end{bmatrix}, \,
          \begin{bmatrix} \xi_{1} \\ \vdots \\ \xi_{d}
\end{bmatrix} \right\rangle =
          \sum_{j=1}^{d}  \xi_{j}^{*} \xi^{\prime}_{j} \notag \\
          & \qquad  \text{for } \xi = \begin{bmatrix} \xi_{1} \\ \vdots \\
          \xi_{d}
         \end{bmatrix}, \, \xi' = \begin{bmatrix}
\xi^{\prime}_{1} \\ \vdots
         \\ \xi^{\prime}_{d} \end{bmatrix} \in E \text{ and }
a \in \cL(\cG).
\label{free}
          \end{align}
          One can then identify $E^{\otimes n}$ with the column space
          $\oplus_{\alpha \colon
          |\alpha|=n} \cL(\cG) $ where $\alpha = i_n\cdots i_1$
is in the free
          semigroup $\cF_{d}$ with notation as in Subsection
\ref{S:free}. Then
          the associated Fock space is
          $$
          \mathcal{F}^2(E)= \oplus_{n=0}^{\infty} E^{\otimes n} =
           \oplus_{n=0}^{\infty}\left[ \oplus_{\alpha \in
\cF_{d}\colon |\alpha| =
          n} \cL(\cG)\right]
          $$
          can equally well be viewed as
          $$
        \mathcal F^{2}(E) = \oplus_{\alpha \in \cF_{d}} \cL(\cG) =:
        \ell^{2}_{\cL(\cG)}(\cF_{d}).
          $$
          Then the analytic Toeplitz algebra $\cF^{\infty}(E)$ can be
          identified as
          $$
          \cF^{\infty}(E) = \mathcal{L}_{d}\otimes\cL(\cG),
          $$
          where $\mathcal{L}_{d}$ is
          the free semigroup algebra discussed by Davidson and Pitts in
          \cite{DP} and is also the ultraweak closure of
Popescu's noncommutative
          disk algebra (see \cite{PopescuNCDA}).

          Just as done in the discussion of the classical case above, we
          now let $\cE_0$ be another Hilbert space and take $\cE$ to be
          $\cE = \cG \otimes \cE_0$. We consider the $*$-representation
          $\sigma$ of $\cA = \cL(\cG)$ to $\cL(\cE)$ given by
          $$
          \sigma(a) = a \otimes I_{\cE_0} \text{ for } a \in \cL(\cG).
          $$
          We compute $\sigma(\cA)'$ as follows:
          \begin{align*}
          \sigma(\cA)' &= \{ b \in \cL(\cE) \colon b \sigma(a) =
\sigma(a) e
          \text{ for all }  a \in \cA \}\\
          & =\{ b \in \cL(\cE) \colon b ( a \otimes I_{\cE_{0}})
=  (a \otimes
          I_{\cE_{0}}) b \text{ for all }  a \in \cL(\mathcal{G})\}\\
          &=\{I_{\mathcal{G}}\otimes b^{0} : b^{0} \in \cL(\cE_0)  \}.
          \end{align*}
          Hence $\sigma(\cA)'$ can be identified with $\cL(\cE_{0})$.

          We next note that
          $$
          E^{\otimes n} \otimes_{\sigma} \cE =
\operatorname{col}_{\alpha \colon
          |\alpha| = n} \cL(\cG) \otimes_{\sigma} \cE \cong
          \operatorname{col}_{\alpha \colon |\alpha| = n} \cE
          $$
          and hence we identify $\cF^{2}(E, \sigma) = \cF^{2}(E)
          \otimes_{\sigma} \cE$ as
          $$
        \cF^{2}(E, \sigma) = \oplus_{\alpha \in \cF_{d}} \cE =
        \ell^{2}_{\cE}(\cF_{d}).
          $$
          The representations $\varphi_{\infty,\sigma} \colon
\cL(\cG) = \cA
          \to \cL(\ell^{2}_{\cE}(\cF_{d}))$ and $\iota_{\infty,
\sigma} \colon
          \cL(\cE_{0}) \cong \sigma(\cA)' \to \cL(\ell^{2}_{\cE}(\cF_{d}))$
          can be seen to be given by
          \begin{align*}
          &  \varphi_{\infty,\sigma}(a) =
I_{\ell^{2}(\cF_{d})} \otimes a
           \otimes I_{\cE_{0}} \text{ for } a \in  \cL(\cG) = \cA, \\
          & \iota_{\infty, \sigma}(b^{0}) =
I_{\ell^{2}(\cF_{d})} \otimes
          I_{\cG} \otimes b^{0} \text{ for } b^{0} \in
\cL(\cE_{0}) \cong
          \sigma(\cA)'.
          \end{align*}

          We now observe that $E \otimes_{\sigma} \cE$ can be
identified with
          $\cE^d$ (the $d$-fold direct sum of $\cE$ with itself) under the
          identification
          $$
          \begin{bmatrix} \xi_{1} \\ \vdots \\ \xi_{d}
\end{bmatrix} \otimes e
          \cong  \begin{bmatrix} (\xi_{1} \otimes
I_{\cE_{0}})  e \\ \vdots \\
           (\xi_{d} \otimes I_{\cE_{0}} ) e \end{bmatrix}
          \text{ for }
        \xi_{1}, \dots, \xi_{d} \in \cL(\cG) \text{ and } e \in \cE.
          $$
          Then $\eta \in (E^{\sigma})^{*}$ means that $\eta$ is a block
          row-matrix $\eta = \begin{bmatrix} \eta_{1} & \cdots & \eta_{d}
          \end{bmatrix}$ mapping $E \otimes_{\sigma} \cE \cong
\cE^{d}$ into
          $\cE$ with the additional property that
          $$
          \begin{bmatrix} \eta_{1} & \cdots & \eta_{d} \end{bmatrix}
          \, \operatorname{diag}\, (a\otimes I_{\cE_0}) =
(a\otimes I_{\cE_0})
          \begin{bmatrix} \eta_{1} & \cdots & \eta_{d} \end{bmatrix}
           \text{ for all } a \in \cL(\cG).
        $$
          It follows that  $\eta_{j} (a\otimes I_{\cE_0}) =(a \otimes
          I_{\cE_0})\eta_j $ and hence that  $\eta_{j} = I_{\cG}\otimes
          \eta^{0}_{j}$ for some $\eta^{0}_{j} \in \cL(\cE_0)$
for $j = 1, \dots,
          d$ and we have an identification
          $$  (E^{\sigma})^{*} \cong \cL( \cE_{0}^{d}, \cE_{0}).
          $$

          One can check that the creation and dual creation operators are
          given by
          \begin{align*}
         &  T_{\xi, \sigma} = \sum_{j=1}^{d} {\mathbf S}_{j}
\otimes \xi_{j}
           \otimes I_{\cE_{0}} \text{ for } \xi
=\begin{bmatrix} \xi_{1} \\
           \vdots \\ \xi_{d} \end{bmatrix} \in \cL(\cG)^{d} = E, \\
           & T_{\mu^{0}, \sigma} = \sum_{j=1}^{d} {\mathbf
S}_{j} \otimes
           I_{\cG} \otimes \mu^{0}_{j} \text{ for }
           \mu^{0} = \begin{bmatrix} \mu^{0}_{1} \\ \vdots \\
\mu^{0}_{d}
          \end{bmatrix} \in \cL(\cE_{0}, \cE^{d}_{0}) \cong E^{\sigma}.
          \end{align*}

         For $\eta^{0} =  \begin{bmatrix} \eta^{0}_{1} & \cdots &\
         \eta^{0}_{d}
          \end{bmatrix} \in {\mathbb B}(\cL( \cE_{0}^{d}, \cE_{0})) \cong
          \mathbb{D}((E^{\sigma})^{*})$ and $b^{0} \in \cL(\cE_{0}) \cong
          \sigma(\cA)'$, we have the bounded point-evaluation operator
          on $\cF^{2}(E, \sigma) = \ell^{2}_{\cE}(\cF_{d})$:
          \begin{equation}  \label{freepteval}
          f=\oplus_{\alpha \in \cF_{d} } f_\alpha
          \mapsto  f^{\wedge}(\eta^{0}, b^{0}) :=
          \sum_{\alpha \in \mathbb{F}_{d}}
          (I_{\cG} \otimes (\eta^{0})^{\alpha} b^{0})f_{\alpha}
          \end{equation}
          where $(\eta^{0})^\alpha = \eta^{0}_{i_{N}} \cdots
\eta^{0}_{i_{1}}$ for
          $\alpha = i_N\cdots i_1 \in \cF_{d}$.  To continue a detailed
          analysis, we now consider in turn two divergent special cases.
          \medskip

          \textbf{Case 1:} {\em $\cE_{0} = {\mathbb C}$ so $\cE = \cG$:}
          In this case we identify ${\mathbb D}((E^{\sigma})^{*}) \cong
          {\mathbb B}(\cL(\cE_{0}^{d}, \cE_{0}))$ with the unit
ball in ${\mathbb C}^{d}$
          $$
         {\mathbb B}^{d}: = \left\{ \lambda = (\lambda_{1},
\dots, \lambda_{d})
         \in {\mathbb C}^{d} \colon \sum_{j=1}^{d} |
\lambda_{j}|^{2} < 1
         \right\}.
          $$
          The point-evaluation map
          \begin{equation}  \label{pteval-free1}
         f^{\wedge}(\lambda, b) = \sum_{n \in {\mathbb
Z}^{d}_{+}} \left[
         \sum_{\alpha  \in \cF_{d}\colon |\alpha| = n }  b
f_{\alpha}  \right]
         \lambda^{n} = b \cdot f^{\wedge}(\lambda, I_{\cE})
          \end{equation}
          where we use the standard $(\cL(\cE), {\mathbb
          C})$-correspondence structure on $\cE$. Here also
        we use the standard commutative multivariable notation
          $$  \lambda^{n} = \lambda_{1}^{n_{1}} \cdots \lambda_{d}^{n_{d}}
          \text{ if } n = (n_{1}, \dots, n_{d}) \in {\mathbb Z}_{+}^{d}.
          $$
           From \eqref{pteval-free1} we see that we are in the situation of
          Remark \ref{R:pos-ker} and completely positivity of the kernel
          $$
          \K(\lambda, \lambda')[b^{*}b'] =\sum_{n \in {\mathbb
         Z}^{d}_{+}}(\lambda^{n}I_{\cE})
         (b^{*}b')(\overline{\lambda}^{\prime})^{n} I_{\cE})
          $$
          associated with $H^{2}(E, \sigma)$ for this case
reduces to classical
          Aronszajn positivity for the Drury-Arveson kernel
          $$
         k(\lambda, \lambda') = \sum_{n \in {\mathbb
Z}^{d}_{+}} \lambda^{n}
         (\overline{\lambda}^{\prime})^{n} = \frac{1}{1 -
         \langle \lambda,
         \lambda' \rangle}.
          $$
          In this case the Fourier transform map $\Phi \colon f \mapsto
          f^{\wedge}$ has a kernel with the cokernel given by the symmetric
          Fock space spanned by symmetric tensors
          $$
        \left\{ \sum_{\alpha \in \cF_{d} \colon |\alpha| = n}
\{ \delta_{\alpha,
        \alpha'}\}_{\alpha' \in \cF_{d}} e \colon e \in \cE \right\}
          $$
          where $\delta_{\alpha, \alpha'}$ is the standard Kronecker delta
          $$
         \delta_{\alpha, \alpha'} = \begin{cases} 1 &\text{if
} \alpha =
         \alpha', \\
         0 &\text{otherwise.}
          \end{cases}
          $$
          Then it is known (see \cite{davidsonpitts,
ariaspopescu, BBF1}) that
          the image of $\Phi$ in this case, i.e., the space
$H^{2}(E, \sigma)$
          of all functions on the ball of
          the form $f^{\wedge}$ for an $f \in \ell^{2}_{\cE}(\cF_{d})$, is
          exactly the Arveson-Drury space and the associated space
          $H^{\infty}(E, \sigma)$ is exactly the multiplier space
$\cM(\cE)$ of
          the Arveson space.
          When we specialize the general Theorem
          \ref{T:MS-Schur} to this case we simply recover Theorem
\ref{T:BTV}
          for the case $\cU = \cY = \cE$.

          \medskip

  \textbf{Case 2:} {\em $\cG = {\mathbb C}$ and $\cE =
   \cE_{0}$ is a separable, infinite-dimensional Hilbert space:}  In this case
   the generalized unit disk ${\mathbb D}((E^{\sigma})^{*}) =
  {\mathbb B}(\cL(\cE^{d}, \cE))$ consists of row contractions
  $$
  \eta = \begin{bmatrix} \eta_{1} & \cdots & \eta_{d} \end{bmatrix}
  \colon \cE^{d} \to \cE.
  $$
  The Fock correspondence $\cF^{2}(E) $ = $\ell^{2}(\cF_{d})$
  has scalar coefficients
   while the Hilbert Fock space  $\cF^{2}(E, \sigma) =
   \ell^{2}_{\cE}(\cF_{d})$ has $\cE$-valued coefficients.
   The point-evaluation map \eqref{freepteval} has the form
  $$
  f = \{f_{\alpha}\}_{\alpha \in \cF_{d}} \mapsto f^{\wedge}(\eta, b)=
  \sum_{\alpha \in \cF_{d}} \eta^{\alpha} b f_{\alpha}.
  $$
  The completely positive kernel associated with $H^{2}(E, \sigma)$ for
  this case is
  $$
   \K_{E, \sigma}(\eta, \zeta)[b] = \sum_{\alpha \in \cF_{d}}
   \eta^{\alpha}[b] \zeta^{\alpha *}.
   $$
 where $b \in \sigma(\cA)' = \cL(\cE)$.

 The analytic Toeplitz algebra $\cF^{\infty}(E)$ is the free
 semigroup algebra $\cL_{d}$ acting on $\cF^{2}(E) =
 \ell^{2}(\cF_{d})$ having noncommutative Toeplitz matrix
 representation
 $$
  R = [R_{\alpha, \beta}]_{\alpha, \beta \in \cF_{d}} \text{ where }
   R_{\alpha, \beta} = R_{\alpha \beta^{-1}, \emptyset}.
 $$
 where the matrix entries $R_{\alpha, \beta}$ are scalars.
 Here $\emptyset$ refers to the empty word in $\cF_{d}$ (the unit
 element for the semigroup $\cF_{d}$) and we use the convention
 \begin{align*}
 & \alpha \beta^{-1} = \begin{cases} \alpha' &\text{if } \alpha =
 \alpha' \beta,  \\ \text{undefined} & \text{otherwise}, \end{cases} \\
 & R_{\text{undefined}} = 0.
 \end{align*}
 Then it is easily seen that $R \otimes_{\sigma} I_{\cE} \in
 \cL(\ell^{2}_{\cE}(\cF_{d}))$ is simply the infinite-multiplicity
 inflation of $R$:
 $$ R \otimes I_{\cE} = [R \otimes I_{\cE}]_{\alpha, \beta} \text{
 where } [R \otimes I_{\cE}]_{\alpha, \beta} = R_{\alpha, \beta}
 \otimes I_{\cE}.
 $$
 The point-evaluation $\eta \mapsto
 (R \otimes_{\sigma} I_{\cE})^{\wedge}(\eta)$ for
 $R \otimes_{\sigma} I_{\cE} \in
 \cF^{\infty}(E, \sigma)$ and $\eta = \begin{bmatrix} \eta_{1} &
 \cdots & \eta_{d} \end{bmatrix} \in {\mathbb B}(\cL(\cE^{d}, \cE))$
 is given by
 $$
  (R \otimes_{\sigma} I_{\cE})^{\wedge}(\eta) =
  \sum_{\alpha \in \cF_{d}} \eta^{\alpha} (R_{\alpha} \otimes I_{\cE}).
 $$
 Viewing the operator $R_{\alpha} \otimes_{\sigma} I_{\cE}$ as simply
 multiplication by the scalar $R_{\alpha}$, we can rewrite this as
 \begin{equation}  \label{pteval-free2}
     (R \otimes_{\sigma} I_{\cE})^{\wedge}(\eta) =
     \sum_{\alpha \in \cF_{d}} R_{\alpha} \eta^{\alpha}.
 \end{equation}
   As a consequence of the fact that there are no
polynomial identities
          valid for matrices of all sizes (see \cite[pp. 22-23]{Rowan}), it
          follows that the point-evaluation map
$$
 R \in \cF^{\infty}(E) = \cL_{d} \mapsto
 \left( \eta \in {\mathbb B}(\cL(\cE^{d}, \cE) \mapsto
    (R \otimes_{\sigma} I_{\cE})^{\wedge}(\eta) \in \cL(\cE) \right)
$$
 is injective.

 For $\pi$ a $*$-representation of $\cL(\cE)$ into $\cL(\cH)$ for a
 Hilbert space $\cH$, one can check that
 \begin{equation}  \label{free2-identify}
  \begin{bmatrix} \mu_{1} \\ \vdots \\ \mu_{d} \end{bmatrix} \otimes
      h \cong \begin{bmatrix} \pi(\mu_{1})h \\ \vdots \\ \pi(\mu_{d})h
      \end{bmatrix}
 \end{equation}
 gives an identification of $E^{\sigma} \otimes \cH$ with $\cH^{d}$.
 For a colligation ${\mathbf U}$ to be of the form
 \eqref{coiscol} and to satisfy \eqref{U-modulemap2} means that there
 is a Hilbert space $\cH$ together with a $*$-representation $\pi
 \colon \cL(\cE) \to \cL(\cH)$ such that, after the identification of
 $E^{\sigma} \otimes_{\pi} \cH$ with $\cH^{d}$ via
 \eqref{free2-identify},
 \begin{equation}  \label{free2-col}
 {\mathbf U} = \begin{bmatrix} A_{1} & B_{1} \\ \vdots & \vdots \\
 A_{d} & B_{d} \\ C & D \end{bmatrix} \colon \begin{bmatrix} \cH \\
 \cE \end{bmatrix} \to \begin{bmatrix} \cH \\ \vdots \\ \cH \\ \cE
 \end{bmatrix}
 \end{equation}
 subject to
 $$
   \begin{bmatrix} \pi(b) & &  & \\ & \ddots &  & \\ & & \pi(b) & \\
   & & & b\end{bmatrix}
       \begin{bmatrix} A_{1} & B_{1} \\ \vdots & \vdots \\ A_{d} &
       B_{d} \\ C & D \end{bmatrix} = \begin{bmatrix} A_{1} &
       B_{1} \\ \vdots & \vdots \\ A_{d} & B_{d} \\ C & D
       \end{bmatrix} \begin{bmatrix} \pi(b) & 0 \\ 0 & b
   \end{bmatrix},
 $$
 or, equivalently,
 \begin{align}
   &  A_{j} \pi(b) = \pi(b) A_{j}, \qquad B_{j} b = \pi(b) B_{j} \text{
     for } j = 1, \dots, d, \notag \\
   & C \pi(b) = b C, \qquad D b = b D,  \label{free2-relations}
  \end{align}
 for all $b \in \cL(\cE)$.
 For $\eta = \begin{bmatrix} \eta_{1} & \cdots & \eta_{d}
\end{bmatrix} \in {\mathbb D}((E^{\sigma})^{*}) = {\mathbb
B}(\cL(\cE^{d}, \cE))$, one can check that the operator
$ L_{\eta^{*}} \colon \cH \to E^{\sigma} \otimes_{\pi} \cH$ given by
$L_{\eta^{*}} \colon h \mapsto \eta^{*} \otimes h$, after
the identification \eqref{free2-identify}, is simply the column
contraction
$$  L_{\eta^{*}} \colon h \mapsto \begin{bmatrix} \eta_{1}^{*} \\
\vdots \\ \eta_{d}^{*} \end{bmatrix} h
$$
with adjoint equal to
$$
  L_{\eta^{*}}^{*} = \begin{bmatrix} \eta_{1} & \cdots & \eta_{d}
\end{bmatrix} \colon \cH^{d} \to \cH.
$$

Suppose that $S \in H^{\infty}(E, \sigma)$ for this example of $(E,
\sigma)$.
Then the realization formula  for $S \in \cF^{\infty}(E, \sigma)$
given by \eqref{cor-realization2} for this case becomes
\begin{equation}  \label{free2-MSreal}
  S(\eta) = D + C (I - \eta A)^{-1} \eta B \text{ for }
  \eta = \begin{bmatrix} \eta_{1} & \cdots & \eta_{d} \end{bmatrix}
  \in {\mathbb B}(\cL(\cE^{d}, \cE))
\end{equation}
where the coisometric ${\mathbf U}$ is as in \eqref{free2-col}.
Using the relations \eqref{free2-relations} and using the expansion
$$
  (I - \eta A)^{-1} = \sum_{n=0}^{\infty} (\eta A)^{n},
$$
we see that \eqref{free2-MSreal} can be rewritten as
$$
 S(\eta) = D +  \sum_{\alpha \in \cF_{d}} \sum_{j=1}^{d} C A^{\alpha}
 B \eta^{\alpha} \eta_{j}.
$$
Moreover, again from the relations \eqref{free2-relations} we see that
$$  (C A^{\alpha} B_{j}) b = b (C A^{\alpha} B_{j}) \text{ and } D b =
b D  \text{ for all }
b \in \cL(\cE),
$$
i.e., {\em $C A^{\alpha} B_{j} =: s_{\alpha j}$ for all
$\alpha \in \cF_{d}$  and $ j =
1, \dots, d$ as well as $D =: s^{0}_{\emptyset}$  are all scalar
operators:}
$$
  s_{\alpha} = s_{\alpha}^{0} I_{\cE} \text{ where } s_{\alpha}^{0}
  \in {\mathbb C}.
$$
From the complete positive kernel condition in Theorem
\ref{T:MS-Schur}, it is easily seen that
$S(\eta) = \sum_{\alpha \in \cF_{d}} s_{\alpha} \eta^{\alpha}$ is
contractive for each row contraction $\eta = \begin{bmatrix}
\eta_{1} & \cdots & \eta_{d} \end{bmatrix} \in {\mathbb
B}(\cL(\cE^{d}, \cE))$.  Thus the formal power series
$$
  S^{0}(z) = \sum_{\alpha \in \cF_{d}} s_{\alpha}^{0} z^{\alpha}
$$
is in the formal noncommutative Schur class with scalar coefficients
${\mathcal S}_{nc, d}({\mathbb C}, {\mathbb C})$ introduced in
Section \ref{S:free}.

Conversely, if $S^{0}(z) = \sum_{\alpha \in \cF_{d}} s^{0}_{\alpha}
z^{\alpha}$ is in the formal noncommutative Schur class ${\mathcal
S}_{nc, d}({\mathbb C}, {\mathbb C})$,  then part (3) of Theorem \ref{T:NC1}
assures us that $S(z)$ has a realization of the form
\begin{equation}  \label{formal-real}
  S^{0}(z) = D^{0} + C^{0} (I - Z(z) A^{0})^{-1} Z(z) B^{0}
\end{equation}
for a coisometric (even unitary) colligation
$$
  {\mathbf U}^{0} = \begin{bmatrix} A^{0}_{1} & B^{0}_{1} \\ \vdots &
  \vdots \\ A^{0}_{d} & B^{0}_{d} \\ C^{0} & D^{0} \end{bmatrix}
  \colon\begin{bmatrix}  \cH^{0} \\ {\mathbb C} \end{bmatrix} \to \begin{bmatrix}
  \cH^{0} \\ \vdots \\ \cH^{0} \\ {\mathbb C} \end{bmatrix}.
$$
Let us form a new tensored colligation ${\mathbf U}$ of the form
$$
  {\mathbf U} = \begin{bmatrix} A_{1} & B_{1} \\ \vdots & \vdots \\
  A_{d} & B_{d} \\ C & D \end{bmatrix} \colon
  \begin{bmatrix} \cH \\ \cE \end{bmatrix} \to
      \begin{bmatrix} \cH \\ \vdots \\ \cH \\ \cE \end{bmatrix}
$$
where we set
$$ \cH = \cH^{0} \otimes \cE, \quad A_{j} = A_{j}^{0} \otimes
I_{\cE}, \qquad B_{j} = B^{0}_{j} \otimes I_{\cE}, \qquad C = C^{0}
\otimes I_{\cE}, \qquad D = D^{0} \otimes I_{\cE}
$$
where $j = 1, \dots, d$.  We may define a $*$-representation $\pi
\colon \cL(\cE) \to \cL(\cH)$ by
$$
  \pi(b) = I_{\cH^{0}} \otimes b.
$$
Then it is easily seen that this ${\mathbf U}$ satisfies
\eqref{free2-col} and \eqref{free2-relations}. Moreover, from these
relations and the realization \eqref{formal-real} for the formal noncommutative
Schur-class function $S^{0}(z)$, we see that we have a realization
for the associated function $\eta \mapsto S^{0}(\eta)$ of the form
\eqref{free2-MSreal}:
$$
 S(\eta): = \sum_{\alpha \in \cF_{d}} s^{0}_{\alpha} \eta^{\alpha} =
         D + C (I - L_{\eta^{*}}^{*} A)^{-1} L_{\eta^{*}}^{*} B.
$$
We conclude: {\em there is a one-to-one correspondence between formal
power series $S^{0}(z) = \sum_{\alpha \in \cF_{d}} s^{0}_{\alpha}
z^{\alpha}$ in the noncommutative scalar-coefficient Schur class
${\mathcal S}_{nc, d}({\mathbb C}, {\mathbb C})$ and functions $\eta
\mapsto S(\eta)$ in the Muhly-Solel
class $H^{\infty}(E, \sigma$ for the particular choice of $(E, \sigma)$}
(described in \eqref{free} with $\cG = {\mathbb C}$ and $\cE_{0} =
\cE$ infinite-dimensional), {\em given by}
$$ S^{0}(z) = \sum_{\alpha \in \cF_{d}} s^{0}_{\alpha} z^{\alpha}
\mapsto
\left( \eta \mapsto S(\eta) : = \sum_{\alpha \in \cF_{d}} \eta^{\alpha}
(s^{0}_{\alpha} I_{\cE}) \right).
$$

Here we have made explicit the correspondence between condition (3) in
Theorem \ref{T:NC1} for $S^{0}(z)$ versus condition (3) in Theorem
\ref{T:MS-Schur} for $S(\eta)$.  An amusing exercise would be to
understand directly the equivalence between any of the other
conditions in Theorem \ref{T:NC1} for $S^{0}(z)$
and the corresponding condition for
$S(\eta)$ in Theorem \ref{T:MS-Schur}.

        \subsection{Analytic crossed-product algebras}
        \label{S:tv}

        We discuss here a particular case of analytic crossed-product
        algebras (see Example 2.6 in \cite{MS-JFA} as well as
the references
        there).   This particular case has strong connections with
        time-varying system theory and was discussed in connection with
        point-evaluation and generalized Nevanlinna-Pick
interpolation in
        \cite{MS-OT} (see Examples 2.5, 2.6 and 2.25 there).
Here we wish
        to draw out the connections between the realization
theorem (the
        equivalence of (1) and (3) in Theorem \ref{T:MS-Schur} for this
        case) and a result from \cite{ABP} that any lower-triangular
        contractive operator on $\ell^{2}({\mathbb Z})$ can
be realized as
        the input-output map of a linear time-varying
input/state/output
        system.
        For simplicity we discuss in detail only the multiplicity-free
        case ($\cU = \cY = {\mathbb C}$).

        We take the algebra $\cA$ to be the algebra
$\ell^{\infty}({\mathbb
        Z})$ with coordinate-wise multiplication with
correspondence $E$
        equal to $\cA = \ell^{\infty}({\mathbb Z})$ as a set.
Let $\alpha$
        be the automorphisms $\alpha(a)(k) = a(k-1)$ ($k \in
{\mathbb Z}$)
        for $a \colon {\mathbb
        Z} \to {\mathbb C}$ in $\cA$.  We consider $E$ as a
        correspondence over $\cA$ with left and right action given by
        $$
         (a \cdot \xi) = (\alpha(a))(k) = a(k-1) \xi(k),
         \qquad (\xi \cdot a)(k) = (\xi a)(k) = \xi(k) a(k)
         $$
        and with the $\cA$-valued inner product
        \begin{equation}  \label{tv-innerprod}
         \langle \xi', \xi \rangle_{E}(k) = \overline{\xi(k)} \xi'(k)
         \end{equation}
         for $d \in {\mathbb Z}$, $a \in \cA =
\ell^{\infty}({\mathbb Z})$
         and $\xi', \xi \in  E = \ell^{\infty}({\mathbb Z})$.
        Then it is easily seen that $E^{\otimes n}$ is the
correspondence
        over $\cA$ identified again with $E =
\ell^{\infty}({\mathbb Z})$ as
        a set with $\cA$-valued inner product as in
\eqref{tv-innerprod} but
        with left and right $\cA$-action given by
        $$
         (a \cdot \xi^{(n)})(k)= (\alpha^{n}(a) \xi^{(n)})(k)
= a(k-n) \xi(k), \qquad
         (\xi^{(n)} \cdot a)(k) =  (\xi^{(n)} a)(k) = \xi^{(n)}(k) a(k)
        $$
        for  $k \in {\mathbb Z}$, $\xi^{(n)} \in E^{\otimes
         n} = \ell^{\infty}({\mathbb Z})$ and $a \in \cA =
         \ell^{\infty}({\mathbb Z})$.
        The Fock space $\cF^{2}(E)$ is then the correspondence
        $\oplus_{n=0}^{\infty} \ell^{\infty}({\mathbb Z})$
with left and right
        $\cA$-action given by
        $$
         a \cdot( \oplus_{n=0}^{\infty} \xi^{(n)}) =
         \oplus_{n=0}^{\infty}\alpha^{n}(a) \xi^{(n)}, \qquad
         (\oplus_{n=0}^{\infty} \xi^{(n)}) \cdot a =
\oplus_{n=0}^{\infty}
         \xi^{(n)}a.
        $$
        More generally, when $\cA$ is a general von Neumann algebra and
        $\alpha$ is an automorphism of $\cA$, this
construction gives rise
        to analytic crossed-product algebras which have been
studied by a
        number of authors over the past several decades (see
\cite{MS-JFA,
        MS-Annalen, MS-OT} and the references therein).

        An appealing alternative representation of the
correspondence, as
          explained in Example 2.6
        in \cite{MS-OT}, is as follows.  View $\cA$ as the
algebra ${\mathcal
        D}$ of all
        diagonal operators acting on $\ell^{2}({\mathbb Z})$,
let $U$ be the
        bilateral shift operator $U e_{k} = e_{k+1}$ (where
$e_{k}(k')$ is
        the Kronecker delta function) and let $E = U \cD
        \subset \cL(\ell^{2}({\mathbb Z}))$.  Then define the
left and right
        actions of $\cA = {\mathcal D}$ on $E = U \cD$ simply
by left and
        right operator multiplications with the inner product given by
        $$
          \langle UD_{1}, UD_{2}\rangle_{\cE} = D_{2}^{*}D_{1} \in \cD.
        $$
        One can check that this $\cA$-correspondence $E$
        is unitarily equivalent to
        the $\cA'$-correspondence $E'$
        given in the previous paragraph with the obvious
        identifications:
        \begin{align*}
          & \{d(k)\}_{k \in {\mathbb Z}} \in \cA' =
\ell^{\infty}({\mathbb Z})
        \mapsto \text{diag}_{k \in {\mathbb Z}} d(k) \in \cD = \cA, \\
        &
        \{\xi(k)\}_{k \in {\mathbb Z}} \in E' =
\ell^{\infty}({\mathbb Z})
        \mapsto U \cdot (\text{diag}_{k \in {\mathbb Z}}
\xi(k)) \in U \cD = E.
        \end{align*}
        One can easily check that the identification map
        $$UD_{2} \otimes UD_{1} \in U \cD \otimes U \cD \mapsto
        UD_{2}UD_{1}=U^{2}(U^{*}D_{2}U)D_{1} \in U^{2} \cD.
        $$
        is unitary from $E \otimes E$ to $U^{2} \cD$.  After this
        identification,  the left and right $\cD$-action on
$E^{\otimes 2}
        \cong U^{2} \cD$ is
         again given by left and right operator multiplication.  More
         generally, we view $E^{\otimes n}$ as $U^{n} \cD$
with left and
         right $\cD$-action given by operator multiplication
and with inner
         product inherited from $\cL(\ell^{2}({\mathbb Z}))$:
         $$
         \langle U^{n}D_{1}, U^{n} D_{2} \rangle_{E^{\otimes n}} =
         D_{2}^{*}U^{n*} U^{n} D_{1} = D_{2}^{*} D_{1} \in \cD.
         $$
         The Fock space $\cF^{2}(E)$ can then be identified with lower
         triangular matrices $T$ with diagonal expansion $T =
         \sum_{n=0}^{\infty} U^{n}D_{n}$ ($D_{n} \in \cD$) such that
         \begin{equation} \label{tvF2}
           \sum_{n=0}^{N}  D_{n}^{*}D_{n} \text{ is bounded above in }
           \cD.
           \end{equation}
           The Toeplitz algebra
         $\cF^{\infty}(E)$ consists of all lower triangular
matrices $R$
         which give rise to bounded operators on
$\ell^{2}({\mathbb Z})$.
         As elements of $\cF^{\infty}(E)$, they act on
$\cF^{2}(E)$ (lower
         triangular matrices satisfying \eqref{tvF2}) via
multiplication on
         the left. We can view this algebra as generated by a single
         creation operator $T_{I}$ (the creation operator
associated with
         the identity matrix $I \in \cD$, namely the bilateral shift
         operator $U$), together with the diagonal operators $\cD$.
         Note that $U$ is really a unilateral shift operator
since it is
         restricted to the space $\cF^{2}(E)$ of lower
triangular matrices
         (with action equal to a shifting of the subdiagonals).

         We now set $\cE = \ell^{2}({\mathbb Z})$ and let
$\sigma$ be the
         identity representation of $\cD$ on $\cE =
\ell^{2}({\mathbb Z})$.
         Then $E^{\otimes n} \otimes_{\sigma} \cE$ can be
identified with
         $\cE = \ell^{2}({\mathbb Z})$ in the natural way
           $$\iota \colon  U^{n}D \otimes e \mapsto  D e.
           $$
           When this is done the left action of $\cA = \cD$ becomes
           $$  d \cdot e = U^{n*}d U^{n} e
           $$
           since
           \begin{align*}
           \iota (d \cdot (U^{n}D \otimes e) )  & =
           \iota (d U^{n} D \otimes e )=
           \iota (U^{n} (U^{n*}d U^{n}) D \otimes e ) \\
         &   =  U^{n*}d U^{n} D e = U^{n*} d U^{n} \iota(
U^{n} D \otimes e).
           \end{align*}
           Hence we identify $\cF^{2}(E, \sigma) = \cF^{2}(E)
           \otimes_{\sigma} \cE$ with
           $$
            \cF^{2}(E, \sigma) = \ell^{2}_{\ell^{2}({\mathbb
Z})}({\mathbb
            Z}_{+})
           $$
           with left action by $\cA = \cD$ given by
           $$
            b \cdot \{ e_{n}\}_{n \in {\mathbb Z}_{+}} =
            \{ U^{*n} b U^{n} e_{n} \}_{n \in {\mathbb Z}_{+}}.
           $$
           One can see that the image of the generating
creation operator
           ${\mathbf T}_{I} = U
           \otimes I_{\cE}$ after these identifications is
           the unilateral shift operator $S \otimes
I_{\ell^{2}({\mathbb
           Z})}$ acting on $\cF^{2}(E, \sigma) =
\ell^{2}_{\ell^{2}({\mathbb
           Z})}({\mathbb Z}_{+})$:
           $$
            {\mathbf T}_{I} = [t_{i,j}]_{i,j \in {\mathbb
Z}_{+}} \text{
            where }
            t_{i,j} = \begin{cases} I_{\ell^{2}({\mathbb Z})}
&\text{if } i
            = j+1, \\
            0 &\text{otherwise.}
            \end{cases}
            $$
            The elements $R$ of $\cF^{\infty}(E, \sigma) \subset
            \cL(\ell^{2}_{\ell^{2}({\mathbb Z})}({\mathbb
Z}_{+}))$ can then
            be identified as the following algebra of sparse
matrices: there
            is a
            sequence $\{d_{n}\}_{n \in {\mathbb Z}} \subset
\cD$ of diagonal
            operators on $\ell^{2}({\mathbb Z})$ so that $R$
has the form
            \begin{equation}  \label{tv-Toeplitz}
             R = [R_{i,j}]_{i,j \in {\mathbb Z}_{+}} \text{ where }
             R_{i,j} = \begin{cases} U^{*j}d_{i-j}U^{j}
&\text{for } i \ge
             j,\\
              0 &\text{otherwise,}
             \end{cases}
           \end{equation}
           or, in block-matrix form,
           $$
           R = \begin{bmatrix} d_{0} & 0 & 0 & \dots \\
           d_{1} & U^{*} d_{0}U & 0 & \dots \\
           d_{2} & U^{*}d_{1}U & U^{*2} d_{0} U^{2} & \dots \\
           \ddots & \ddots & \ddots & \ddots  \end{bmatrix}.
           $$

           We identify $(E^{\sigma})^{*}$ for this example as follows.
           The space $(E^{\sigma})^{*}$ consists of operators
$\eta \colon E
           \otimes_{\sigma} \cE \to \cE$ such that $\eta
(\varphi(a) \otimes
           I_{\cE}) = \sigma(a) \eta$.  For the present
situation, both $E
           \otimes_{\sigma} \cE$ and $\cE$ are identified with
           $\ell^{2}({\mathbb Z})$ but the left action by an
element $d \in
           \cA = \cD$ is given by multiplication by $U^{*}dU$
in the first
           case and by multiplication by $d$ in the second.  Thus the
           operator $\eta \in \cL( \ell^{2}({\mathbb Z}))$ is
required to
           satisfy
           $$
             \eta U^{*} D U = D \eta
           $$
           which means that $\eta U^{*}$ is diagonal, so
$(E^{\sigma})^{*}$
           is identified with weighted shift operators
           \begin{equation} \label{tvEsigma*}
           (E^{\sigma})^{*} \cong \{\eta = D_{\eta} U = U
(U^{*}D_{\eta}U)
           \in \cL(\ell^{2}({\mathbb Z}))\colon
           D_{\eta} \in \cD \} = U \cD.
          \end{equation}
          Recall that there is a representation of
         $\cF^{\infty}(E^{\sigma})$ on $\cF^{2}(E, \sigma)$ (where
         $E^{\sigma}$ is viewed as a $\sigma(\cA)'$-correspondence).
         For our situation here,
         $\sigma(\cA)' = \cA = \cD$ considered as acting on $\cE =
         \ell^{2}({\mathbb Z})$ and the representation of
$\sigma(\cA)'$ on
         $\cF(E, \sigma) = \ell^{2}_{\ell^{2}({\mathbb Z})}({\mathbb
         Z}_{+})$ turns out to be the diagonal action:
         \begin{equation}  \label{tv-b-action}
         b \cdot (\oplus_{n=0}^{\infty} e_{n}) =
\oplus_{n=0}^{\infty} b
         e_{n} \text{ for } b \in \cD,\,
\oplus_{n=0}^{\infty} e_{n} \in
         \ell^{2}_{\ell^{2}({\mathbb Z})}({\mathbb Z}_{+}).
         \end{equation}
         For purposes of getting a generating set for
         $\cF^{\infty}(E^{\sigma})$,
         it suffices to consider the single creation operator
associated
          with
         $\eta^{*} = U^{*}$: the associated action on
$\cF^{2}(E, \sigma)$
         turns out to be
         \begin{equation}  \label{tv-creation-rep}
        T_{U^*,\sigma}^d = [t'_{i,j}] \text{ where }
         t'_{i,j} = \begin{cases} U^{*} &\text{if } j = i+1, \\
             0 &\text{otherwise.}
             \end{cases}
         \end{equation}
         According to the duality result from
\cite{MS-Annalen}, an operator
         $R$ on $\ell^{2}_{\ell^{2}({\mathbb Z})}({\mathbb Z}_{+})$
         is of the form \eqref{tv-Toeplitz} if and only if
$R$ commutes with
         the scalar-diagonal operators \eqref{tv-b-action} and the
         $E^{\sigma}$-creation operator
\eqref{tv-creation-rep}; an amusing
         exercise for the reader is to verify this fact
directly for this
         example.

         We now identify the $Z$-transform and compute the
function spaces
         $H^{2}(E, \sigma)$ and $H^{\infty}(E, \sigma)$ as follows. By
          \eqref{tvEsigma*}
         we have an identification of $(E^{\sigma})^{*}$
         with the space of weighted
         shift operators $U \cD$ in $\cL(\ell^{2}({\mathbb Z}))$.
          After carrying out the
         identifications $E^{\otimes n} \otimes_{\sigma} \cE =
         \ell^{2}({\mathbb Z})$, one can check that the
generalized power
         $\eta^{n} \colon E^{\otimes n} \otimes_{\sigma} \cE
\to \cE$ of an
         $\eta \in (E^{\sigma})^{*} = U \cD$
         coincides with the usual power $\eta^{n}$ as an element of the
         operator algebra
         $\cL(\ell^{2}({\mathbb Z}))$.  Therefore, for $f =
\{f_{n}\}_{n \in
         {\mathbb Z}_{+}} \in \ell^{2}_{\ell^{2}({\mathbb Z})}({\mathbb
         Z}_{+})$ and $\eta = D_{\eta} U \in {\mathbb
D}((E^{\sigma})^{*})$
         (with $D_{\eta} \in \cD$),
         we have
         $$
            f^{\wedge}(\eta,b) = \sum_{n=0}^{\infty} \eta^{n} b f_{n} =
           \sum_{n=0}^{\infty} (D_{\eta} U)^{n} b f_{n}.
         $$
         If we restrict the second variable $b \in
\sigma(\cA)' = \cD$ to be
         $b = I_{\ell^{2}({\mathbb Z})}$, we have the
restricted Fourier
         transform
         $$
           \Phi^{1} \colon f = \{f_{n}\}_{n \in {\mathbb
Z}_{+}} \mapsto
           f^{\wedge}(\eta, I_{\ell^{2}({\mathbb Z})}) =
\sum_{n=0}^{\infty}
           \eta^{n}f_{n} = \sum_{n=0}^{\infty} (D_{\eta}U)^{n} f_{n}.
         $$
           We assert that {\em the restricted $Z$-transform
$\Phi^{1} \colon
           f \mapsto f^{\wedge 1}: = f^{\wedge}(\cdot,
I_{\ell^{2}({\mathbb
           Z})})$
           is injective.}
         Indeed,
         if $f^{\wedge 1}(\eta) = 0$ for all $\eta$,
evaluating at $\eta =
         0$ gives that $f_{0} = 0$ and hence $ F^{\wedge
1}(\eta) = \eta \cdot
         \sum_{n=0}^{\infty} f_{n+1} \eta^{n} = 0$.  Choosing $\eta$
         invertible and premultiplying by $\eta^{-1}$ then gives that
         \begin{equation} \label{=zero}
         \sum_{n=0}^{\infty} f_{n+1} \eta^{n} = 0
         \end{equation}
         for all invertible $\eta$.  By approximating a
noninvertible $\eta$
         by invertible
         $\eta$'s, we see that \eqref{=zero} actually holds
for all $\eta \in
         {\mathbb D}((E^{\sigma})^{*})$. Iteration of the
same argument now
         gives that $f_{n} = 0$ for all $n \in {\mathbb
Z}_{+}$, i.e., $f =
         \{f_{n}\}_{n \in {\mathbb Z}_{+}}$ is the zero element of
         $\ell^{2}_{\ell^{2}({\mathbb Z})}({\mathbb Z}_{+})$, and the
         assertion follows.  Note that the $\sigma(\cA)' = \cD$-action
         on $\cF^{2}(E, \sigma)$ is given by
         $$ d \cdot \{f_{n}\}_{n \in {\mathbb Z}_{+}} = \{ d
f_{n}\}_{n \in
         {\mathbb Z}_{+}} \text{ for } d \in \cD.
         $$

         The completely positive kernel $\K$ associated with the
          reproducing kernel Hilbert
         correspondence $H^{2}(E, \sigma)$ $=$
         $\Phi(\ell^{2}_{\ell^{2}({\mathbb Z})}({\mathbb Z}_{+}))$ is
         $$
         \K(\eta, \zeta)[b] = \sum_{n=0}^{\infty} \eta^{n} b
         \zeta^{*n} \text{ for } \eta, \zeta \in {\mathbb
         D}((E^{\sigma})^{*}) = {\mathbb B}(U \cD) \text{ and
} b \in \cD.
         $$
         Note that $\Phi^{*} \K(\cdot, \zeta)[b] e = b
k_{\zeta} e$ where
         $$
          b k_{\eta} e = \{ b \zeta^{*n} e \}_{n \in
{\mathbb Z}_{+}}
         \in \ell^{2}_{\ell^{2}({\mathbb Z})}({\mathbb Z}_{+}).
         $$
         We conclude that the subcollection
         $$
           \{ b k_{\zeta} e \colon b \in \cD, \, \zeta \in
{\mathbb B}(U \cD),
           \, e \in \ell^{2}({\mathbb Z}) \}
         $$
         has dense span in $\ell^{2}_{\ell^{2}({\mathbb Z})}({\mathbb
         Z}_{+})$.

         An element $R$ of
         $\cF^{\infty}(E)$ is identified with a lower triangular matrix
         representing a bounded operator on
$\ell^{2}({\mathbb Z})$; it is
         convenient to represent such a matrix via a
generalized Fourier
         series along subdiagonals:
         \begin{equation}  \label{Fourier-series}
           R  \sim  \sum_{n=0}^{\infty} U^{n} d_{n} \text{
where } d_{n} \in
          \cD.
         \end{equation}
         (The Caesaro averages of the partial sums of the
series converges to
         $R$ in the weak-$*$ topology but we shall not need this.)
         Then $R \otimes I_{\ell^{2}({\mathbb Z})}$, after the
         identification of $\cF^{2}(E) \otimes_{\sigma}
\ell^{2}({\mathbb
         Z})$ with $\ell^{2}_{\ell^{2}({\mathbb Z})}({\mathbb
Z}_{+})$, is
         identified with the operator acting on
         $\ell^{2}_{\ell^{2}({\mathbb Z})}({\mathbb Z}_{+})$
with the sparse
         matrix representation \eqref{tv-Toeplitz}.  For $\eta \in
         {\mathbb D}((E^{\sigma})^{*}) = {\mathbb B}(U \cD)$,
the associated
         point evaluation of
         $R \otimes I_{\cE}$  is then given by
         \begin{align}
         (R \otimes I_{\ell^{2}({\mathbb Z})})^{\wedge}(\eta) & =
         \sum_{n=0}^{\infty} \eta^{n} R_{n,0} \notag \\
         & = \sum_{n=0}^{\infty} \eta^{n} d_{n}
         \label{tv-MS-evaluation}
         \end{align}
         if $R$ is given by \eqref{Fourier-series}. In
particular, formally
         we recover $R$ from $(R \otimes I_{\ell^{2}({\mathbb
Z})})^{\wedge}$ as
         \begin{equation}  \label{recover'}
           R = (R \otimes I_{\ell^{2}({\mathbb Z})})^{\wedge}(U).
         \end{equation}
         More precisely, we interpret the right-hand side of
\eqref{recover'}
         as
         $$
         (R \otimes I_{\ell^{2}({\mathbb Z})})^{\wedge}(U) =
         \lim_{r \uparrow 1} (R \otimes I_{\ell^{2}({\mathbb
Z})})^{\wedge}(rU).
         $$

         The realization theorem (the equivalence of (1) and
(3) in Theorem
         \ref{T:MS-Schur}) assures us that any function of the form $(R
         \otimes I_{\ell^{2}({\mathbb Z})})^{\wedge}$ with
$\|R\| \le 1$ can
         be realized as follows.  Suppose first that $\cH$ is a
         $(\sigma(\cA)' = \cD, {\mathbb C})$-correspondence,
i.e., $\cH$ is
         a Hilbert space and there is a $*$-representation $\pi$ of
         $\sigma(\cA)' = \cD$ with values in $\cL(\cH)$.
         Noting that
         $$
          U^{*}d  \otimes_{\pi} h = U^{*} \otimes_{\pi}
\pi(d) h \cong \pi(d) h
          $$
         for $\mu = U^{*}d \in E^{\sigma} = U^{*} \cD$
         (so $d \in \cD = \sigma(\cA)'$) and $h
         \in \cH$, we see that $E^{\sigma} \otimes_{\pi} \cH$
can be identified with
         $\cH$, but at the price that the left ($\sigma(\cA)'
= \cD$)-action
         on $\cH$ is given by $\pi^{(1)} \colon b \mapsto
\pi(U b U^{*})$
         rather than by $\pi$.  With this identification, we see that
         the unitary colligation ${\mathbf U}$ in
         \eqref{coiscol} and \eqref{U-modulemap2} for this
case has the form
         $$
         {\mathbf U} = \begin{bmatrix} A & B \\ C & D
\end{bmatrix} \colon
         \begin{bmatrix} \cH \\ \ell^{2}({\mathbb Z}) \end{bmatrix} \to
             \begin{bmatrix} \cH \\ \ell^{2}({\mathbb Z}) \end{bmatrix}
          $$
          subject to
         \begin{equation*}
          \begin{bmatrix} A & B \\ C & D \end{bmatrix} \begin{bmatrix}
         \pi(b) & 0 \\ 0 & b \end{bmatrix} =
         \begin{bmatrix} \pi^{(1)}(b) & 0 \\ 0 & b \end{bmatrix}
             \begin{bmatrix} A & B \\ C & D \end{bmatrix}
          \text{ for all $b \in \sigma(\cA)' = \cD$,}
          \end{equation*}
      or, equivalently,
      \begin{equation}  \label{tv-modulemap}
  A \pi(b) = \pi(U b U^{*}) A, \qquad B b = \pi(UbU^{*}) B, \qquad
  C \pi(b) = b C, \qquad D b = b D \text{ for all } b \in \cD.
\end{equation}
          The realization theorem then tells us that any $(R \otimes
          I_{\ell^{2}({\mathbb Z})})^{\wedge}$ (where $R \in
          \cL(\ell^{2}({\mathbb Z}))$  is lower-triangular
and contractive)
          can be realized as
          \begin{equation}  \label{tv-realization}
          (R \otimes I_{\ell^{2}({\mathbb Z})})^{\wedge}(\eta) =
          D + C (I - \pi(\eta U^{*}) A)^{-1} \pi( \eta U^{*}) B.
          \end{equation}

          Let us now consider a time-varying input/state/output linear
          system of the form
          \begin{equation}  \label{tv-sys}
           \Sigma:  \left\{ \begin{array}{rcl}
           x(n+1) & = & A(n) x(n) + B(n) u(n) \\
              y(n) & =  & C(n) x(n) + D(n) u(n).
              \end{array} \right.
          \end{equation}
          determined by the time-varying system matrix
          $${\mathbf U}(n) = \begin{bmatrix} A(n) & B(n) \\ C(n) & D(n)
          \end{bmatrix} \colon
          \begin{bmatrix} \cH(n) \\ {\mathbb C} \end{bmatrix} \to
          \begin{bmatrix} \cH(n+1) \\ {\mathbb C} \end{bmatrix}.
           $$
           We say that the system is {\em conservative}
(respectively, {\em
           dissipative}) if each $U(n)$ is
           unitary (respectively, contractive).
           Let us assume that we have a dissipative time-varying linear
           system with time-varying system matrix ${\mathbf U}(n) =
           \sbm{A(n) & B(n) \\ C(n) & D(n)}$.
           Then it can be shown that, given an input string
$\{u(n)\}_{n \in
             {\mathbb Z}}$ in $\ell^{2}({\mathbb Z})$, there
is a unique system
             trajectory $(u(n), x(n), y(n))$, i.e., solution
of the system
             equations \eqref{tv-sys}, such that $\lim_{n \to
-\infty} x(n) = 0$
             with the resulting output string $\{y(n) \}_{n
\in {\mathbb Z}} \in
             \ell^{2}({\mathbb Z})$.  In this way there is
defined an input-output map
             $T_{\Sigma}$ on $\ell^{2}({\mathbb Z})$ such
that $T_{\Sigma} \colon
             \{u(n)\}_{n \in {\mathbb Z}} \mapsto
\{y(n)\}_{n \in {\mathbb
             Z}}$.

           Let us introduce an aggregate
           state space
           \begin{equation}  \label{agg-cH}
            \cH = \oplus_{n \in {\mathbb Z}} \cH(n)
           \end{equation}
           and an aggregate system matrix
           \begin{equation}  \label{agg-col}
           {\mathbf U} = \begin{bmatrix} \boldsymbol \cA &
\boldsymbol \cB
           \\ \boldsymbol \cC & \boldsymbol \cD \end{bmatrix}
           \begin{bmatrix} \cH \\ \ell^{2}({\mathbb Z})
\end{bmatrix} \to
          \begin{bmatrix} \cH \\ \ell^{2}({\mathbb Z}) \end{bmatrix}
          \end{equation}
         with ${\boldsymbol \cA}$, ${\boldsymbol \cB}$, ${\boldsymbol \cC}$ and
         ${\boldsymbol \cD}$ specified by block-matrix entries
         \begin{align}
        &  [\boldsymbol \cA]_{i,j} = A(j) \delta_{i,j+1},\qquad
         [\boldsymbol \cB]_{i,j} = B(j) \delta_{i,j+1},  \notag \\
         & [\boldsymbol \cC]_{i,j} = C(j) \delta_{i,j}, \qquad
         [\boldsymbol \cD]_{i,j} = D(j) \delta_{i,j}.
         \label{boldABCD}
         \end{align}
         If the operator $\boldsymbol \cA$ has spectral radius
        strictly less than 1 as an operator on $\cH$, then one can
        compute that $T_{\Sigma}$ is given by
        \begin{equation}  \label{TSigmareal}
         T_{\Sigma} = \boldsymbol \cD + \boldsymbol \cC (I -
\boldsymbol
         \cA)^{-1}  \boldsymbol \cB \in
         \cL(\ell^{2}({\mathbb Z})).
         \end{equation}
          Even if $ \boldsymbol \cA$ does not have spectral
          radius strictly less than $1$, there are various ways whereby
          one can still make sense of the formula
          \eqref{TSigmareal}; one such is via a limit
          $$
        T_{\Sigma} = \lim_{r \uparrow 1}
        \boldsymbol \cD + \boldsymbol \cC(I - r  \boldsymbol
          \cA))^{-1} (r \boldsymbol \cB).
          $$
  From the representation \eqref{TSigmareal} for $T_{\Sigma}$ one can
  compute that $T_{\Sigma}$ has the diagonal decomposition
  $$
  T_{\Sigma} = \sum_{n=0}^{\infty} U^{n} d_{n} \text{ where }
  d_{0} = \boldsymbol \cD
  \text{ and }
  d_{n} = U^{*n} \boldsymbol \cC \boldsymbol \cA^{n-1} \boldsymbol
  \cB \text{ for } n \ge 1.
  $$
  Hence an application of \eqref{tv-MS-evaluation} gives us
  \begin{equation}  \label{eval-TSigma}
      (T_{\Sigma} \otimes I_{\ell^{2}({\mathbb Z})})^{\wedge}(\eta) =
   \boldsymbol \cD +  \sum_{n=1}^{\infty} \eta^{n} U^{*n} \boldsymbol \cC \boldsymbol
      \cA^{n-1} \boldsymbol \cB.
  \end{equation}

      Given $\cH$ in the form \eqref{agg-cH}, we may define a
      representation $\pi$ of $\cD$ by
          $$
         \pi(b) \colon \oplus_{n \in {\mathbb Z}} h(n)
\mapsto \oplus_{n \in
         {\mathbb Z}} b(n) h(n) \text{ for } b = \text{\rm
         diag}_{n \in {\mathbb Z}} \{ b(n)\} \in \cD.
          $$
          Note that if $\boldsymbol \cA, \boldsymbol \cB, \boldsymbol \cC,
          \boldsymbol \cD$ are as in \eqref{boldABCD}, then
         ${\mathbf U}$ as in
          \eqref{agg-col} satisfies the $\cD$-module property
         \eqref{tv-modulemap} (with $\boldsymbol \cA, \boldsymbol
     \cB, \boldsymbol \cC, \boldsymbol \cD$ in place of $A,B,C,D$).
       By a careful induction argument making using of these
       relations, one can show that
       $$ \boldsymbol \cC ( \pi(\eta U^{*}) \boldsymbol \cA)^{k-1}
       \pi(\eta U^{*}) = \eta^{k} U^{*k} \boldsymbol \cC \boldsymbol
       \cA^{k-1} \text{ for } k = 1,2, \dots.
       $$
   One can then show that
  \begin{align*}
      \boldsymbol \cD + \boldsymbol \cC ( I -  \pi(\eta U^{*})
  \boldsymbol \cA)^{-1} \pi(\eta U^{*}) \boldsymbol \cB  & =
  \boldsymbol \cD + \sum_{n=1}^{\infty} \boldsymbol \cC ( \pi(\eta
  U^{*}) \boldsymbol \cA)^{n-1} \pi(\eta U^{*}) \boldsymbol \cB \\
  & = \boldsymbol \cD + \sum_{n=1}^{\infty} \eta^{n} U^{*n}
  \boldsymbol \cC \boldsymbol \cA^{n-1} \boldsymbol \cB \\
  & = (T_{\Sigma} \otimes I_{\ell^{2}({\mathbb Z})})^{\wedge}(\eta),
  \end{align*}
  i.e., {\em the aggregate colligation ${\mathbf U} = \sbm{ \boldsymbol
  \cA & \boldsymbol \cB \\ \boldsymbol \cC & \boldsymbol \cD}$
  arising from the realization of $T_{\Sigma}$ as the input-output
  map for the time-varying linear system \eqref{tv-sys} gives rise to
  a realization of the form \eqref{tv-realization} for the function
  $(T_{\Sigma}\otimes I_{\ell^{2}({\mathbb Z})})^{\wedge}$
  in the Muhly-Solel Schur class for this
  special setting.}

 This suggests a different approach to the realization theorem (the
 equivalence of (1) and (3) in Theorem \ref{T:MS-Schur}) for this
 particular case.  Given a contractive lower-triangular operator $R$ on
 $\ell^{2}({\mathbb Z})$, it is known (see \cite[Theorem 6.2]{ABP}) that one can
realize $R$ as the input-output map $R = T_{\Sigma}$ of a conservative
time-varying input/state/output linear system as in \eqref{tv-sys}; the solution
 in \cite{ABP} is given via a time-varying analogue of the Pavlov
 functional model, or, alternatively, via a time-varying
 analogue of the Sz.-Nagy-Foias or de Branges-Rovnyak functional model.
Once we have realized $R$ as $R=T_{\Sigma}$ with $\Sigma$ as in
\eqref{tv-sys}, we get $(R \otimes I_{\ell^{2}({\mathbb Z})})^{\wedge}$
realized in the form \eqref{tv-realization} and hence we have recovered the
implication (1)  $\Longrightarrow$ (3) of Theorem \ref{T:MS-Schur}.  We
conclude that the Muhly-Solel realization theorem for this case, after some
 translation, has essentially the same content as the conservative realization
 theorem for linear time-varying systems in \cite{ABP}.


\begin{thebibliography}{10}


%\bibitem{agler2}
%J.~Agler, Some interpolation theorems of Nevanlinna--Pick type,
%Preprint, 1988.



\bibitem{agler-hellinger}
J.~Agler, On the representation of certain holomorphic functions
defined on a polydisk, in {\em Topics in Operator Theory: Ernst D.
Hellinger memorial Volume} (L.~de~Branges, I.~Gohberg and J.~Rovnyak,
eds.), pp. 47--66,  \textbf{OT 48},
Birkh\"auser Verlag, Basel, 1990.

\bibitem{aglmccar-poly} J.~Agler and J.E.~McCarthy,
   Nevanlinna-Pick interpolation on the bidisk, {\em J.~Reine
   Angew.~Math.} {\bf 506} (1999), 191--204.



\bibitem{aglmccar} J.~Agler and J.E.~McCarthy,  Complete
Nevanlinna-Pick kernels, {\em J.~Functional Analysis} {\bf  175} (2000),
111--124.

\bibitem{ABP} D.~Alpay, J.A.~Ball and Y.~Peretz, System theory,
operator models and scattering: the time-varying case, {\em
J.~Operator Theory} \textbf{47} (2002), 245-286.

%\bibitem{ADD} D.~Alpay, P.~Dewilde and H.~Dym, Lossless inverse
%scattering and reproducing kernels for upper triangular operators, in
%{\em Extensions and Interpolation of Linear Operators and Matrix
%Functions} (Ed.~I.~Gohberg), pp. 61--135,  \textbf{OT 47},
%Birkh\"auser-Verlag,
%Basel-Boston, 1990.

%\bibitem{AK1} D.~Alpay and D.S.~Kalyuzhny\u{\i}-Verbovetzki\u{\i},
%On the intersection of null spaces for matrix substitutions in a
%non-commutative rational formal power series, {\em Comptes rendus
%Mathematiques Acad. Sci. Paris I} {\bf 339} (2004) 533--538.

%\bibitem{AK2} D.~Alpay and D.S.~Kalyuzhny\u{\i}-Verbovetzki\u{\i},
%\emph{Matrix-$J$-unitary non-commutative rational formal power
%series}, in {\em The State Space Method Generalizations and Applications}
%(Ed. D. Alpay and I. Gohberg), pp. 49--113,
%\textbf{OT 161}, Birkh\"auser-Verlag, Basel-Boston-Berlin, 2006


%\bibitem{AmRmk}
%C.-G. Ambrozie, Remarks on the operator-valued interpolation for
%multivariable bounded analytic functions. {\em Indiana Univ. Math. J.},
%53(6):1551--1576, 2004.

%\bibitem{AE}  C.-G.~Ambrozie and J.~Eschmeier, A commutant
%lifting theorem on analytic polyhedra, {\em Proceedings of Operator
%Theory Conference Dedicated to Prof.  Wieslaw Zelazko}, Banach Center
%publ., Warszawa, to appear.


\bibitem{at}
C.-G.~Ambrozie and D.~Timotin,
A von Neumann type inequality for certain domains in
${\mathbb C}^{n}$, {\em Proc.~Amer.~Math.~Soc.}, \textbf{131} (2003), 859--869.


\bibitem{ariaspopescu}
A.~Arias and G.~Popescu, {\em Noncommutative interpolation and Poisson
transforms}, Israel J. Math. {\bf 115} (2000), 205--234.

\bibitem{aron}
N.~Aronszajn, Theory of reproducing kernels,
{\em Trans.~Amer.~Math.~Soc.}, \textbf{68} (1950), 337--404.

%\bibitem{arovgros1}
%D.Z.~Arov and L.~Z.~Grossman,
%Scattering matrices in the theory of unitary
%extensions of isometric operators,
%{\em Soviet Math. Dokl.} {\bf 270} (1983), 17--20, MR0705184 (85c:47008),
%Zbl 0543.47010.

%\bibitem{arovgros2}
%D.Z.~Arov and L.~Z.~Grossman,
%Scattering matrices in the theory of unitary
%extensions of isometric operators,
%{\em Math. Nachr.} {\bf 157} (1992), 105--123.


%\bibitem{arv}
%W.~Arveson, Subalgebras of $C^*$-algebras. III. Multivariable
%operator theory, {\em Acta Math.} \textbf{181} (1998), no. 2,
%159--228.





\bibitem{Ball-Winnipeg} J.A. Ball, Linear systems, operator model
theory and
scattering: multivariable generalizations, in {\em  Operator Theory
and
Its Applications (Winnipeg, MB, 1998)} (Ed. A.G. Ramm,
P.N. Shivakumar and A.V. Strauss), Fields Institute Communications
Vol. 25, Amer. Math. Soc., Providence, 2000, pp. 151-178.

%\bibitem{bb3}
%J.A.~Ball and V.~Bolotnikov,  A tangential interpolation problem
%on the distinguished boundary of the polydisk for the Schur-Agler class,
%{\em J.~Mathematical Analysis and Applications} {\bf 273} (2002),
no. 2, 328--348.

\bibitem{bb4}
J.A.~Ball and V.~Bolotnikov, Realization and interpolation
for Schur-Agler-class functions on domains with matrix polynomial
defining function in ${\mathbb C}^n$, {\em J.~Functional Analysis}
\textbf{213} (2004), 45--87.


%\bibitem{bb5}
%J.A.~Ball and V.~Bolotnikov, Interpolation problems with
%operator argument for contractive-valued functions on general domains
%in ${\mathbb C}^{n}$, {\em New York J.~Math.} \textbf{11} (2005), 1-44.

%\bibitem{BB-noncomint}  J.A.~Ball and V.~Bolotnikov, Interpolation in
%the noncommutative Schur-Agler class, {\em J.~Operator Theory}, to
%appear.

\bibitem{BBF1} J.A.~Ball, V.~Bolotnikov and Q.~Fang, Multivariable
backward-shift invariant subspaces and observability operators, {\em
Multidimensional Systems and Signal Processing}, to appear.

\bibitem{BBF2a} J.A.~Ball, V.~Bolotnikov and Q.~Fang,
Transfer-function realization for multipliers of the Arveson space,
{\em J.~Mathematical Analysis and Applications}, to appear.


%\bibitem{BGK1} J.A.~Ball, I.~Gohberg and M.A.~Kaashoek,
%Nevanlinna-Pick
%interpolation for time-varying input-output maps: the discrete case,
%in {\em Time-Variant Systems and Interpolation} (Ed. I. Gohberg),
%pp.1--51, OT \textbf{56} Birkh\"auser
%Verlag, Basel, 1992.

%\bibitem{BGK2}  J.A.~Ball, I.~Gohberg and M.A.~Kaashoek, A frequency
%response
%function for linear, time-varying systems, {\em  Math. Control
%Signals Systems} \textbf{8} (1995), 334-351.

%\bibitem{BGR} J.A.~Ball, I.~Gohberg and L.~Rodman, \emph{Interpolation
%of Rational Matrix Functions}, \textbf{OT 45} Birkh\"auser-Verlag, 1990.

%\bibitem{BGM1} J.A.~Ball, G.~Groenewald and T.~Malakorn,
%Structured noncommutative multidimensional linear systems,
%{\em SIAM J.~Control and Optimization}  \textbf{44}
%No. 4 (2005), 1474-1528.

\bibitem{BGM2} J.A.~Ball, G.~Groenewald and T.~Malakorn, Conservative
structured noncommutative multidimensional linear systems,
in {\em  The State Space Method Generalizations and Applications}
(Ed. D. Alpay and I. Gohberg),
pp. 179-223, \textbf{OT 161}, Birkh\"auser-Verlag,
Basel-Boston-Berlin, 2006.

\bibitem{B-KV} J.A.~Ball and D.S.~Kaliuzhnyi-Verbovetskyi,
Conservative dilations of dissipative multidimensional systems:  the
commutative and non-commutative settings, {\em Multidimensional
Systems and Signal Processing}, to appear.


%\bibitem{BLTT}
%J.A.~Ball, W.~S.~Li, D.~Timotin and T.~T.~Trent, A
%commutant lifting theorem on the polydisc,
%{\em Indiana University Math. J.} {\bf 48} (1999), 653--675.

%\bibitem{BallMal}  J.A.~Ball and T.~Malakorn, Multidimensional
%linear feedback control systems and interpolation problems for
%multivariable holomorphic functions,
%{\em Multidimensional Systems and Signal Processing} \textbf{15} (2004),
%7--36.

%\bibitem{BSV} J.A.~Ball, C.~Sadosky and V.~Vinnikov, Scattering
%systems with several evolutions and multidimensional
%input/state/output systems,  {\em Integral Equations and Operator
%Theory}  \textbf{52} (2005), 323-393.

\bibitem{BT}
J.A.~Ball and T.~Trent,
Unitary colligations, reproducing kernel Hilbert
spaces and Nevanlinna--Pick interpolation in several variables,
{\em J. Functional Analysis} {\bf 157} (1998), no.1, 1--61.


\bibitem{BTV}
J.A.~Ball, T.~T.~Trent and V.~Vinnikov, Interpolation and
commutant lifting for multipliers on reproducing kernels Hilbert
spaces, in  {\em Operator Theory and Analysis: The M.A. Kaashoek
Anniversary Volume (Workshop in Amsterdam, Nov. 1997)}, OT 122,
Birkh\"auser-Verlag, Basel-Boston, 2001, pp. 89--138.

\bibitem{NFRKHS} J.A.~Ball and V.~Vinnikov, Formal reproducing kernel
Hilbert spaces: the commutative and noncommutative settings, in {\em
Reproducing Kernel Spaces and Applications} (Ed.~D.~Alpay), pp.
77--134, \textbf{OT 143}, Birkh\"auser-Verlag, Basel-Boston, 2003.

\bibitem{Cuntz-scat} J.A.~Ball and V.~Vinnikov, \emph{Lax-Phillips
scattering and conservative linear systems: a Cuntz-algebra
multidimensional setting},  Mem.~Amer.~Math.~Soc.  \textbf{178} (2005),
no. \textbf{837}.

\bibitem{BBLS} S.D.~Barreto, B.V.R.~Bhat, V.~Liebscher and M.~Skeide,
Type I product systems of Hilbert modules, {\em J. Functional
Analysis} \textbf{212} (2004), 121-181.

%\bibitem{BES} T.~Bhattacharyya, J.~Eschmeier and J.~Sarkar,
      %   Characteristic function of a pure commuting contractive tuple,
    %{\em Integral Equations and Operator Theory} {\bf 53} (2005),
    %no.1, 23--32.

\bibitem{dBR}
L.~de Branges and J.~Rovnyak,
   Canonical models in quantum scattering theory,
    in: {\em Perturbation Theory and its
    Applications in Quantum Mechanics} (C.~Wilcox, ed.) pp.
295--392, Holt, Rinehart and Winston, New York, 1966.


%\bibitem{CJ}  T.~Constantinescu and J.L.~Johnson, A note on
%noncommutative interpolation, {\em Canadian Math.~Bull.} \textbf{\bf
%46} (2003)
%no. 1
%59--70.


%\bibitem{CH} R.E.~Curto and D.A.~Herrero, On closures of joint
%similarity orbits,  {\em Integral Equations and Operator Theory}
%\textbf{8} (1985), 489--556.


\bibitem{Davidson} K.~Davidson, {\em $C^{*}$-Algebras by EXample},
Fields Institute Monographs Vol. 6, American Mathematical Society,
Providence, 1996.


\bibitem{DP} K. Davidson and D. Pitts, The algebraic structure
of non-commutative analytic Toeplitz algebras,
{\em Math.~Ann.} \textbf{311} (1998),
275-303.


\bibitem{davidsonpitts}
K.R.~Davidson and D.R.~Pitts, Nevanlinna--Pick interpolation for
non-commutative analytic Toeplitz algebras, {\em Integral
Equations Operator Theory} \textbf{31} (1998), no. 3, 321--337.

%\bibitem{DD} P.~Dewilde and H.~Dym, Interpolation for upper
%triangular operators,
%in {\em Time-Variant Systems and Interpolation} (Ed. I. Gohberg), pp.
%153--260.
%\textbf{OT56} Birkh\"auser-Verlag, Basel, 1992.


\bibitem{DMM} M.~Dritschel, S.~Marcantognini, S.~McCullough, {\em
Interpolation in Semigroupoid Algebras}, preprint.

\bibitem{Drury} S.W.~Drury, A generalization of von Neumann's
inequality to the complex ball, {\em Proc.~Amer.~Math.~Soc.} {\bf 68}
(1978), 300--304.

\bibitem{EP} J.~Eschmeier and M.~Putinar, Spherical contractions and
interpolation problems on the unit ball, {\em J.~Reine und Angewandte
Mathematik} \textbf{542} (2004), 219-236.

%\bibitem{FFGK} C.~Foia\c{s}, A.~Frazho, I.~Gohberg and M.A.~Kaashoek,
%{\em Metric Constrained Interpolation, Commutant Lifting and
%Systems}, \textbf{OT100}, Birkh\"auser-Verlag, Boston-Basel, 1998.

%\bibitem{FM}
%E.~Fornasini and G.~Marchesini, {\it State-space
%realization theory of two-dimensional filters}, IEEE Trans. Automat.
%Contr. {\bf AC-21}, No. 4, 1976, 484--492.

\bibitem{Schur-survey} I.~Gohberg (ed.), {\em I.~Schur Methods in
Operator Theory and Signal Processing}, \textbf{OT 18}
Birkh\"auser-Verlag, Basel-Boston, 1986.


%\bibitem{HMcCV} J.W.~Helton, S.~McCullough and V.~Vinnikov,
%\emph{Noncommutative convexity arises from linear matrix
%inequalities}, {\em J.~Functional Analysis} \textbf{240} (2006) no.
%1, 105--191.


%\bibitem{Kac}
%T.~Kaczorek, {\it Two-Dimensional Linear Systems},
%Lecture Notes in Control and Information Sciences {\bf 68}
%Springer-Verlag, Berlin, 1985.

%\bibitem{KV} D.S.~Kalyuzhny\u{\i}-Verbovetzki\u{\i} and V.~Vinnikov,
%Non-commutative positive kernels and their matrix functions,
%{\em Proceedings of the American Mathematical Society} \textbf{134}
%(2006) no. 3, 805--816.

\bibitem {KP}D. Kribs and S. Power, Free semigroupoid
algebras, {\em J. Ramanujan Math.~Soc.}, \textbf{19} (2004), 75-117

\bibitem{KP06} D.~Kribs, S.~Power, The $H^\infty$ algebras of higher
rank graphs, {\em Math. Proc. Royal Irish Acad.}
\textbf{106} (2006), 199-218.

%\bibitem{mccull2}
%S.~McCullough, The local de {B}ranges-{R}ovnyak construction and
%complete Nevanlinna-Pick kernels, in
%{\em Algebraic methods in operator theory} (Ed. R.~Curto and
%P.E.T.~Jorgensen), Birkh\"auser--Verlag, Boston, 1994, pp. 15--24.

%\bibitem{Mal}  T.~Malakorn, {\em Multidimensional Linear Systems and
%Robust Control}, Dissertation, Department of Electrical and Computer
%Engineering, Virginia Tech (April, 2003).

\bibitem{MT} V.M.~Manuilov and E.V.~Troitsky, {\em Hilbert
$C^{*}$-Modules}, Translations of Mathematical Monographs Volume
\textbf{226}, American Mathematical Society, Providence (2005).

\bibitem {pM97}P.S. Muhly, A finite-dimensional introduction to
operator algebra, in {\em Operator algebras and applications (Samos, 1996)},
pp. 313--354,
NATO Adv. Sci. Inst. Ser. C Math. Phys. Sci., 495, Kluwer Acad. Publ.,
Dordrecht, 1997.

\bibitem{MS-JFA} P.S.~Muhly and B.~Solel, Tensor algebras over
$C^{*}$-correspondences:  representations, dilations and
$C^{*}$-envelopes, {\em J. Functional Analysis} \textbf{158} (1998),
no. 2, 389--457.

\bibitem{MS-CanJ} P.S.~Muhly and B.~Solel, Tensor algebras, induced
representations, and the Wold decomposition, {\em Canadian J.~Math.}
\textbf{51} (4) (1999), 850--880.

\bibitem{MS-Annalen} P.S.~Muhly and B.~Solel, Hardy algebras,
$W^{*}$-correspondences and interpolation theory, {\em Math.~Annalen}
\textbf{330} (2004), 353--415.

\bibitem{MS-models} P.S.~Muhly and B.~Solel, Canonical models for
representations of Hardy algebras, {\em Integral Equations and
Operator Theory} \textbf{53} (2995), 411-452.

\bibitem{MS-OT} P.S.~Muhly and B.~Solel, Hardy algebras associated
with $W^{*}$-correspondences (point evaluation and Schur class
functions), in {\em Operator Theory, systems theory and scattering
theory: multidimensional generalizations} (Ed. D.~Alpay and
V.~Vinnikov), pp. 221--241, \textbf{OT 157} Birkh\"auser-Verlag,
Basel-Boston, 2005.


\bibitem{MS-Schur} P.S.~Muhly and B.~Solel, {\em Schur class operator
functions and automorphisms of Hardy algebras}, preprint.




\bibitem{Popescu-CLT1} G.~Popescu, Isometric dilations for infinite
sequences of noncommuting operators, {\it Trans.~Amer.~Math.~Soc.}
316 (1989), 523--536.


%\bibitem{Popescu-multi} G.~Popescu, Multi-analytic operators on
%Fock spaces, {\em Math.~Ann.} \textbf{303} (1995), 31--46.

\bibitem{PopescuNCDA} G.~Popescu, Non-commutative disc algebras
and their representations, {\em Proceedings Amer.~Math.~Soc.}
\textbf{124} (1996), 2137-2148


%\bibitem{Popescu-Poisson}  G.~Popescu, Poisson transforms on some
%$C^{*}$-algebras generated by isometries,  {\em J.~Functional Analysis}
%\textbf{161} (1999), 27--61.

\bibitem{PopescuNF1}
G.~Popescu, Multi-analytic operators and some factorization
theorems,  {\em Indiana Univ.~Math.~J.} \textbf{38} (1989), no.~3,
693-710.

\bibitem{PopescuNF2} G.~Popescu, Characteristic functions for
infinite sequences of noncommuting operators, {\em J.~Operator
Theory} \textbf{22} (1989) no. 1, 51-71.

%\bibitem{popescu1}
%G.~Popescu, Interpolation problems in several variables,
%{\em J.~Math.~Anal.~Appl.}, \textbf{227} (1998), 227--250.

\bibitem{Popescu-Nehari} G.~Popescu, Multivariable Nehari
problem and interpolation, {\em J.~Functional Analysis} \textbf{200}
(2003), no.~2, 536--581.

\bibitem{Popescu-Memoir} G.~Popescu, {\em Entropy and multivariable
interpolation},  Memoir Amer.~Math.~Soc. \textbf{184} (2006), No.868.


%\bibitem{quig}
%P.~Quiggin, For which reproducing kernel {H}ilbert
%spaces is Pick's theorem true? {\em Integral Equations Operator Theory}
%\textbf{16} (1993), no. 2, 244--266.

\bibitem{RW} I.~Raeburn and D.P.~Williams, {\em Morita Equivalence
and Continuous-Trace $C^{*}$-Algebras},  Mathematical Surveys and
Monographs Volume \textbf{60}, American Mathematical Society,
Providence, 1998.

%\bibitem{Roesser}
%R.P.~Roesser, A discrete state-space model for
%linear image processing, {\em IEEE Trans.~Automat.~Control} \textbf{AC-20}
%(1975), no. 1, 1--10.

%\bibitem{RR} M.~Rosenblum and J.~Rovnyak, \emph{Hardy classes and
%operator theory}, Oxford Mathematical Monographs, Oxford University
%Press, New York, 1985; Dover republication, New York, 1997. MR
%\textbf{97j:47002}

\bibitem{Rowan} L.H.~Rowan, {\em Polynomial identities in ring
theory}, Volume 84 of {\em Pure and Applied Mathematics}, Academic
Press (Harcourt Brace Jovanovich Publisher), New York, 1980.


%\bibitem{Rudin-ball}
%W.~Rudin, {\em Function theory in the unit ball of ${\mathbb C}^{n}$},
%Springer-Verlag, New York, 1980.

%\bibitem{sarasonbook}
%D.~Sarason, \emph{ Sub-{H}ardy {H}ilbert spaces in the unit disk},
%John Wiley and Sons Inc., New York, 1994.

\bibitem{schur} I.~Schur,  \"Uber Potenzreihen die im Innern
des Einheitskreises  beschr\"ankt sind, {\em J.~Reine Angew.~Math.}
\textbf{14}, 7 (1917),  205--232.

\bibitem{Solel} B.~Solel, Representations of product systems over
semigroups and dilations of commuting CP maps, {\em J.~Functional
Analysis} \textbf{235} (2006) no. 2, 593-618.

\bibitem{NF} B.~Sz.-Nagy and C.~Foia\c{s}, {\em Harmonic Analysis of
Operators on Hilbert Space}, North-Holland, Amsterdam-London, 1970.

\bibitem{Rudin} W.~Rudin, {\em Functional Analysis}, Second Edition,
McGraw Hill, 1993




\end{thebibliography}
\end{document}